\documentclass[a4paper,11pt,fullpage]{article}
\usepackage{amsthm}
\usepackage{amsfonts}
\usepackage{xr}
\usepackage{graphicx}
\usepackage{amsmath}
\usepackage{epsfig,subfigure}
\usepackage{color}
\usepackage{hyperref} 

\def\R{\mathbb{R}}
\def\N{\mathbb{N}}

\def\epsilon{\varepsilon}
\def\hat{\widehat}
\def\tilde{\widetilde}

\newcommand{\me}{\mathrm{e}}

\newcommand{\be}{\begin{equation}}
\newcommand{\ee}{\end{equation}}
\newcommand{\baa}{\begin{array}}
\newcommand{\eaa}{\end{array}}
\newcommand{\ba}{\begin{eqnarray}}
\newcommand{\ea}{\end{eqnarray}}

\newtheorem{theo}{\bf Theorem}[section]
\newtheorem{lem}[theo]{\bf Lemma}
\newtheorem{pro}[theo]{\bf Proposition}

\newtheorem{cor}[theo]{\bf Corollary}

\newtheorem{rem}[theo]{\bf Remark}

\numberwithin{equation}{section}

\begin{document}
\date{}
\title{\bf{Bistable pulsating fronts in slowly oscillating one-dimensional environments }\thanks{This work has received funding from Excellence Initiative of Aix-Marseille Universit\'e~-~A*MIDEX, a French ``Investissements d'Avenir'' programme, and from the French ANR RESISTE (ANR-18-CE45-0019) and ReaCh (ANR-23-CE40-0023-02) projects. The first author is partly supported by NSFC (12471197) and GuangDong Basic and Applied Basic Research Foundation (2023B1515020034). The second author acknowledges support of the Institut Henri Poincaré (UAR 839 CNRS-Sorbonne Universit\'e), and LabEx CARMIN (ANR-10-LABX-59-01). The third author is partially supported by the NSFC (12331006).}}
\author{Weiwei Ding\footnote{School of Mathematical Sciences, South China Normal University, Guangzhou 510631, China (dingweiwei@m.scnu.edu.cn).} \qquad Fran{\c c}ois Hamel\footnote{Aix Marseille Univ, CNRS, I2M, Marseille, France (francois.hamel@univ-amu.fr). Corresponding author.} \qquad Xing Liang\footnote{School of Mathematical Sciences and Wu Wen-Tsun Key Laboratory of Mathematics, University of Science and Technology of China, Hefei, Anhui, 230026, China (xliang@ustc.edu.cn).}}
\maketitle

\begin{abstract}
We consider reaction-diffusion fronts in spatially periodic bistable media with large periods. Whereas the homogenization regime associated with small periods had been well studied for bistable or Fisher-KPP reactions and, in the latter case, a formula for the limit minimal speeds of fronts in media with large periods had also been obtained thanks to the linear formulation of these minimal speeds and their monotonicity with respect to the period, the main remaining open question is concerned with fronts in bistable environments with large periods. In bistable media the unique front speeds are not linearly determined and are not monotone with respect to the spatial period in general, making the analysis of the limit of large periods more intricate. We show in this paper the existence of and an explicit formula for the limit of bistable front speeds as the spatial period goes to infinity. We also prove that the front profiles converge to a family of front profiles associated with spatially homogeneous equations. The main results are based on uniform estimates on the spatial width of the fronts, which themselves use zero number properties and intersection arguments.
\end{abstract}


\section{Introduction and main results}\label{intro}

In this paper, we study the following reaction-diffusion equation
\be\label{eqL}
u_t=(a_L(x)u_x)_x+f_L(x,u),\ \ t\in\R,\ x\in\R,
\ee
with $L>0$, where $u_t=\partial_tu=\frac{\partial u}{\partial t}$ and $(a_L(x)u_x)_x=\partial_x(a_L(x)\partial_xu)$ and $\partial_x=\frac{\partial}{\partial x}$. The diffusion and reaction coefficients $a_L$ and $f_L$ are given by
$$a_L(x)=a\Big(\frac{x}{L}\Big)\ \hbox{ and }\ \ f_L(x,u)=f\Big(\frac{x}{L},u\Big),$$
where the function $a:\R\to\R$ is positive, of class $C^{2}(\R)$, $1$-periodic, that is, $a(x+1)=a(x)$ for all $x\in\R$, and the function $f:\R\times[0,1]\to\R,\ (x,u)\mapsto f(x,u)$ is of class $C^1$, $1$-periodic in $x$, and $(x,u)\mapsto\partial_uf(x,u)$ is Lipschitz-continuous in $u\in[0,1]$ uniformly with respect to $x\in\R$. We also assume the following conditions on~$f$:
\begin{itemize}
\item[(A1)] there is a continuous function $b:\R\to (0,1)$, $x\mapsto b(x)$, such that for every $x\in\R$, the profile $f(x,\cdot)$ is of the bistable type in the following sense
\begin{equation*}
f(x,0)=f(x,1)=f(x,b(x))=0,\ f(x,\cdot)<0\hbox{ on }(0,b(x)),\ f(x,\cdot)>0\hbox{ on }(b(x),1);
\end{equation*}
\item[(A2)] $0$ and $1$ are uniformly (in $x$) stable zeroes of $f(x,\cdot)$, in the sense that there exist $\gamma_0>0$ and $\delta_0\in(0,1/2)$ such that
\begin{equation}\label{asspars}\left\{\baa{ll}
f(x,u)\le-\gamma_0 u & \hbox{for all }(x,u)\in\R\times[0,\delta_0],\vspace{5pt}\\
f(x,u)\ge\gamma_0(1-u) & \hbox{for all }(x,u)\in\R\times[1-\delta_0,1],\eaa\right.
\end{equation}
and $\partial_uf(\cdot,0)$ and $\partial_uf(\cdot,1)$ are assumed to be of class $C^1$ in $\R$.
\end{itemize}
Notice that the assumptions (A1)-(A2) imply in particular that $\max(\partial_uf(x,0),\,\partial_uf(x,1))\leq -\gamma_0$ for all $x\in\R$. Hence, $0$ and $1$ are two linearly stable $L$-periodic steady states of~\eqref{eqL}. Here, an $L$-periodic steady state $\bar{u}$ is said to be linearly stable (resp. linearly unstable) if $\lambda_1(\bar{u})<0$ (resp.~$\lambda_1(\bar{u})>0$), where $\lambda_1(\bar{u})$ is the principal eigenvalue of  the operator $\varphi\mapsto\mathcal{L}\varphi:=\partial_x(a_L(x)\partial_x\varphi)+\partial_uf_L(x,\bar{u}(x))\varphi$ in the space of $L$-periodic functions $\varphi\in C^2(\R)$.  

In the paper, for mathematical convenience, we extend $f$ in $\R\times(\R\backslash [0,1])$ as follows:
\begin{equation}\label{extf}
\left\{\baa{ll}
f(x,u)=\partial_uf(x,0)u & \hbox{for }(x,u)\in\R\times(-\infty,0),\vspace{5pt}\\
f(x,u)=\partial_uf(x,1)(u-1) & \hbox{for }(x,u)\in\R\times(1,+\infty).
\eaa\right.
\end{equation}
Thus, it is clear that $f$ is of class $C^1$ in $\R^2$, $1$-periodic in $x$, $\min_{x\in\R}f(x,u)>0$ for all $u<0$ and $\max_{x\in\R}f(x,u)<0$ for all $u>1$, and $f(x,u)$, $\partial_uf(x,u)$ are globally Lipschitz-continuous in $u$ uniformly in $x\in\R$. It is also easily seen that this extension does not affect the behavior of the pulsating front connecting $0$ and $1$ defined below. 
 
By {\it a pulsating front of~\eqref{eqL} connecting $0$ and $1$}, we mean a classical entire solution $U_L:\R\times\R\to(0,1)$ of~\eqref{eqL} for which there exist a real number $c_L$ and a function $\phi_L:\R\times \R\to (0,1)$ satisfying
\begin{equation}\label{de-pulsating}
\left\{\baa{l}
\displaystyle U_L(t,x)=\phi_L\Big(x-c_Lt,\frac{x}{L}\Big)\ \hbox{ for all }(t,x)\in\R\times\R,\vspace{5pt}\\
\phi_L(\xi,y)\hbox{ is }1\hbox{-periodic in }y,\vspace{5pt}\\
\phi_L(-\infty,y)=1,\ \phi_L(+\infty,y)=0\hbox{ uniformly in }y\in\R.\eaa\right.
\end{equation}
The constant $c_L$ is called the {\it front speed} and $\phi_L$ is the {\it front profile}. Clearly, if $c_L\neq 0$, then
\begin{equation}\label{UL-period}
U_{L}(t,x)=U_{L}\Big(t+\frac{L}{c_L},x+L\Big)\ \hbox{ for all }(t,x)\in\R\times\R,
\end{equation}
that is, the spatial profiles $x\mapsto U_L(t,x+c_Lt)$ of the solution in the frame moving with speed~$c_L$ repeat periodi\-cally in time with period $L/c_L$. In the case $c_L=0$,~\eqref{de-pulsating} simply means that~$U_L$ does not depend on $t$ (that is, it is a steady solution) and $U_L(-\infty)=1$, $U_L(+\infty)=0$. {\it Throughout this paper, when we refer to a pulsating front for~\eqref{eqL}, we mean a solution $U_L:\R\times\R\to(0,1)$ satisfying~\eqref{de-pulsating} or, equivalently, a pair $(\phi_L,c_L)$ given as in~\eqref{de-pulsating}, with $\phi_L:\R\times\R\to(0,1)$ and~$U_L$ solving~\eqref{eqL}.}

Fronts are of great importance in reaction-diffusion equations, which are the most-used equations in population dynamics models in biology and ecology, as fronts are the main mathematical notion for the description of the invasion of a state (here, the trivial state $0$) by another one (here, the constant state $1$). The notion of pulsating front in spatially periodic media was first introduced in~\cite{skt,x1}. It is a natural generalization of the classical notion of traveling front in spatially homogeneous media where the coefficients $a$ and $f$ are independent of $x$. In that case, the function $x\mapsto b(x)$ is a constant which is the only unstable zero of $f$ in $(0,1)$. For the homogeneous bistable equation 
$$u_t=au_{xx}+f(u),$$
it is well known~\cite{aw,fm} that there exist a unique speed $c\in\R$ and a unique traveling front $u(t,x)=\phi(x-ct)$ with $0<\phi<1$ in $\R$ and $\phi(-\infty)=1$, $\phi(+\infty)=0$. The speed $c$ has the sign of the integral $\int_0^1\!f$, and the front profile $\phi$ is decreasing, unique up to shifts, and globally asymptotically stable, as in Theorem~\ref{th-known}~(iii) below. The global stability of the traveling fronts reinforces their fundamental role in reaction-diffusion equations. 

Media are however rarely homogeneous, and reaction-diffusion equations with spatially periodic coefficients, such as~\eqref{eqL}, are some of the most important classes of non-homogeneous equations. Their study has attracted much attention in the mathematical literature. One of the reasons for that interest is that the question of the existence of pulsating fronts is rather subtle. Indeed, the possible presence of multiple ordered $L$-periodic steady states may prevent the existence of pulsating fronts connecting the extremal steady states $0$ and $1$, see \cite{dg,dgm,gm,gr}. On the other hand, it is known from \cite{dgm,fz} that such a possibility does not happen provided that equation~\eqref{eqL} admits a bistable structure in the sense that any $L$-periodic steady state strictly between $0$ and $1$ is linearly unstable (we recall that $0$ and $1$ are linearly stable steady states). With such a bistable structure, pulsating fronts have been proved to exist by using different approaches \cite{du,dgm,fz}. Furthermore, in an earlier paper by the first two authors~\cite{dhz1}, the bistable structure was verified under various explicit conditions on the functions $a(x)$ and $f(x,u)$ when the spatial period $L$ is large or small. In particular, if in addition to (A1)-(A2) the function $f$ is assumed to satisfy 
\begin{itemize}
\item[(A3)]
$\qquad\qquad\qquad\displaystyle\int_0^1\!\!f(x,u)\,du>0\ \hbox{ and }\ \partial_uf(x,b(x))>0\ \hbox{ for all }x\in\R,\footnote{Together with (A1) and the $C^1$ smoothness of $f$, condition~(A3) implies that $b\in C^1(\R)$.}$
\end{itemize} 
then a pulsating front with a positive speed $c_L$ exists when $L$ is large (i.e., the medium oscillates slowly). On the other hand, if the function $x\mapsto\int_0^1f(x,u)du$ is not of a constant sign, it was shown in~\cite{dhz2} that, when $L$ is large, non-stationary pulsating fronts do not exist, while multiple stationary fronts do exist (see also~\cite{dhz2,hfr,k,x2,xz} for further results on the existence of stationary fronts for bistable equations of the type~\eqref{eqL}). We point out that, apart from the aforementioned works, the existence of pulsating fronts was investigated in~\cite{m,nr,x1,x3,z}. References~\cite{m,nr,z} were also concerned with generalized transition fronts in the sense of~\cite{bh07,bh12}, for which the media can be non-periodic (see~\cite{amo,cs,dll,e1,e2,lk,n2,p,st,ss} for further propagation or blocking results in non-periodic bistable media).

For clarity, we collect some known results from \cite{dhz2,dhz1} concerning the existence and qualitative properties of pulsating front which will be frequently used in the present work.

\begin{theo}\label{th-known}{\rm{\cite{dhz2,dhz1}}}
Let {\rm (A1)-(A3)} hold. Then there exists $L_*>0$ such that, for any $L\geq L_*$, problem~\eqref{eqL} admits a pulsating front $U_L(t,x)=\phi_L(x-c_Lt,x/L)$ with speed $c_L$. Furthermore:
\begin{itemize}
\item [{\rm (i)}] $c_L>0$, $c_L$ is the unique front speed and the map $L\mapsto c_L$ is bounded in $[L_*,+\infty)$;
\item [{\rm (ii)}] the front $U_L$ is unique up to shifts in time $t$, and $\partial_tU_L(t,x)>0$ for all $(t,x)\in\R^2$;
\item [{\rm (iii)}] the front $U_L$ is globally asymptotic stable in the sense that, for any $u_0\in C(\R)\cap L^\infty(\R)$ with
\begin{equation}\label{initial-g}
\liminf_{x\to-\infty} u_0(x) >1-\delta_0 \quad\hbox{and}\quad \limsup_{x\to+\infty} u_0(x) <\delta_0,
\end{equation}
where $\delta_0>0$ is the constant given in {\rm (A2)}, we have
$$\lim_{t\to+\infty} \| u(t,\cdot;u_0)-U_L(t+\tau_0,\cdot)\|_{L^{\infty}(\R)} = \lim_{t\to+\infty} \| u(t,\cdot;u_0)-\phi_L(\cdot-c_Lt-c_L\tau_0,\cdot/L)\|_{L^{\infty}(\R)}=0,$$
where $(t,x)\mapsto u(t,x;u_0)$ is the solution of the Cauchy problem for~\eqref{eqL} with initial condition~$u_0$, and $\tau_0\in\R$ is a constant depending on $u_0$;
\item [{\rm (iv)}] the front speed $c_L$ and the front profile $\phi_L$ are continuous with respect to $L\geq L_*$ in the sense that $c_L\to c_{L_0}$ as $L\to L_0$ in $[L_*,+\infty)$ and, up to translations of $\phi_L$ in such a way that $\phi_L(0,0)=\phi_{L_0}(0,0)$, there holds $\phi_L-\phi_{L_0}\to 0$ in $H^1(\R\times (0,1))$ and $U_L\to U_{L_0}$, $\phi_L\to\phi_{L_0}$ locally uniformly in~$\R^2$. 
\end{itemize}
\end{theo}

More precisely, we refer to \cite[Theorems~1.1,~1.5,~1.9]{dhz1} for the existence, uniqueness, monotonicity and global stability of pulsating fronts as well as the sign property of $c_L$. The continuity of $(c_L,\phi_L,U_L)$ with respect to $L$ follows directly from \cite[Theorem~1.8]{dhz1} and its proof, together with the uniform boundedness of $\|\partial_tU_L\|_{L^\infty(\R^2)}+\|\partial_xU_L\|_{L^\infty(\R^2)}$ with respect to $L\in[L_*,+\infty)$ (which itself follows from standard parabolic estimates). The boundedness of $c_L$ is an easy consequence of \cite[Theorem 1.4]{dhz2} which established a uniform bound for the propagation rates of all generalized transition fronts. 

In this paper, we are chiefly interested in the asymptotic behavior of the front speeds $c_L$ and the front profiles $\phi_L$ as $L\to+\infty$. Before stating our main results, let us first recall some related results in the spatially periodic Fisher-KPP case,\footnote{This case is named after the seminal works~\cite{f,kpp} on homogeneous equations.} i.e., when, instead of~(A1)-(A2), $f(x,u)$ is periodic in $x$ with $f(\cdot,0)=f(\cdot,1)=0$ in $\R$ and
\be\label{fkpp}
0<f(x,u)\le\partial_uf(x,0)u\ \hbox{ for all $(x,u)\in\R\times(0,1)$}.
\ee
In this case, for each $L>0$, there is a minimal front speed~$c^*_L$ of pulsating fronts~\cite{bh02,w} (namely, pulsating fronts exist for all speeds $c\ge c^*_L$, and do not exist if~$c<c^*_L$), this minimal speed~$c^*_L$ is positive and non-decreasing with respect to~$L>0$~\cite{n1}, the limit of $c^*_L$ as $L\to+\infty$ was determined in~\cite{hfr,hnr},\footnote{In multidimensional periodic media, even if the minimal speed $c^*_L(\mathrm{e})$ of pulsating fronts in a given direction~$\mathrm{e}$ is not monotone in general with respect to $L$ for equations with periodic drifts, the existence of a finite limit of~$c^*_L(\mathrm{e})$ as $L\to+\infty$ is still known, see~\cite{hnr,n1}. We also refer to~\cite{ag,dg,gr,gu} for further results on the dependence of the minimal or unique speeds of multidimensional pulsating monostable or bistable fronts with respect to the direction~$\mathrm{e}$.} and the homogenization limit of~$c^*_L$ as~$L\to0$ was characterized in~\cite{ehr,kks,sk}. The proofs in the aforementioned  papers strongly rely on the fact that the minimal speed is linearly determined and has a variational characterization in terms of principal eigenvalues of parameterized elliptic operators \cite{bhn,bhr,lz1,w}. However, for equation~\eqref{eqL} with the bistable assumptions~(A1)-(A2) or~(A1)-(A3), the unique front speed~$c_L$ is not linearly determined and it is not monotone in $L$ in general. More precisely, first of all, the limit of $c_L$ as $L\to0$ was determined in~\cite{dhz1,hps} under certain assumptions (see also Remark~\ref{remc0} below). Moreover, in the case where the diffusion rate~$a$ is a constant and~$b$ in~(A1) is piecewise constant, numerical results~\cite{hfr} showed that~$c_L$ can in some cases be increasing in small~$L$ and then decreasing for larger~$L$. On the other hand, if $a(x)$ truly depends on $x$ while the reaction is homogeneous and has the special form~$f(u)=u(1-u)(u-b)$  with $0<b<1/2$, it was proved in~\cite{ds} that the speed~$c_L$ for~$L\simeq0$ is smaller than its limit as $L\to0$.

In short, the question whether, in the bistable case~(A1)-(A3), the speeds $c_L$ of the pulsating fronts of~\eqref{eqL} (which exist for all large $L$ by Theorem~\ref{th-known}) converge as $L\to+\infty$ is the main remaining open problem in the theory of bistable pulsating fronts. This question, which is definitely more challenging than in the Fisher-KPP case~\eqref{fkpp}, is solved in this paper (see the main Theorem~\ref{speed-limit} below). Furthermore, a formula for the limit speed is established, as well as further properties on the front profiles as $L\to+\infty$.


\subsection*{Convergence of the front speeds $c_L$ as $L\to+\infty$}

In this section, we state our main result on the convergence of $c_L$ as $L\to+\infty$. To do so, we need some notations. For each $y\in\R$, let $u(t,x)=\psi(x-c(y)t,y)$ be the traveling front connecting~$0$ and~$1$, with speed~$c(y)$, of the following equation with coefficients frozen at $y$:
\begin{equation}\label{eq-homo}
u_t=a(y)u_{xx}+f(y,u)\,\,\hbox{ for $t\in\R$ and $x\in\R$}.
\end{equation}
Namely, $(\psi(\xi,y),c(y))$ satisfies
\begin{equation}\label{traveling-wave}
\left\{\baa{l}
a(y)\partial_{\xi\xi}\psi(\xi,y)+c(y)\partial_{\xi}\psi(\xi,y)+f(y,\psi(\xi,y))=0 \,\hbox{ for }\xi\in\R,\vspace{5pt}\\
\psi(-\infty,y)=1,\ \psi(+\infty,y)=0.\eaa\right.
\end{equation}
In~\eqref{eq-homo}-\eqref{traveling-wave}, $y$ is viewed as a parameter, since the derivatives only concern the variables~$(t,x)$ or~$\xi$. Under conditions~(A1)-(A3), it is well known~\cite{aw,fm} that the function $c:y\mapsto c(y)$ is 1-periodic and positive, and that for each $y\in\R$, the wave profile $\psi(\xi,y)$ is unique up to shifts in $\xi$ (see Section~\ref{sc-pre} below for more properties of the speed $c(y)$ and the wave profile $\psi(\xi,y)$, in particular the $C^1$ smoothness of $y\mapsto c(y)$). Throughout this paper, we normalize $\psi(\xi,y)$ uniquely in such a way that
\begin{equation}\label{normal-homo-wave}
\psi(0,y)=\frac{1}{2} \,\,\hbox{ for each }y\in\R,
\end{equation}
and we denote by $c_*$ the harmonic mean of the $1$-periodic function $y\mapsto c(y)$, that is,
\begin{equation}\label{limit-speed}
c_*= \left(\int_0^1  c^{-1}(y) dy \right)^{-1}.
\end{equation}

\begin{theo}\label{speed-limit}
Let {\rm (A1)-(A3)} hold, and for any $L\geq L_*$ $($with $L_*>0$ as in Theorem~$\ref{th-known}$$)$, let $c_L$ be the unique front speed of the pulsating front of~\eqref{eqL}. Then, $c_L \to c_*$ as $L\to+\infty$.
\end{theo}

Several comments are in order. Let us first point out that Theorem~\ref{speed-limit} coincides with the known result for the minimal front speeds $c^*_L$ of equations with Fisher-KPP reactions~\eqref{fkpp} in the special case where the function $r:=\partial_uf(x,0)$ is a constant, see \cite[Corollary 2.4]{hnr}. Indeed, in the latter case, the limit of the minimal front speeds~$c^*_L$ as $L\to+\infty$ is given by
$$2 \sqrt{r}\left ( \int_0^1 (a(y))^{-1/2} dy \right)^{-1},$$
which is equal to the harmonic mean of the $1$-periodic function $y\mapsto w^*(y):=2\sqrt{r\,a(y)}$, where, for each $y\in\R$, $2\sqrt{r\,a(y)}$ is the minimal front speed of the associated homogeneous equation~\eqref{eq-homo}. In the case where $\partial_uf(x,0)$ truly depends on $x$, the limit of the minimal front speeds $c^*_L$ has a variational formula in terms of the underlying parameters (see \cite[Theorem 2.3]{hnr}), but it may be not equal to the harmonic mean of the function $w^*$ (see \cite{hfr} for the case where $\partial_uf(x,0)$ is a piecewise constant function). It is noteworthy that, in our bistable case~(A1)-(A3), Theorem~\ref{speed-limit} shows that the unique wave speed $c_L$ converges to the harmonic mean $c_*$ of the function $y\mapsto c(y)$, no matter whether $f$ depends on $x$ or not.    

Theorem~\ref{speed-limit} is proved by taking the large-time and large-space scaling, i.e., $t\to Lt$, $x\to Lx$, and   determining the singular limit of the Cauchy problem of the rescaled equation
\begin{equation}\label{eq-singular}  
\partial_t u_L= L^{-1} \partial_x(a(x)\partial_x u_L)+Lf(x,u_L),\,\,  t\in\R,\,x\in\R.
\end{equation}
For a given initial condition independent of $L$, the singular limit of the solutions $u_L$ of~\eqref{eq-singular} as~$L \to +\infty$ has been widely investigated, even in the general case where the coefficients are allowed to be non-periodic. Due to the smallness of the diffusion and the largeness of the reaction in~\eqref{eq-singular}, an interface separating the regions where $u_L(t,x)\to 1$ and $u_L(t,x)\to 0$ is developed in a short time. This property has been proved by different approaches, see e.g.~\cite{g} for a probabilistic approach and~\cite{bbs,bss,bs,ds1} for a  viscosity solution approach. Particularly, in our spatially periodic case, G\"artner's result~\cite{g} suggested that, as $L\to+\infty$, the interface propagates at a mean speed equal to $c_*$. Thanks to this observation and further crucial uniform estimates on the boundedness of the wave length (more precisely, see Theorem~\ref{th-wave-length} below), we show that $c_*$ is the limit of $c_L$ by constructing suitable sub- and super-solutions of~\eqref{eq-singular}.  

Notice that, if the assumption $\int_{0}^{1}f(x,u)du>0$ in (A3) if replaced by $\int_{0}^{1}f(x,u)du<0$, equation~\eqref{eqL} admits a pulsating front with negative speed $c_L<0$ for all large $L$, see~\cite{dhz1}. In this situation, the quantity $c_*$ is negative, and similarly to Theorem~\ref{speed-limit}, we can prove that~$c_L\to c_*$ as $L\to+\infty$.  Furthermore, as an immediate consequence of this observation and \cite[Theorem~1.7]{dhz2}, we have the following classification of the asymptotic behavior of $c_L$ as $L\to+\infty$.

\begin{cor}\label{class-sign}
Let {\rm (A1)-(A2)} hold. Then, we have
\begin{itemize}
\item [{\rm (i)}] if {\rm (A3)} holds, then $c_L>0$ for all large $L$, and $c_L\to c_*$ as $L\to+\infty$;
\item [{\rm (ii)}] if $\int_0^1 f(x,u) du<0$ and $\partial_uf(x,b(x))>0$ for all $x\in\R$, then $c_L<0$ for all large $L$, and $c_L\to c_*$ as $L\to+\infty$;
\item [{\rm (iii)}] if the functions $a$, $\partial_uf(\cdot,0)$ and  $\partial_uf(\cdot,1)$ are constants and if
$$\min_{x\in\R}\int_0^1 f(x,u) du<0<\max_{x\in\R} \int_0^1 f(x,u) du,$$
then~\eqref{eqL} admits pulsating fronts for all large $L$, with necessarily $c_L=0$, that is, the fronts are stationary.
\end{itemize}
\end{cor}

\begin{rem}\label{reverse}
By {\rm{\cite{dhz1}}}, the conditions {\rm{(A1)-(A3)}} not only imply the existence of pulsating fronts of~\eqref{eqL} satisfying~\eqref{de-pulsating}, but they also provide the existence of pulsating fronts propagating in the opposite direction for the same range of periods $L$. Namely, for any $L\geq L_*$,~\eqref{eqL} admits a pulsating front of the type $\tilde{U}_L(t,x)=\tilde{\phi}_L(x+\tilde{c}_Lt,x/L)$, where $\tilde{\phi}_L$ is $1$-periodic in the second variable, and satisfies the reversed limiting conditions: 
\begin{equation}\label{limit-reverse}
\tilde{\phi}_L(-\infty,y)=0,\,\,\, \tilde{\phi}_L(+\infty,y)=1 \,\,\hbox{ uniformly in } y\in\R.
\end{equation}
Moreover, both $c_L$ and $\tilde{c}_L$ are positive. It is worthy pointing out that if the functions~$a(x)$ and~$f(x,u)$ are not even in~$x$, the speeds~$c_L$ and~$\tilde{c}_L$ may be dramatically different~{\rm{\cite{dg}}}. Nevertheless, since the function $y\mapsto c(y)$ given in~\eqref{traveling-wave} is $1$-periodic and $c_*$ in~\eqref{limit-speed} can be expressed~as
$$c_*= \left(\int_{-1/2}^{1/2} c^{-1}(y) dy \right)^{-1}= \left(\int_{-1/2}^{1/2} c^{-1}(-y) dy \right)^{-1},$$
it follows from Theorem~$\ref{speed-limit}$ that $c_L$ and $\tilde{c}_L$ have the same limit $c_*$ as $L\to+\infty$.
\end{rem}

\begin{rem}\label{remc0}
For Fisher-KPP functions $f$ satisfying~\eqref{fkpp}, for each $L>0$, pulsating fronts of~\eqref{eqL} exist, with speeds not smaller than a minimal speed $c^*_L>0$, and the map $L\mapsto c^*_L$ is non-decreasing~{\rm\cite{n1}}, hence $\lim_{L\to0}c^*_L\le\lim_{L\to+\infty}c^*_L$ $($and the inequality is strict if the coefficients $a(x)$ and $\partial_uf(x,0)$ are not constant in $x$$)$.

In the bistable case with assumptions~{\rm{(A1)-(A3)}}, the situation is in general different and more complicated. Consider for instance a function~$f$ defined by $f(x,u)=u(1-u)(u-b(x))$ for $(x,u)\in\R\times[0,1]$, where $b:\R\to(0,1/2)$ is a given $1$-periodic function of class $C^1$. Such a function~$f$ satisfies~{\rm{(A1)-(A3)}}.

On the one hand, if $a$ is constant, then it follows from~{\rm{\cite[Theorems~1.2 and~1.4]{dhz1}}} that equation~\eqref{eqL} admits a pulsating front with speed $c_L$ for every small $L>0$, and that $c_L\to c_0$ as $L\to0$, where $c_0=\sqrt{2a}\,(1/2-\bar{b})$ is the unique speed of a traveling front connecting $0$ and $1$ for the homogeneous equation $u_t= au_{xx} + u(1-u)(u-\bar{b})$, with $\bar{b}=\int_0^1b(x)\,dx\in(0,1/2)$. With the notations~\eqref{traveling-wave}, one has $c(y)=\sqrt{2a}\,(1/2-b(y))$ for each $y\in\R$, and Theorem~$\ref{speed-limit}$ implies that $c_L\to c_*=\sqrt{2a}\,\big(\int_0^1(1/2-b(y))^{-1}dy\big)^{-1}$ as $L\to+\infty$. Therefore, the Cauchy-Schwarz inequality yields $c_0\ge c_*$ and even $c_0>c_*$ as soon as $b$ is not constant. This is in sharp contrast with the inequality $\lim_{L\to0}c^*_L\le\lim_{L\to+\infty}c^*_L$ under the Fisher-KPP assumption~\eqref{fkpp}.

On the other hand, if $b\equiv\bar{b}\in(0,1/2)$ is constant, it follows from~{\rm{\cite[Theorems~1.2 and~1.4]{dhz1}}} that equation~\eqref{eqL} admits a pulsating front with speed $c_L$ for every small $L>0$,\footnote{Notice that, in this case, pulsating fronts are even known to exist for all $L>0$ if the function~$a$ is uniformly close to a constant~{\rm\cite{fz,x1}}.} and that $c_L\to c_0$ as $L\to0$, where $c_0=\sqrt{2a_H}\,(1/2-\bar{b})$ is the unique speed of a traveling front connecting $0$ and $1$ for the homogeneous equation $u_t= a_Hu_{xx} + u(1-u)(u-\bar{b})$, with $a_H=(\int_0^1a(y)^{-1}dy)^{-1}$ being the harmonic mean of the diffusion coefficient $a$. With the notations~\eqref{traveling-wave}, one has $c(y)=\sqrt{2a(y)}\,(1/2-\bar{b})$ for each $y\in\R$, and Theorem~$\ref{speed-limit}$ implies that $c_L\to c_*=\sqrt{2}\,(1/2-\bar{b})\,\big(\int_0^1\sqrt{a(y)}^{-1}dy\big)^{-1}$ as $L\to+\infty$. Therefore, the Cauchy-Schwarz inequa\-lity yields $c_0\le c_*$ and even $c_0<c_*$ as soon as $a$ is not constant.

To sum up, the influence of the oscillations of the heterogeneities in~\eqref{eqL} under the bistable assumptions~{\rm{(A1)-(A3)}} turns out to be more complicated than under the Fisher-KPP assumption~\eqref{fkpp}. Moreover, the oscillations in the diffusion and reaction terms can lead to opposite effects on the bistable front speed. 
\end{rem}


\subsection*{Convergence of the front profiles $\phi_L$ as $L\to+\infty$}
 
The second main result of the present paper is concerned with the convergence as $L\to+\infty$ of the front profiles $\phi_L(\xi,y)$ given in~\eqref{de-pulsating}. It is difficult to address this problem by considering the singular equation~\eqref{eq-singular}, since the solutions $\phi_L(L(x-c_Lt),x)$ of~\eqref{eq-singular} converge to a step function as $L\to+\infty$ for any $t\in\R$ (this will actually be a consequence of Theorem~\ref{wave-length-1} below), hence the information of $\phi_L(\xi,y)$ is lost when passing to the limit as $L\to+\infty$. We instead show the convergence of $\phi_L(\xi,y)$ by constructing sub- and super-solutions for the original equation~\eqref{eqL} when $L$ is large (as a matter of fact, this approach also provides the convergence of $c_L$ as $L\to +\infty$, as will follow from the proof of Theorem~\ref{front-limit} below).

Before stating the result, we point out that, since equation~\eqref{eqL} is invariant by translation in time $t$, we have to suitably shift the fronts in order to pass to the limit as $L\to+\infty$. Namely, with $L_*>0$ as in Theorem~\ref{th-known}, since for each $L\ge L_*$ the pulsating front $U_L$ is unique up shifts in time, with $U_L(-\infty,x)=0$ and $U_L(+\infty,x)=1$ for all $x\in\R$, we normalize for definiteness $U_L$ by requiring 
\begin{equation}\label{normal-U}
U_L(0,0)=\frac{1}{2}.
\end{equation}
Clearly, under such a normalization, $U_L(t,x)$ is uniquely determined, and the front profile $\phi_L(\xi,y)$ is uniquely determined and satisfies $\phi_L(0,0)=1/2$. Now, for each $L\ge L_*$, since $\partial_{\xi}\phi_L<0$ in~$\R^2$ and $\phi_L$ satisfies~\eqref{de-pulsating} and is at least of class $C^1(\R^2)$, the implicit function theorem yields the existence of a uniquely determined~$C^1(\R)$ function $y\mapsto\zeta_L(y)$ such that
\begin{equation}\label{normmal-front-1}
\phi_L(\zeta_L(y),y)=\frac{1}{2}.
\end{equation}
Clearly, $\zeta_L(0)=0$ and  $\zeta_L(y)$ is $1$-periodic in $y\in\R$. Moreover, it follows from Theorem~\ref{th-known}~(ii) and~(iv) that the function $(L,y)\mapsto\zeta_L(y)$ is continuous in $[L_*,+\infty)\times\R$.

For the convergence of the front profiles as $L\to+\infty$, in addition to (A1)-(A3), the following condition will also be assumed:
\begin{itemize}
\item [(A4)]
the function $a$ is a constant, and there exists a constant $\delta_0' \in (0,1/2)$ such that $\partial_xf(x,u)=0$ for all $x\in\R$ and $u\in [0,\delta_0']\cup[1-\delta_0',1]$, and then for all $x\in\R$ and $u\in [-\delta_0',\delta_0']\cup[1-\delta_0',1+\delta_0']$ by~\eqref{extf}.
\end{itemize}

\begin{theo}\label{front-limit}
Let {\rm (A1)-(A4)} hold. For each $L\geq L_*$, let $U_L(t,x)=\phi_L(x-c_Lt,x/L)$ be the pulsating front of~\eqref{eqL} given by Theorem~$\ref{th-known}$, and let $\zeta_L$ be the function given by~\eqref{normmal-front-1}. Let also~$(\xi,y)\mapsto\psi(\xi,y)$ and $y\mapsto c(y)$ be given in~\eqref{traveling-wave}-\eqref{normal-homo-wave}. Then,
\begin{equation}\label{eq-front-limit}
\sup_{y\in\R,\,\xi\in\R}\left|\phi_L \left(\xi+ \zeta_L(y),\frac{\xi}{L}+y \right) -\psi(\xi,y) \right|\mathop{\longrightarrow}_{L\to+\infty}0.
\end{equation}
Moreover,
\begin{equation}\label{asym-zeta}
\sup_{y\in\R,\,x\in[-A,A]}\left| \zeta_L\left(y+\frac{x}{L}\right) -  \zeta_L(y) -L\int_{y}^{y+x/L} \left( 1-\frac{c_L}{c(s)}\right)ds \right|\mathop{\longrightarrow}_{L\to+\infty}0
\end{equation}
for every $A>0$, and
\begin{equation}\label{zetaL-zetaLx}
\sup_{L\geq L_*,\, y\in\R,\,x\in[-L,L]}\left| \zeta_L\left(y+\frac{x}{L}\right) -  \zeta_L(y) -L\int_{y}^{y+x/L} \left( 1-\frac{c_L}{c(s)}\right)ds \right|<+\infty.
\end{equation}
\end{theo}

Several remarks are in order. First of all, by choosing $y=0$ and $x=y_0 L$ for some $y_0\in (0,1]$ in~\eqref{zetaL-zetaLx}, one has
$$\sup_{L\geq L_*}\left| \zeta_L(y_0) -L\int_{0}^{y_0} \left( 1-\frac{c_L}{c(s)}\right)ds \right| <+\infty.$$
In the special case $y_0=1$, since $\zeta_L(1)=\zeta_L(0)=0$, we rediscover that $c_L\to c_*=(\int_0^1c(s)^{-1}ds)^{-1}$ as $L\to+\infty$, here under the additional assumption~(A4). Furthermore, since the function $y\mapsto c(y)$ is in general not constant even under the assumption (A4), it then follows that $\zeta_L(y_0)$ is unbounded in $L\geq L_*$ for each $y_0\in(0,1)$ such that $y_0\int_0^1c(s)^{-1}ds\neq\int_0^{y_0}c(s)^{-1}ds$. On the other hand,~\eqref{zetaL-zetaLx} implies that $\zeta_L(y_0)=O(L)$ as $L\to+\infty$ for all $y_0\in\R$ (that will coincide with the observation in Theorem~\ref{wave-length-1}~(iii) below).

Secondly, we point out that~\eqref{eq-front-limit} can be equivalently rephrased as
$$\sup_{y\in\R,\,\xi\in\R}\left|U_L\left(\frac{Ly-\zeta_L(y)}{c_L},\xi+Ly\right) -\psi(\xi,y) \right|\mathop{\longrightarrow}_{L\to+\infty}0.$$
As we will see, this last formulation implies that, as $L\to+\infty$, the spatial profiles $U_L(t,\cdot)$ look, at any time~$t$ and uniformly in space, like shifts of the profiles $\psi(\cdot,y)$ of traveling fronts of the homogeneous equations~\eqref{traveling-wave} for some values of~$y$. Indeed, by~\eqref{zetaL-zetaLx} and $\lim_{L\to+\infty}c_L\!=\!c_*$, the continuous functions $z\mapsto z-\zeta_L(z)/L$ converge as $L\to+\infty$ locally uniformly in $z\in\R$ to the function $z\mapsto c_*\int_0^zc(s)^{-1}ds$ (by applying~\eqref{zetaL-zetaLx} with $y=0$ and $x=Lz$). Therefore, for every~$L>0$ large enough and for every $t\in[0,L/c_L]$, there is $z_{L,t}\in\R$ such that $z_{L,t}-\zeta_L(z_{L,t})/L=c_Lt/L$. Hence $U_L(t,\cdot)=U_L((Lz_{L,t}-\zeta_L(z_{L,t}))/c_L,\cdot)$ and $\|U_L(t,\cdot)-\psi(\cdot-Lz_{L,t},z_{L,t})\|_{L^\infty(\R)}\to0$ as $L\to+\infty$ uniformly in $t\in[0,L/c_L]$. Since the profiles $x\mapsto U_L(t,x+c_Lt)$ in the frames moving with speeds $c_L$ are periodic in time with period $L/c_L$, one finally gets that
\be\label{eq-front-limitbis}
\sup_{t\in\R}\ d\Big(U_L(t,\cdot),\big\{\psi(\cdot+a,y):a\in\R,\,y\in\R\big\}\Big)\mathop{\longrightarrow}_{L\to+\infty}0,
\ee
where $d$ is the distance associated with the norm in $L^\infty(\R)$. Actually, the convergence of the profiles of the pulsating fronts to profiles of homogeneous traveling fronts as $L\to+\infty$ is quite natural since equation~\eqref{eqL} looks like, as $L\to+\infty$, families of locally spatially homogeneous equations. But the difficulty in Theorem~\ref{front-limit} and in the consequent formula~\eqref{eq-front-limitbis} is to show rigorously this convergence and its uniformity with respect to time and space.

The strategy of the proof of Theorem~\ref{front-limit} can be briefly summarized as follows. First, by virtue of Theorem~\ref{th-known}~(iii) on the global stability of the pulsating front $\phi_L$, we show that, up to some shift in $\xi$, $\phi_L(\cdot,y)$ is close to $\psi(\cdot,y)$ when $L$ is sufficiently large.  Next, we use Theorem~\ref{wave-length-1} below on the uniform boundedness of the front width to determine the limit shifts appearing in~\eqref{eq-front-limit}. The assumption (A4) is used in the construction of suitable sub- and super-solutions (see Lemmas~\ref{super-solu}-\ref{sub-solu} below), which is the key point of the proof of Theorem~\ref{front-limit}. We believe that (A4) is only a technical assumption, as one can see from Theorem \ref{speed-limit} that the wave speed $c_L$ converges to $c_*$ without this assumption. However, at the moment, it is unclear how to remove this assumption in Theorem~\ref{front-limit}.


\subsection*{Uniform boundedness of front width}

As a matter of fact, Theorems~\ref{speed-limit} and~\ref{front-limit}  on the limit speeds and profiles of the pulsating fronts as $L\to+\infty$ are based on the uniform boundedness of the width of the front profiles, in a sense to be made more precise below. Under conditions~(A1)-(A3), with $L_*>0$ as in Theorem~\ref{th-known}, it follows from the definition~\eqref{de-pulsating} of pulsating fronts that, for any $\delta\in(0,1/2]$ and $L\ge L_*$, there exists a constant $C>0$ such that 
\begin{equation}\label{eq-wave-lengh-t}
\sup_{x\in\R}\ {\rm{diam}}\,\big(\left\{t\in\R: \, \delta\leq U_L(t,x) \leq 1-\delta \right\}\big)\le C
\end{equation} 
and 
\begin{equation}\label{eq-wave-lengh-x}
\sup_{t\in\R}\ {\rm{diam}}\,\big(\left\{x\in\R: \, \delta\leq U_L(t,x) \leq 1-\delta \right\}\big)\le C, 
\end{equation}
where ${\rm{diam}}(E)$ denotes the diameter of a set $E\subset\R$. Notice that, for each $\delta\in(0,1/2]$, $L\ge L_*$ and $x\in\R$, the set $\{t\in\R: \delta\leq U_L(t,x) \leq 1-\delta\}$ is a compact interval, from~\eqref{de-pulsating} and the continuity of $U_L$ and its monotonicity with respect to $t$. Hence the diameter of the set in~\eqref{eq-wave-lengh-t} is equal to its length. We also point out that property~\eqref{eq-wave-lengh-x} could be rephrased equivalently by saying that, for any $\lambda\in(0,1)$, the distance between $c_Lt$ and the level set $\{x\in\R:U_L(t,x)=\lambda\}$ is uniformly bounded in $t$.

Our next result, which is a key stone in determining the convergences of front profiles $\phi_L$ and front speeds $c_L$ as $L\to+\infty$, and which has its own interest, shows that the constants $C$ in~\eqref{eq-wave-lengh-t}-\eqref{eq-wave-lengh-x} can be taken independently of $L\ge L_*$.

\begin{theo}\label{th-wave-length}
Let {\rm (A1)-(A3)} hold. For each $L\geq L_*$, with $L_*$ as in Theorem~$\ref{th-known}$, let $U_L$ be a pulsating front of~\eqref{eqL}. Then, for any $\delta\in (0,1/2]$, there exist constants $C>0$ and $\beta>0$ $($both are independent of $L$$)$ such that~\eqref{eq-wave-lengh-t} and~\eqref{eq-wave-lengh-x} hold true for all $L\ge L_*$, and 
\begin{equation}\label{derivative-t} 
\partial_t U_L(t,x) \geq \beta \,\,\hbox{ for all }L\ge L_*\hbox{ and }(t,x)\in\R^2\hbox{ such that }\delta\leq U_L(t,x) \leq 1-\delta.
\end{equation}
\end{theo}

The following theorem is a consequence of Theorem~\ref{th-wave-length}. 

\begin{theo}\label{wave-length-1}
Let {\rm (A1)-(A3)} hold. For each $L\geq L_*$, with $L_*$ as in Theorem~$\ref{th-known}$, let $\phi_L$ be a front profile of~\eqref{eqL} and let $\zeta_L$ be defined by~\eqref{normmal-front-1} and $\zeta_L(0)=0$. Then,
\begin{itemize}
\item [{\rm (i)}] the following two convergences hold uniformly in $L\geq L_*$ and $y\in\R$:
\begin{equation*}
\lim_{\xi\to-\infty}\phi_L(\xi+\zeta_L(y),y)=1\quad\hbox{and}\quad \lim_{\xi\to+\infty}\phi_L(\xi+\zeta_L(y),y)=0;
\end{equation*}
\item [{\rm (ii)}] for any $\delta\in (0,1/2]$, there exists $\beta> 0$ $($independent of $L$$)$ such that
\begin{equation*}
\partial_{\xi} \phi_L(\xi,y)\leq -\beta \,\,\hbox{ for all }L\ge L_*\hbox{ and }(\xi,y)\in\R^2\hbox{ such that }\delta\leq \phi_L(\xi,y) \leq 1-\delta;
\end{equation*}
\item [{\rm (iii)}] the function $(L,y)\mapsto \zeta_L(y)/L$ is continuous and bounded in $[L_*,+\infty)\times\R$.  
\end{itemize}
\end{theo}

Notice that the continuity of $(L,y)\mapsto \zeta_L(y)/L$ in $[L_*,+\infty)\times\R$ follows from the comments after~\eqref{normmal-front-1}, while the quantities $\zeta_L(y)$ are in general unbounded with respect to $L\ge L_*$, from the discussion after Theorem~\ref{front-limit}.

Theorems~\ref{th-wave-length} and~\ref{wave-length-1} mean that the fronts do not flatten as $L\to+\infty$. Their proofs are based on some intersection number arguments inspired by~\cite{a} and from the comparison of the pulsating fronts with some well chosen stationary solutions. 


\subsection*{Outline of the paper}

Section~\ref{sec2} is devoted to the proof of Theorems~\ref{th-wave-length} and~\ref{wave-length-1} on the uniform boundedness of the wave lengths of the pulsating fronts $U_L$ with respect to $L\ge L_*$. The proof of Theorem~\ref{speed-limit} on the limit of the speeds $c_L$ as $L\to+\infty$ is done in Section~3, and that of Theorem~\ref{front-limit} on the limit of the profiles is done in Section~\ref{sec4}. Lastly, Section~\ref{sec5} is an appendix on some properties of the profiles $\psi(\xi,y)$ and speeds $c(y)$ of the family of homogeneous equations~\eqref{eq-homo}-\eqref{traveling-wave}.


\section{Uniform boundedness of front width: proofs of Theorems~\ref{th-wave-length} and~\ref{wave-length-1}}\label{sec2}

This section is devoted to the proof of Theorems~\ref{th-wave-length} and~\ref{wave-length-1}. Our proof strongly relies on a zero number argument, and the classification of all solutions to the ordinary differential equations
\begin{equation}\label{ode}
(a_Lw')'+f_L(x,w)=0  \hbox{ in }\R,\qquad  w(0) \in (0,1),
\end{equation}
according to the number of intersections of their graphs with those of the constant solutions $0$ and $1$. First, after recalling in Section~\ref{sec21} some properties on intersection number arguments, we investigate in Section~\ref{sec22}, for any given $L\geq L_*$ with $L_*>0$ as in Theorem~\ref{th-known} and for any time $t\in\R$, the number and type of intersection points of the spatial profiles $U_L(t,\cdot)$ of the pulsating fronts $U_L$ with all solutions to~\eqref{ode}. Then, using a passage to the limit as $L\to+\infty$ and some properties of the stationary solutions to the $y$-frozen equation~\eqref{eq-homo}, we prove Theorem~\ref{th-wave-length} by contradiction in Section~\ref{sec23}. Lastly, we show in Section~\ref{sec24} that Theorem~\ref{th-wave-length} implies Theorem~\ref{wave-length-1}. 


\subsection{Zero number properties}\label{sec21}

We begin with the collection of some zero number properties which will be used later. We refer the reader to \cite{a,dm,dgm,po} for a more detailed overview of the general arguments and their applications.

Let us first introduce some fundamental notions. 
For any sign-changing function $w:\R\to\R$, let $\mathcal{Z}[w]$ denote the number of sign changes of $w$, and let $SGN[w]$ (which is defined when $\mathcal{Z}[w]<+\infty$) be a word consisting of $+$ and $-$, which describes the sign changes of $w$. More precisely,  $\mathcal{Z}[w]$ is the supremum of all $k\in\N^*$ such that there exist real numbers $x_1 <x_2 < \cdots <x_{k+1}$ with
$$w(x_i) \cdot  w(x_{i+1}) < 0\,\,\hbox{ for all }  i =1, 2, \cdots,k,$$
and if $\mathcal{Z}[w]<+\infty$, $SGN[w]=[sgn(w(x_1)),\cdots, sgn(w(x_{k+1}))]$, where $sgn(w)\in\{-,+\}$ is the classical sign function, and $x_1 < \cdots <x_{k+1}$ is the sequence that appears in the definition of $\mathcal{Z}[w]$ with maximal $k$. When $w$ does not change sign, for convenience, we set $\mathcal{Z}[w]=0$ and $SGN[w]=[sgn(w)]$ if $w\not\equiv 0$, and set $\mathcal{Z}[w]=-1$ and $SGN[w]=[\,\,]$ (the empty word) if $w\equiv 0$.     
By definition, the length of the word  $SGN[w]$ is equal to $\mathcal{Z}[w]+1$.

For any two words $A,B$ consisting of $+$ and $-$, we write $A\rhd B$ (or equivalently, $B\lhd A$) if $B$ is a subword of $A$. For example, any element in $\{[+-],[+],[-],[\,\,]\}$  is a subword of $[+-]$.

It follows from the above definitions that $\mathcal{Z}$ is semi-continuous with respect to the pointwise convergence, as stated below.

\begin{lem}\label{semi-continuity}
Let $(w_n)_{n\in\N}$ be a sequence of real-valued functions converging to $w$ pointwise on~$\R$. Then,
$$\mathcal{Z}[w] \leq \liminf_{n\to+\infty} \mathcal{Z}[w_n]  ,\quad\hbox{and}\quad SGN[w] \lhd \liminf_{n\to+\infty} SGN[w_n].$$
\end{lem}

The following lemma is an easy application of the zero number properties for solutions to a linear parabolic equation of the form
\begin{equation}\label{eq-linear}
\partial_t w=a(t,x)\partial_{xx}w+b(t,x)\partial_xw+ c(t,x)w \,\,\hbox{ in }  (t_1,t_2)\times\R,
\end{equation}
where $-\infty\le t_1<t_2\le+\infty$ and the coefficients $a>0,\,a^{-1},a_t,\,a_x,\,a_{xx},\,b,\,b_t,\,b_x,\,c$ belong to~$L^{\infty}((t_1,t_2)\times\R)$.

\begin{lem}\label{zero-decrease} {\rm{\cite{a,dm,dgm}}}
Let $L>0$, let $u_1(t,x)$ be a bounded solution of the Cauchy problem associated with~\eqref{eqL} in $(0,+\infty)\times\R$, with a piecewise continuous bounded initial condition $u_1(0,\cdot)$, and let $u_2(x)$ be a $C^2(\R)$ stationary solution of~\eqref{eqL}. Assume also that $u_1(0,\cdot)-u_2$ changes sign at most finitely many times on $\R$, that is, $\mathcal{Z}[u_1(0,\cdot)-u_2]<+\infty$. Then,
\begin{itemize}
\item [{\rm (i)}] For any $0\leq t<t'<+\infty$, we have
$$\mathcal{Z}[u_1(t,\cdot)-u_2] \geq  \mathcal{Z}[u_1(t',\cdot)-u_2]$$
and
$$SGN[u_1(t,\cdot)-u_2] \rhd  SGN[u_1(t',\cdot)-u_2].$$
\item [{\rm (ii)}] If $u_1(t_*,x_*)=u_2(x_*)$ and $\partial_x u_1(t_*,x_*)=u_2'(x_*)$ for some $t_*>0$ and $x_*\in\R$, and if~$u_1(t,x)\not \equiv u_2(x)$ in $(0,+\infty)\times\R$, then
$$\mathcal{Z}[u_1(t,\cdot)-u_2]-2\geq \mathcal{Z}[u_1(s,\cdot)-u_2] \geq 0 \,\,\hbox{ for all } t\in (0,t_*)\hbox{ and }s\in (t_*,+\infty).$$
\end{itemize}
\end{lem}

Since the function $(x,u)\mapsto f(x,u)$ is globally Lipschitz-continuous with respect to $u\in\R$ uniformly in $x\in\R$, it is easily seen that $w(t,x):=u_1(t,x)-u_2(x)$ satisfies an equation of the form~\eqref{eq-linear} in $(0,+\infty)\times\R$ with a bounded continuous coefficient $c$ defined by $c(t,x):=(f_L(x,u_1(t,x))-f_L(x,u_2(x)))/(u_1(t,x)-u_2(x))$ if $u_1(t,x)\neq u_2(x)$ and $c(t,x):=\partial_uf_L(x,u_1(t,x))$ if $u_1(t,x)=u_2(x)$. Then, the proof of Lemma~\ref{zero-decrease} follows from that of \cite[Lemma 2.4]{dgm} with some obvious modifications; therefore, we omit the details.


\subsection{Intersections of pulsating fronts with stationary solutions}\label{sec22}

Now, for any given $L\geq L_*$, with $L_*>0$ as in Theorem~\ref{th-known}, we investigate $\mathcal{Z}[U_L(t,\cdot)-w]$ and~$SGN[U_L(t,\cdot)-w]$ at any time $t\in\R$, where $U_L:\R^2\to(0,1)$ is the pulsating front of~\eqref{eqL} and~$w$ is an arbitrary $C^2(\R)$ solution to~\eqref{ode}. Let us first classify the solutions of~\eqref{ode} as follows.

\begin{lem}\label{classfy-w}
Let $L\geq L_*$ be fixed, and let $w\in C^2(\R)$ be a solution of~\eqref{ode}. Then, $w$ must be one of the following types:
\begin{itemize}
\item [{\rm (a)}] $\mathcal{Z}[w-1]=1$ and $\mathcal{Z}[w]=1$;
\item [{\rm (b)}] $\mathcal{Z}[w-1]=0$, $\mathcal{Z}[w]=0$, and $0<w<1$ in $\R$;
\item [{\rm (c)}] $\mathcal{Z}[w-1]=0$,  $\mathcal{Z}[w]=1$, $w<1$ in $\R$, and $SGN[w]=[-,+]$;
\item [{\rm (d)}] $\mathcal{Z}[w-1]=0$,  $\mathcal{Z}[w]=1$, $w<1$ in $\R$, and $SGN[w]=[+,-]$;
\item [{\rm (e)}] $\mathcal{Z}[w-1]=0$, $\mathcal{Z}[w]=2$, $w<1$ in $\R$, and $SGN[w]=[-,+,-]$.
\end{itemize}
\end{lem}

\begin{proof}
First of all, since $w(0)\in(0,1)$, we have $\mathcal{Z}[w-1]\ge0$ and $\mathcal{Z}[w]\ge0$. Suppose now that there exists some $x_1\in\R$ such that $w(x_1)=1$.  It is then clear that $w'(x_1)\neq 0$ (otherwise, since~$f_L(\cdot,1)\equiv 0$, the Cauchy-Lipschitz theorem would imply that $w\equiv 1$ in $\R$, a contradiction with~$w(0) \in (0,1)$).

Then, we observe that, if $w'(x_1)>0$, then $w(x)>1$ for all $x>x_1$. Indeed, otherwise, there would exist some $x_2>x_1$ such that $w(x_2)=1$, $w'(x_2)<0$ and $w(x)>1$ for all $x\in (x_1,x_2)$. Since  $\max_{\R}f_L(\cdot,u)<0$ for all $u>1$, one would have $(a_Lw')'>0$ in $(x_1,x_2)$, and hence, $a_L(x_2)w'(x_2) > a_L(x_1)w'(x_1) >0$. Since $a$ is a positive function, there holds $w'(x_2)>0$, which is a contradiction. Similarly, we can conclude that, if $w'(x_1)<0$, then $w(x)>1$ for all $x<x_1$. As a consequence, we obtain that $\mathcal{Z}[w-1]\leq 2$, and if $\mathcal{Z}[w-1]=2$, then $SGN[w-1]=[+,-,+]$.

Since $\min_{\R}f_L(\cdot,u)>0$ for all $u<0$, with similar arguments as above, one gets that $\mathcal{Z}[w]\leq 2$, and if $\mathcal{Z}[w]=2$, then $SGN[w]=[-,+,-]$.

Consequently, we always have
$$\mathcal{Z}[w-1]+ \mathcal{Z}[w] \leq 2.$$
Observe that the Cauchy-Lipschitz theorem also implies that, if $w\le1$ in $\R$ (respectively $w\ge0$ in $\R$), then $w<1$ in $\R$ (respectively $w>0$ in $\R$). Therefore, since $0<w(0)<1$, one gets that~$w<1$ in~$\R$ if $\mathcal{Z}[w-1]=0$ (respectively $w>0$ in $\R$ if $\mathcal{Z}[w]=0$). 

Now, to complete the proof, it remains to exclude the following three cases:
\begin{itemize}
\item $\mathcal{Z}[w-1]=1$, $\mathcal{Z}[w]=0$ and $SGN[w-1]=[+,-]$;
\item $\mathcal{Z}[w-1]=2$, $\mathcal{Z}[w]=0$ and $SGN[w-1]=[+,-,+]$;
\item $\mathcal{Z}[w-1]=1$, $\mathcal{Z}[w]=0$ and $SGN[w-1]=[-,+]$.
\end{itemize}
In the first two cases, one can find a continuous function $g:\R\to[0,1]$ satisfying~\eqref{initial-g} such that~$g< w$ in $\R$. Let $u(t,x)$ be the solution of the Cauchy problem associated with~\eqref{eqL} with initial condition $u(0,x)=g(x)$. Then, the comparison principle immediately implies that
\begin{equation}\label{comp-w-u}
u(t,x) < w(x)\,\,\hbox{ for all }t\geq 0\hbox{ and }x\in\R.
\end{equation}
On the other hand, it follows from Theorem~\ref{th-known} (iii) that $\|u(t,\cdot)-U_L(t+\tau_0,\cdot)\|_{L^{\infty}(\R)}\to0$ as~$t\to+\infty$, for some $\tau_0\in\R$. Remember that the wave speed $c_L$ is positive. This implies in particular that $u(t,0)\to 1$ as~$t\to+\infty$. Since $w(0)\in (0,1)$, we have $u(t,0)>w(0)$ for all large~$t$, which is a contradiction with~\eqref{comp-w-u}. 

Finally, if the last case occurs, then one reaches a similar contradiction by noticing that equation~\eqref{eqL} admits a pulsating front of the type $\tilde{U}_L(t,x)=\tilde{\phi}_L(x+\tilde{c}_Lt,x/L)$ satisfying the asymptotic conditions~\eqref{limit-reverse}, and that such a front is also globally asymptotic stable, with positive speed $\tilde{c}_L$ (see Remark~\ref{reverse}). The proof of Lemma~\ref{classfy-w} is thus complete.
\end{proof}

\begin{lem}\label{intersection}
Let $L\geq L_*$ be fixed, let $U_L$ be a pulsating front of~\eqref{eqL}, and let $w\in C^2(\R)$ be a solution of $(a_Lw')'+f_L(x,w)=0$ in $\R$. Then $\mathcal{Z}[U_L(t,\cdot)-w]\leq 2$ for all $t\in\R$. Furthermore, if~$\mathcal{Z}[U_L(t,\cdot)-w]= 2$ for some $t\in\R$, then there must hold
\begin{equation}\label{sgn-UL-w}
SGN[U_L(t,\cdot)-w]=[+,-,+].
\end{equation}
\end{lem}

\begin{proof}
Let $t_*\in\R$ be an arbitrary time such that $U_L(t_*,x_*)=w(x_*)$ for some $x_*\in\R$. Denote $\alpha:=w(x_*)$. Notice that $\alpha\in(0,1)$ since $0<U_L<1$ in $\R^2$, hence $w(\cdot+x_*)$ solves~\eqref{ode} and Lemma~\ref{classfy-w} applies to $w(\cdot+x_*)$ and then to $w$ as well. To show Lemma~\ref{intersection}, it suffices to prove that either $\mathcal{Z}[U_L(t_*,\cdot)-w|\le1$, or $\mathcal{Z}[U_L(t_*,\cdot)-w]= 2$ and~\eqref{sgn-UL-w} holds with $t=t_*$.

For any $z\in\R$, let $\hat{u}(t,x;z)$ be the solution of the Cauchy problem associated with~\eqref{eqL} with initial condition
$$\hat{u}(0,x;z)=H(z-x),$$
where $H$ denotes the Heaviside function defined by $H(y)=0$ if $y<0$ and $H(y)=1$ if $y\geq 0$. When $z<x_*$, we define
\begin{equation}\label{de-tau} 
\tau(z):=\min\big\{t>0: \, \hat{u}(t,x_*;z)=\alpha\big\}. 
\end{equation}
Let us first observe that $\tau(z)$ is well defined and $0<\tau(z)<+\infty$. Indeed, since $c_L>0$, it follows from Theorem~\ref{th-known} (iii) that $\lim_{t\to+\infty}\hat{u}(t,x_*;z)=1>\alpha$. This, together with the continuity of~$\hat{u}(\cdot,x_*;z)$ in $[0,+\infty)$ and the fact that $\hat{u}(0,x_*;z)=0<\alpha$, immediately implies that the minimum in~\eqref{de-tau} is well defined and $0<\tau(z)<+\infty$.  Moreover, since $\hat{u}(t,x;z)\to0$ as~$z\to -\infty$ locally uniformly in $(t,x)\in[0,+\infty)\times\R$, we have
$$\lim_{z\to-\infty}\tau(z)= +\infty.$$
Furthermore, by the proof of \cite[Lemma 3.1]{dgm}, there is an entire solution $\hat{u}_{\infty}:\R^2\to(0,1)$ of~\eqref{eqL} that is steeper than any other entire solution between $0$ and $1$,\footnote{The steepness is understood in the following sense: for any two entire solutions $u_1:\R\times\R\to(0,1)$ and $u_2:\R\times\R\to(0,1)$ of~\eqref{eqL}, we say that $u_1$ is steeper than $u_2$ provided that $SGN[u_1(t+t',\cdot)-u_2(t,\cdot)] \lhd [+,-]$ for any $t$ and $t'$ in $\R$.} and such that the following limit exists in the topology of $C^{1;2}_{t;x;loc}(\R^2)$, up to extraction of a subsequence:
\be\label{de-uinfty}
\lim_{k\to+\infty} \hat{u}(t+\tau(-kL),x;-kL)=\hat{u}_{\infty}(t,x).
\ee
Next, we claim that
\begin{equation}\label{id-hatu-UL}
\hat{u}_{\infty}(t,x)\equiv U_L(t+t_*,x).
\end{equation}
Indeed, since the pulsating front of~\eqref{eqL} is unique up to shifts in time (see Theorem~\ref{th-known} (ii)) and since $U_L$ is steeper than any other entire solution between $0$ and $1$ by \cite[Theorem 1.8]{gm}, the function~$U_L$ is then identically equal to $\hat{u}_{\infty}$ up to a shift in time. Since $\hat{u}_{\infty}(0,x_*)= U_L(t_*,x_*)=\alpha$, we obtain~\eqref{id-hatu-UL}.

For clarity, we divide the remaining arguments into two parts, according to the behavior of~$w$. First, if $w$ is one of the types (a)-(c) listed in Lemma~\ref{classfy-w}, then there exists $k_*\in\N$ sufficiently large such that
$$\mathcal{Z}[\hat{u}(0,\cdot;-kL)-w] =1 \,\,\hbox{ for all }k\geq k_*.$$
It then follows from Lemma~\ref{zero-decrease} that
$$\mathcal{Z}[\hat{u}(t+\tau(-kL),\cdot;-kL)-w]\le1 \,\,\hbox{ for all }t \geq -\tau(-kL)\hbox{ and }k\geq k_*.$$
Hence, by Lemma~\ref{semi-continuity} and~\eqref{de-uinfty}-\eqref{id-hatu-UL}, we obtain that $\mathcal{Z}[U_L(t,\cdot)-w]\leq 1$ for all $t\in\R$. Notice actually that, since $U_L(t_*,x_*)=w(x_*)=\alpha$ and $U_L$ is not stationary, Lemma~\ref{zero-decrease} implies that~$\mathcal{Z}[U_L(t_*,\cdot)-w]=1$.

Next, we consider the case where $w$ is one of the types (d)-(e) in Lemma~\ref{classfy-w}. In this situation, it is also easily seen that
$$\mathcal{Z}[\hat{u}(0,\cdot;-kL)-w] =2,\,\, \hbox{ and }  \,\, SGN[\hat{u}(0,\cdot;-kL)-w]=[+,-,+] \,\,\hbox{ for all  large }k\in \N.$$
By using Lemmas~\ref{semi-continuity} and~\ref{zero-decrease} again, we have
$$\mathcal{Z}[U_L(t_*,\cdot)-w]\leq 2, \quad\hbox{and}  \quad SGN[U_L(t_*,\cdot)-w]\lhd  [+,-,+].$$
In particular, this implies that if $\mathcal{Z}[U_L(t_*,\cdot)-w]=2$, then $SGN[U(t_*,\cdot)-w]= [+,-,+]$. The proof of Lemma~\ref{intersection} is thus complete.
\end{proof}

In the following lemma, we give a special solution of~\eqref{ode} that decays to $0$ as $x\to+\infty$, with a decay rate that is controlled uniformly with respect to $L\geq L_*$. 

\begin{lem}\label{special-w}
For any $L\geq L_*$ and any $\delta\in (0,\delta_0]$, where $\delta_0\in(0,1/2)$ is the constant provided by {\rm (A2)}, there exists a solution $w\in C^2(\R)$ of~\eqref{ode} such that 
\begin{equation}\label{property-w}
w(0)=\delta,\quad \mathcal{Z}[w-1]=0,\quad w<1\hbox{ in $\R$},\quad\mathcal{Z}[w]\leq 1,
\end{equation}
and
\begin{equation}\label{uni-decay}
0<w(x)\leq \delta\me^{-\mu x} \,\,\hbox{ for all }x\geq 0,
\end{equation}
for some constant $\mu>0$ independent of $L$.
\end{lem}

\begin{proof}
For any $L\geq L_*$ and any $\delta\in (0,\delta_0]$, we show the existence of the desired solution by an approximation argument. For each $n\in\N$, let us first consider the following problem in bounded interval:
\begin{equation}\label{bd-ode} 
\left\{\baa{l}
(a_Lw')'+f_L(x,w)=0 \,\, \hbox{ in }[0,n],\vspace{5pt}\\
w(0) =\delta,\, \,w(n) =0.\eaa\right.
\end{equation}
Since $f_L(x,0)\equiv 0$ and since $f_L(x,\delta)<0$ for all $x\in\R$, it follows that $w=0$ and $w=\delta$ are, respectively, a sub-solution and a super-solution of~\eqref{bd-ode}.  Then, a classical iteration argument implies that there exists a solution $w_n \in C^2([0,n])$ such that $0\leq w_n\leq \delta$ in $[0,n]$. Furthermore, $0<w_n<\delta$ in $(0,n)$ from the strong maximum principle and, from the sliding method~\cite{bn}\footnote{For example, to prove $w_n\le w_m$ in $[0,n]$ when $n\le m$, the sliding method is used in the following way: one first proves that $w_n(\cdot+\tau)\leq w_m(\cdot)$ in $[0,n-\tau]$ when $\tau>0$ is close to $n$, and then slides $\tau$ to $0$.},  each~$w_n$ is unique and the sequence $(w_n)_{n\in\N}$ is nondecreasing in $n\in\N$, in the sense that $w_n\le w_m$ in $[0,n]$ if $n\le m$. Define now 
$$w(x):=\lim_{n\to+\infty} w_n(x)\,\,\hbox{ for }x\in [0,+\infty).$$
It then follows from the standard elliptic estimates that $w$ is a $C^2([0,+\infty))$ solution of 
$$(a_Lw')'+f_L(x,w)=0\hbox{ in }[0,+\infty),\ \ w(0)=\delta,\ \ 0\le w\le\delta\hbox{ in }[0,+\infty),$$
hence $0<w< \delta$ in $(0,+\infty)$ from the strong maximum principle.  Since $f_L(x,u)$ is continuous in $\R\times\R$ and globally Lipschitz-continuous in $u\in\R$ uniformly in $x\in\R$, the function $w$ can be extended to $x\in\R$ so that it is a $C^2(\R)$ solution of~\eqref{ode} with $w(0)=\delta$. It further follows from Lemma~\ref{classfy-w} that this solution must satisfy $\mathcal{Z}[w-1]=0$, $w<1$ in $\R$, and $\mathcal{Z}[w]\leq 1$ (more precisely, only cases~(b) or~(c) are possible). 

Now, to complete the proof, it remains to show the estimate~\eqref{uni-decay}. Since $0\le w_n\le\delta$ in $[0,n]$, by the first line of~\eqref{asspars}, we have, for each $n\in\N$,
$$(a_Lw'_n)'-\gamma_0 w_n \geq 0 \,\,\hbox{ for }x\in[0,n]. $$
Since $a_L(x)=a(x/L)$ is $L$-periodic and $a$ is positive and at least of class $C^1(\R)$, we can find a small constant $\mu>0$ (independent of $L\geq L_*$) such that $a_L(x)\mu^2-a_L'(x)\mu-\gamma_0\leq 0$ for all $x\in\R$. Letting $\bar{w}(x):=\delta\me^{-\mu x}$ for $x\geq 0$, we compute that 
$$ (a_L\bar{w}')'(x)-\gamma_0\bar{w}(x)= (a_L(x)\mu^2-a_L'(x)\mu-\gamma_0)\bar{w}(x) \leq 0 \,\,\hbox{ for all }x\geq 0.$$ 
It then follows from the elliptic weak maximum principle that $w_n(x)\le\bar{w}(x)=\delta\me^{-\mu x}$ for all $x\in[0,n]$ and $n\in\N$, hence $w(x)\le\delta\me^{-\mu x}$ for all $x\ge0$. This ends the proof of Lemma~\ref{special-w}.
\end{proof}

In the last lemma of this subsection, we consider the intersection of the pulsating front with the stationary solution of~\eqref{ode} obtained in Lemma~\ref{special-w}. 

\begin{lem}\label{inter-U-w}
Let $L\geq L_*$ be fixed, let $U_L$ be a pulsating front of~\eqref{eqL}, let $\delta\in(0,\delta_0]$ and let $w\in C^2(\R)$ be the solution of~\eqref{ode} provided by Lemma~$\ref{special-w}$. Then, $\mathcal{Z}[U_L(t,\cdot)-w]\leq  1$ for every $t\in\R$. Furthermore, if $\mathcal{Z}[U_L(t,\cdot)-w]= 1$ for some $t\in\R$, then $SGN[U_L(t,\cdot)-w]=[+,-]$.
\end{lem}

\begin{proof}
The proof is similar to that of Lemma~\ref{intersection}; therefore, we only give its outline. Let $t_*\in\R$ be arbitrary. Without loss of generality, we suppose that there exists some $x_*\in\R$ such that  $U_L(t_*,x_*)=w(x_*)$, hence $w(x_*)\in(0,1)$. To show the lemma, it suffices to prove that $\mathcal{Z}[U_L(t_*,\cdot)-w]\le 0$, or $\mathcal{Z}[U_L(t_*,\cdot)-w]=1$ and $ SGN[U_L(t_*,\cdot)-w]=[+,-]$. 

For any $z<x_*$, let $\hat{u}(t,x;z)$ be the solution of~\eqref{eqL} with Heaviside type initial condition $H(z-\cdot)$ and let $\tau(z)$ be defined as in~\eqref{de-tau} with $\alpha=w(x_*)$. By the proof of Lemma~\ref{intersection}, we see that $0<\tau(z)<+\infty$, $\lim_{z\to-\infty} \tau(z)=+\infty$, and that, up to extraction of a subsequence,
$$\hat{u}(t+\tau(-kL),x;-kL)\to U_L(t+t_*,x) \,\,\hbox{ as }k\to+\infty  \,\, \hbox{ in }C^{1;2}_{t;x;loc}(\R^2).$$
Since $w$ satisfies~\eqref{property-w} and $\lim_{x\to+\infty}w(x)=0$, it is straightforward to check that 
$$\mathcal{Z}[\hat{u}(0,\cdot;-kL)-w]=1 \quad\hbox{and}\quad SGN[\hat{u}(0,\cdot;-kL)-w]=[+,-]$$
for all large $k\in\N$.  Passing to the limit as $k\to+\infty$, it follows from Lemmas~\ref{semi-continuity} and~\ref{zero-decrease} that~$\mathcal{Z}[U_L(t_*,\cdot)-w] \leq 1$ and $SGN[U_L(t_*,\cdot)-w]\lhd[+,-]$, yielding the desired results.
\end{proof}


\subsection{Proof of Theorem~\ref{th-wave-length}}\label{sec23}

For clarity, we divide the proof into three steps. In our arguments below, $\delta\in (0,1/2]$ is given. 

\medskip 
\noindent{\it Step 1:~\eqref{eq-wave-lengh-t} holds true with $C>0$ independent of $L\ge L_*$}.  For any $L\geq L_*$ and $x\in\R$, since $U_L(t,x)$ is continuous and increasing in $t$ with $U_L(-\infty,x)=0$ and $U_L(+\infty,x)=1$, the set
$$I_{L}(x):=\{t\in\R:\, \delta\leq U_{L}(t,x)\leq 1-\delta\}$$
is a compact interval in $\R$. Denote by $m(I_{L}(x))$ the length of this interval. It follows from Theorem~\ref{th-known}~(ii) and~(iv) that $m(I_{L}(x))$ is continuous with respect to $(L,x)\in [L_*,+\infty)\times\R$. 

Assume by contradiction that the desired result~\eqref{eq-wave-lengh-t} is not true. Then, there exist sequences $(L_n)_{n\in\N} \subset [L_*,+\infty)$ and $(x_n)_{n\in\N}\subset \R$ such that 
\begin{equation}\label{unbound-I}
m(I_{L_n}(x_n))=m\big(\{t\in\R:\, \delta\leq U_{L_n}(t,x_n)\leq 1-\delta\}\big)\to+\infty \,\,\hbox{ as } n\to+\infty.
\end{equation}
From~\eqref{UL-period}, one easily sees that $m(I_{L_n}(x_n))=m(I_{L_n}(x_n+L_n))$. Therefore, without loss of generality, we may assume that $(x_n)_{n\in\N} \subset [0,L_n)$. Moreover, one observes that $L_n\to+\infty$ as $n\to+\infty$. Otherwise, there would exist subsequences $(L_n')\subset (L_n)$ and $(x_n')\subset (x_n)$ such that as $n\to+\infty$, $L_n'\to L_*'\in [L_*,+\infty)$ and $x_n'\to x_*' \in  [0,L_*']$, whence by the continuity of $m(I_{L}(x))$ in $(L,x)$, there would hold $m(I_{L'_n}(x_n')) \to m(I_{L'_*}(x_*'))$ as $n\to\infty$, which would contradict \eqref{unbound-I}.

For each $n\in\N$, since $U_{L_n}(\cdot,x_n)$ is increasing and continuous, with $U_{L_n}(-\infty,x_n)=0$ and $U_{L_n}(+\infty,x_n)=1$, there exists a unique $t_n\in I_{L_n}(x_n)$ such that $U_{L_n}(t_n,x_n)=1/2$.
Define
$$\tilde{U}_{n}(t,x):=U_{L_n}(t+t_n,x+x_n) \,\,\hbox{ for }(t,x)\in\R^2.$$
Clearly, each function $\tilde{U}_{n}$ is an entire solution of
\begin{equation}\label{eq-shift}
\partial_t\tilde{U}_{n}=\partial_x\left(a\left(\frac{x+x_n}{L_n}\right)\partial_x\tilde{U}_{n}\right)+f\left(\frac{x+x_n}{L_n},\tilde{U}_{n}\right) \,\,\hbox{ in }\R^2,
\end{equation}
with $\tilde{U}_{n}(0,0)=1/2$. Up to extraction of some subsequence, one can assume that $x_n/L_n \to x_{\infty}\in [0,1]$, and that, from standard parabolic estimates, there exists a function $\tilde{U}_{\infty}\in C^{1;2}_{t;x}(\R^2,[0,1])$ such that
\begin{equation}\label{ULn-Uinfty}
\tilde{U}_{n} \to \tilde{U}_{\infty}\,\,\hbox{ in }C^{1;2}_{t;x;loc}(\R^2)\hbox{ as }n\to+\infty.
\end{equation}
Moreover, $\tilde{U}_{\infty}$ satisfies $\tilde{U}_{\infty}(0,0)=1/2$ and it is a solution to
\begin{equation}\label{homo-parabolic}
\partial_t\tilde{U}_{\infty} =a(x_{\infty}) \partial_{xx}\tilde{U}_{\infty}+ f(x_{\infty},\tilde{U}_{\infty})\,\, \hbox{ in }\R^2.
\end{equation}
By the strong maximum principle, we have $0<\tilde{U}_{\infty}<1$ in $\R^2$. Furthermore, since $\partial_t\tilde{U}_{n}>0$ in~$\R^2$, we have $\partial_t \tilde{U}_{\infty}\geq 0$ in $\R^2$. Finally, one infers from~\eqref{unbound-I} that
\begin{equation}\label{inifite-lengh}
m\big(\{t\in\R:\, \delta\leq \tilde{U}_{\infty}(t,0)\leq 1-\delta\}\big)=+\infty.
\end{equation}

By the monotonicity of $\tilde{U}_{\infty}$ in $t$ and standard parabolic estimates applied to equation~\eqref{homo-parabolic}, it follows that there exist two steady states $0\leq p_-\leq p_+ \leq 1$ of~\eqref{homo-parabolic} such that
\begin{equation}\label{Uinfty-ppm}
\tilde{U}_{\infty}(t,\cdot) \to p_{\pm}\,\,\hbox{ as }t\to \pm\infty\hbox{ in } C^2_{loc}(\R).
\end{equation}
It is clear that $p_{\pm}$ satisfy
\begin{equation}\label{homo-elliptic}
a(x_{\infty}) p_\pm''+ f(x_{\infty}, p_\pm)=0\,\,\hbox{ in }\R,
\end{equation}
and $0\le p_{-}(0)\leq 1/2\leq p_{+}(0)\le 1$.
Furthermore, one can conclude that either $\delta\leq p_{-}(0)\leq 1/2$ or $1/2\leq p_{+}(0)\leq 1-\delta$. (Otherwise, one would have $0\leq p_{-}(0)< \delta$ and $1-\delta< p_{+}(0)\leq 1$, whence by \eqref{Uinfty-ppm}, there would hold $m\big(\{t\in\R:\, \delta\leq \tilde{U}_{\infty}(t,0)\leq 1-\delta\}\big)<+\infty$, which would contradict  \eqref{inifite-lengh}.) Then, by the strong maximum principle, we have
$$\hbox{either }0<p_-<1\hbox{ in }\R, \,\,\hbox{ or }\,\, 0<p_+<1\hbox{ in }\R.$$
Without of loss generality, we assume that the former case happens, since the latter case can be handled in a similar way. Then, thanks to the assumptions (A1) and (A3), according to the phase diagrams of equation~\eqref{homo-elliptic}, the solution $p_-$ can only be one of the following three types: either a constant function, or a non-constant periodic function, or a ground state solution such that $p_-(x)\to 0$ as $x\to\pm\infty$. We will derive a contradiction in each of these three cases.

{\it Case 1: $p_-$ is a constant solution, that is, $p_-\equiv b(x_{\infty})$ in $\R$.}  In this case, for any non-constant periodic solution $0<q<1$ of~\eqref{homo-elliptic}, one finds some points $y_1<z_1<y_2<z_2$ and some constant $\epsilon_0>0$ such that
\begin{equation}\label{com-p-q}
q(y_i)- p_-(y_i) \geq \epsilon_0, \quad\hbox{and} \quad q(z_i)- p_-(z_i) \leq -\epsilon_0 \,\,\hbox{ for }i=1,2.
\end{equation}
Remember that function $u\mapsto f(x,u)$ is globally Lipschitz-continuous in $u\in\R$ uniformly in $x\in\R$. Then, for each $n\in\N$, it follows from the classical ODE theory that the following problem
\begin{equation}\label{elliptic-n}
\left(a\left(\frac{x+x_n}{L_n}\right)w'\right)'+f\left(\frac{x+x_n}{L_n},w\right)=0 \,\,\hbox{ in }\R,
\end{equation}
with $w(y_1)=q(y_1)$ and $w'(y_1)=q'(y_1)$, admits a unique solution $w_n\in C^2(\R)$. Moreover, by standard elliptic estimates, we have $w_n\to q$ as $n\to+\infty$ in $C^2_{loc}(\R)$. This implies in particular that there exists $N\in\N$ such that
$$|w_n(x)- q(x)| \leq \frac{\epsilon_0}{4}  \,\,\hbox{ for all }x\in [y_1,z_2]\hbox{ and }n\geq N.$$
On the other hand, it follows from~\eqref{ULn-Uinfty} and~\eqref{Uinfty-ppm} that there exists $s_*<0$, with $|s_*|$ sufficiently large, such that, replacing $N$ by a larger integer if necessary,
$$\left| \tilde{U}_{n}(s_*,x)-p_-(x)\right|\leq  \frac{\epsilon_0}{4}  \,\,\hbox{ for all }x\in [y_1,z_2]\hbox{ and }n\geq N.$$
Combining the above,  we obtain that, for all large $n$,
$$w_n(y_i)- \tilde{U}_{n}(s_*,y_i) \geq \frac{\epsilon_0}{2},\quad\hbox{and} \quad w_n(z_i)- \tilde{U}_{n}(s_*,z_i) \leq -\frac{\epsilon_0}{2} \,\,\hbox{ for }i=1,2.$$
By continuity, the function $x\mapsto  \tilde{U}_{n}(s_*,x)-w_n(x)$ has at least three sign changes in $[y_1,z_2]$, which is a contradiction with the fact that $\mathcal{Z}[\tilde{U}_{n}(s_*,\cdot)-w_n] \leq 2$ (by Lemma~\ref{intersection} applied with~$\tilde{U}_n(\cdot,\cdot-x_n)$ and $w_n(\cdot-x_n)$). Hence, Case 1 is ruled out.

{\it Case 2: $p_-$ is a non-constant periodic solution.}  In this case,  letting $q(x)\equiv b(x_{\infty})$, we see that $q$ is a constant solution of~\eqref{homo-elliptic}. Then, we can find some points $y_1<z_1<y_2<z_2$ and a constant $\epsilon_0>0$ such that~\eqref{com-p-q} holds true, and hence, the same reasoning as in Case 1 yields a contradiction. Therefore, Case 2 is ruled out too.

{\it Case 3: $p_-$ is a ground state solution with $p_-(x)\to 0$ as $x\to\pm\infty$.} In this case, $p_-(x)$ is symmetrically decreasing and there exist some points $\tilde{y}_1<\tilde{y}_0<\tilde{y_2}$ and a constant $\tilde{\epsilon}_0>0$ such that
\begin{equation*}
b(x_{\infty})-p_-(\tilde{y}_i) \geq \tilde{\epsilon}_0  \hbox{ for } i=1,2,\quad\hbox{and}\quad b(x_{\infty})-p_-(\tilde{y}_0) \leq -\tilde{\epsilon}_0. 
\end{equation*}
Similarly as in Case 1, for each $n\in\N$, problem~\eqref{elliptic-n} with $w(\tilde{y}_1)=b(x_{\infty})$ and $w'(\tilde{y}_1)=0$ has a unique solution $\tilde{w}_n\in C^2(\R)$, and from standard elliptic estimates, $\tilde{w}_n \to b(x_{\infty})$ as $n\to+\infty$ in~$C^2_{loc}(\R)$. As a consequence, there exists a large integer $\tilde{N}\in \N$ such that
$$\left|\tilde{w}_n(x)- b(x_{\infty})\right| \leq \frac{\tilde{\epsilon}_0}{4}  \,\,\hbox{ for all }x\in [\tilde{y}_1,\tilde{y}_2]\hbox{ and }n\geq \tilde{N}.$$
On the other hand, by using~\eqref{ULn-Uinfty} and~\eqref{Uinfty-ppm} again and by making some adjustment to $\tilde{N}$ if necessary, one finds some negative time $\tilde{s}_*<0$, with $|\tilde{s}_*|$ large enough, such that
$$\left| \tilde{U}_{n}(\tilde{s}_*,x)-p_-(x)\right|\leq  \frac{\tilde{\epsilon}_0}{4}  \,\,\hbox{ for all }x\in [\tilde{y}_1,\tilde{y}_2]\hbox{ and }n\geq \tilde{N}.$$
Therefore, we obtain that, for each $n\geq \tilde{N}$,
\begin{equation}\label{tilde-wn-ULn}
\tilde{w}_n(\tilde{y}_i)- \tilde{U}_{n}(\tilde{s}_*,\tilde{y}_i) \geq \frac{\tilde{\epsilon}_0}{2}  \hbox{ for }i=1,2,\quad\hbox{and}\quad \tilde{w}_n(\tilde{y}_0)- \tilde{U}_{n}(\tilde{s}_*,\tilde{y}_0) \leq -\frac{\tilde{\epsilon}_0}{2}.
\end{equation}
Remember that by Lemma~\ref{intersection} (applied with $\tilde{U}_n(\cdot,\cdot-x_n)$ and $\tilde{w}_n(\cdot-x_n)$), we have  $\mathcal{Z}[\tilde{U}_{n}(\tilde{s}_*,\cdot)-\tilde{w}_n] \leq 2$. Then,~\eqref{tilde-wn-ULn} implies that
$$\mathcal{Z}[\tilde{U}_{n}(\tilde{s}_*,\cdot)-\tilde{w}_n] = 2,\quad\hbox{and}\quad SGN[\tilde{U}_{n}(\tilde{s}_*,\cdot)-\tilde{w}_n]=[-,+,-].$$
This last property contradicts property~\eqref{sgn-UL-w} of Lemma~\ref{intersection}, hence Case 3 is ruled out too.

Consequently, the case where $0<p_-<1$ cannot happen. Furthermore, with similar arguments, one can exclude the case where $0<p_+<1$. Therefore, our assumption~\eqref{unbound-I} at the beginning of the proof is unreasonable. This ends the proof of Step 1.

\medskip
\noindent{\it Step 2:~\eqref{derivative-t} holds true with $\beta>0$ independent of $L\ge L_*$}. Assume by contradiction that there exist sequences $(\bar{L}_n)_{n\in\N} \subset [L_*,+\infty)$ and $(\bar{t}_n,\bar{x}_n)_{n\in\N}\subset \R^2$ such that $U_{\bar{L}_n}(\bar{t}_n,\bar{x}_n) \in [\delta,1-\delta]$ for each $n\in\N$, and
$$\partial_t U_{\bar{L}_n}(\bar{t}_n,\bar{x}_n) \to 0\,\,\hbox{ as }n\to+\infty. $$
Due to~\eqref{UL-period}, we may assume without loss of generality that $0\leq \bar{x}_n<\bar{L}_n$. By Theorem~\ref{th-known}~(ii)~(iv) and standard parabolic estimates, it follows that $\bar{L}_n\to+\infty$ as $n\to+\infty$.

Define now
$$\bar{U}_{n}(t,x):=U_{\bar{L}_n}(t+\bar{t}_n,x+\bar{x}_n) \,\,\hbox{ for } (t,x)\in\R^2.$$
Up to extraction of some subsequence, we may assume that the sequence $(\bar{x}_n/\bar{L}_n)_{n\in\N}\subset [0,1]$ converges as $n\to+\infty$. By a slight abuse of notation, we still denote by $x_\infty$ this limit. Then, standard parabolic estimates imply that, possibly up to extraction of a further subsequence, $\bar{U}_{n}\to \bar{U}_{\infty}$ in $C^{1;2}_{t;x;loc}(\R^2)$, where $\bar{U}_{\infty}$ is an entire solution of~\eqref{homo-parabolic}. It is clear that $\partial_t \bar{U}_{\infty}\geq 0$ in $\R^2$ and $\partial_t \bar{U}_{\infty}(0,0)=0$. Furthermore, the strong maximum principle applied to the equation satisfied by $\partial_t \bar{U}_{\infty}$ immediately gives that $\partial_t \bar{U}_{\infty}\equiv 0$ in~$\R^2$. In other words, $\bar{U}_{\infty}(t,x)=\bar{U}_\infty(x)$ is independent of $t$ and it obeys $a(x_\infty)\bar{U}_{\infty}''+f(x_\infty,\bar{U}_{\infty})=0$ in $\R$. Notice that  $0\leq \bar{U}_{\infty} \leq 1$ in $\R$ and $\delta\leq \bar{U}_{\infty}(0)\leq 1-\delta$. It then follows from the strong maximum principle that $0<\bar{U}_{\infty}<1$ in $\R$. Furthermore, it follows from~(A1) and~(A3) that either $\bar{U}_{\infty}$ is a constant solution (i.e., $\bar{U}_{\infty}\equiv b(x_{\infty})$), or it is a non-constant periodic solution, or it is a ground state solution decaying to $0$ as $x\to\pm\infty$.  Proceeding similarly as in Step 1, we can find a contradiction in each of these cases: namely, for each large $n\in\N$, there exists then a solution $\bar{w}_n$ of~\eqref{elliptic-n} (with $\bar{L}_n$ and $\bar{x}_n$ instead of $L_n$ and $x_n$) such that either
$$\mathcal{Z}[\bar{U}_{n}(0,\cdot)-\bar{w}_n] \geq 3,$$
or
$$\mathcal{Z}[\bar{U}_{n}(0,\cdot)-\bar{w}_n] =2\ \ \hbox{and}\ \ SGN[\bar{U}_{n}(0,\cdot)-\bar{w}_n]=[-,+,-]$$
(the proof is even simpler than that in Step 1, as the limit function $\bar{U}_{\infty}$ is stationary; therefore we do not repeat the details). Both situations are impossible, due to Lemma~\ref{intersection} (applied with~$\bar{U}_n(\cdot,\cdot-\bar{x}_n)$ and $\bar{w}_n(\cdot-\bar{x}_n)$). As a consequence, the proof of Step~2 is complete.  

\medskip 
\noindent{\it Step 3:~\eqref{eq-wave-lengh-x} holds true with $C>0$ independent of $L\ge L_*$}. Assume by contradiction that there does not exist such a bound  independent of $L$. Then, we can find sequences $(L_n)_{n\in \N} \subset [L_*, +\infty)$, $(s_n)_{n\in\N}\subset \R$, $(\alpha_{n})_{n\in\N},\,(\gamma_{n})_{n\in\N} \subset [\delta,1-\delta]$ and $(x_n)_{n\in\N},\, (y_n)_{n\in\N}\subset \R$ such that   
$$U_{L_n}(s_n,x_n)=\alpha_n\,\, \hbox{ and }\,\,U_{L_n}(s_n,y_n)=\gamma_n\,\,\hbox{ for each }n\in\N, $$ 
and
$$|x_n-y_n| \to +\infty\,\, \hbox{ as }n\to+\infty.$$ 
Without loss of generality, for each $n\in\N$, we may assume that $y_n>x_n$,  and thanks to~\eqref{UL-period}, we may also assume that $0\leq x_n<L_n$. Moreover, by Theorem~\ref{th-known}~(ii)~(iv) again, we have 
$L_n\to+\infty$ as $n\to+\infty$.

Define now
$$\tilde{U}_{n}(t,x):=U_{L_n}(t+s_n,x+x_n) \,\,\hbox{ for }(t,x)\in\R^2.$$
Clearly, each function $\tilde{U}_{n}(t,x)$ is an entire solution of~\eqref{eq-shift} with $\tilde{U}_{n}(0,0)=\alpha_n$ and $\tilde{U}_{n}(0,z_n)=\gamma_n$, where $0<z_n=y_n-x_n\to+\infty$ as $n\to+\infty$. Up to extraction of some subsequence, we have 
$$\frac{x_n}{L_n}\to x_{\infty} \in [0,1],\quad \alpha_{n} \to \alpha_{\infty} \in [\delta,1-\delta],\quad \gamma_{n} \to \gamma_{\infty} \in [\delta,1-\delta],  \,\,\,\, \hbox{ as }n\to+\infty,$$
and by standard parabolic estimates, $\tilde{U}_{n} \to \tilde{U}_{\infty}$ in $C^{1;2}_{t;x;loc}(\R^2)$ as $n\to+\infty$, where $\tilde{U}_{\infty}$ is a solution of~\eqref{homo-parabolic} with $\tilde{U}_{\infty}(0,0)=\alpha_{\infty}$. 

Set $\delta_*:=\min\{\delta_0,\delta\}$, where $\delta_0$ is the constant provided by (A2). It follows from Lemma~\ref{special-w} that, for each $n\in\N$,~\eqref{ode} with $L=L_n$ admits a solution $w_n\in C^2(\R)$ such that $w_n(0)=\delta_*$,  $\mathcal{Z}[w_n-1]=0$, $\mathcal{Z}[w_n]\leq 1$ and $0<w_n(x)\leq \delta_*\me^{-\mu x}$ for all $x\geq 0$, where $\mu>0$ is independent of~$n\in\N$. Notice that for each $n\in\N$, $\tilde{U}_n(t,x)$ is increasing and continuous in $t\in\R$, with $\tilde{U}_{n}(-\infty,x)=0$ and $\tilde{U}_{n}(+\infty,x)=1$ locally uniformly in $x\in\R$. This implies in particular that there exists a unique time $\tau_n \in\R $ such that $\tilde{U}_{n}(\tau_n,z_n)=w_n(z_n)$.  Since $w_n(z_n)<\delta_*\leq\delta\leq  \gamma_n=\tilde{U}_{n}(0,z_n)$ for each $n\in\N$, it is clear that $\tau_n<0$. 

Now, we claim that 
\begin{equation}\label{sn-infty}
\tau_n\to-\infty\,\,\hbox{ as }n\to+\infty.  
\end{equation}
Suppose to the contrary that $(\tau_n)_{n\in\N}$ converges, up to extraction of a subsequence, to $\tau_\infty\in(-\infty,0]$. For each $n\in\N$, let us write 
$$z_n:=z_n'+z_n'',\quad \hbox{and}\quad  V_n(t,x):=\tilde{U}_n(t,x+z_n) \,\,\hbox{ in }\R^2, $$
where $z_n'\in L_n\N$ and $z_n'' \in [0,L_n)$. It is easily checked that each $V_n$ is an entire solution of 
$$ \partial_tV_n=\partial_x\left(a\left(\frac{x+x_n+z''_n}{L_n}\right)\partial_xV_n\right)+f\left(\frac{x+x_n+z''_n}{L_n},V_n\right) \,\,\hbox{ in }\R^2, $$
with $V_n(0,0)=\gamma_n$ and $V_n(\tau_n,0)=w_n(z_n)$.  Since $z_n\to+\infty$ as $n\to+\infty$ and since $w_n(x)$ decays to $0$ as $x\to+\infty$ uniformly in $n\in\N$, it follows that $\lim_{n\to+\infty} V_n(\tau_n,0)=0$.  Next, up to extraction of another subsequence, we may assume that $z_n''/L_n \to z_{\infty} \in [0,1]$ as $n\to+\infty$, and that, by standard parabolic estimates, $V_n\to V_{\infty}$ as $n\to+\infty$ in $C^{1;2}_{t;x;loc}(\R^2)$, where $V_{\infty}$ is an entire solution of  
$$\partial_t V_{\infty} =a(x_{\infty}+z_{\infty}) \partial_{xx} V_{\infty}+ f(x_{\infty}+z_{\infty}, V_{\infty})\,\, \hbox{ in }\R^2.$$
Since $0\leq V_{\infty}\leq 1$ and since $V_{\infty}(\tau_{\infty},0)=\lim_{n\to+\infty}V_n(\tau_n,0)=0$, it follows from the strong maximum principle that $V_{\infty}\equiv 0$ in $\R^2$. This contradicts the fact that $V_{\infty}(0,0)=\lim_{n\to+\infty}\gamma_n=\gamma_{\infty}\geq \delta$.  Therefore, our claim~\eqref{sn-infty} is proved.

We are now ready to derive a contradiction.  By the proof of Step 1, it follows that the sequence $\big(m(\{t\in\R: \delta_* \leq \tilde{U}_n(t,0)\leq 1-\delta_*\})\big)_{n\in\N}$ is bounded uniformly in $n\in\N$. This together with $\tilde{U}_n(0,0)=\alpha_n\in[\delta,1-\delta]\subset[\delta_*,1-\delta_*]$ and $\partial_t\tilde{U}_n>0$ yields the existence of $A>0$ such that~$\{t\in\R: \delta_* \leq \tilde{U}_n(t,0)\leq 1-\delta_*\}\subset[-A,A]$ for all $n\in\N$. As a consequence, one finds~$t_*<0$ such that  
$$\tilde{U}_n(t_*,0)< \delta_* =w_{n}(0) \,\,\hbox{ for all }n\in\N. $$
Furthermore, by~\eqref{sn-infty} and the positivity of $\partial_t\tilde{U}_n$, there exists $n_*\in\N$ sufficiently large such that 
$$\tilde{U}_{n_*}(t_*,z_{n_*}) > w_{n_*}(z_{n_*}).$$
Combining the above, we immediately obtain that 
$$\mathcal{Z}[\tilde{U}_{n_*}(t_*,\cdot)- w_{n_*}] \geq 1\quad\hbox{and}\quad SGN[\tilde{U}_{n_*}(t_*,\cdot)- w_{n_*}]\rhd [-,+],$$
which is a contradiction with the conclusions of Lemma~\ref{inter-U-w} (applied with $\tilde{U}_{n_*}(\cdot,\cdot-x_{n_*})$ and $w_{n_*}(\cdot-x_{n_*})$).

This ends the proof of Step 3, and the proof of Theorem~\ref{th-wave-length} is thus complete.\hfill$\Box$


\subsection{Proof of Theorem~\ref{wave-length-1}}\label{sec24}

We first show statements (i). From Theorem~\ref{th-known}~(ii), it suffices to prove that for any $\delta \in (0,1/2]$, the interval
$$E_{L}(y):=\big\{\xi\in\R:\, \delta\leq \phi_L(\xi+\zeta_L(y),y) \leq 1-\delta\big\}$$
is bounded uniformly with respect to $y\in\R$ and $L\geq L_*$. To do so, observe that the formula $U_L(t,x)=\phi_L(x-c_Lt,x/L)$ implies that
\begin{equation}\label{redefine-UL}
U_L\left(-\frac{\xi}{c_L}+\frac{Ly-\zeta_L(y)}{c_L}, Ly \right)= \phi_L(\xi+\zeta_L(y),y)\,\,\hbox{ for all }\xi\in\R\hbox{ and }y\in \R.
\end{equation}
Clearly, $U_L((Ly-\zeta_L(y))/c_L, Ly)=1/2$ for all $y\in\R$, by~\eqref{normmal-front-1}. Then, by Theorem~\ref{th-wave-length}, the following interval
$$\tilde{I}_{L}(Ly):=\left\{t\in\R:\, \delta\leq U_L\left(t+\frac{Ly-\zeta_L(y)}{c_L},Ly\right)\leq 1-\delta   \right\} $$
is bounded uniformly in $L\geq L_*$ and $y\in \R$. Remember that $c_L$ is bounded in $L\geq L_*$, by Theorem~\ref{th-known}~(i). This together with~\eqref{redefine-UL}  immediately implies that the interval $E_{L}(y)$ is bounded uniformly in $y\in\R$ and $L\geq L_*$. The proof of statements~(i) of Theorem~\ref{wave-length-1} is thus complete.

Next, we show statement (ii). It follows from~\eqref{redefine-UL} that for any $(\xi,y)\in \R\times \R$ with $\xi \in E_{L}(y)$, we have $-\xi/c_L \in\tilde{I}_{L}(Ly)$ and
$$\partial_{\xi} \phi_L(\xi+\zeta_L(y),y)=-\frac{1}{c_L} \partial_t U_L \left(-\frac{\xi}{c_L}+\frac{Ly-\zeta_L(y)}{c_L},Ly\right). $$
Then by Theorem~\ref{th-wave-length} and the fact that $c_L>0$ is bounded in $L\geq L_*$, we immediately obtain statement~(ii) (with a constant $\beta>0$ which is in general different from the one appearing in~\eqref{derivative-t}). 

Finally, we prove statement (iii). Since $\zeta_L(y)$ is continuous in $(L,y)\in [L_*,+\infty)\times\R$, it is clear that $(L,y)\mapsto \zeta_L(y)/L$ is also continuous. It remains to show that this function is bounded. Since $\zeta_L(y)$ is $1$-periodic in $y\in\R$, we only need to show that $\zeta_L(y)/L$ is bounded in $(L,y)\in [L_*,+\infty)\times [0,1]$. It follows from~\eqref{normal-U} and Theorem~\ref{th-wave-length} that the set $\{x\in\R: U_L(0,x)=1/2\}$ is bounded uniformly in $L\geq L_*$, that is, there exists $\bar{L}\geq L_*$ such that 
$$\big\{x\in\R: U_L(0,x)=1/2\big\} \subset [-\bar{L},\bar{L}] \,\,\hbox{ for all }L\geq L_*.$$
Remember that $\lim_{x\to+\infty}U_L(0,x)=0$ for each $L\geq L_*$. This implies that $U_L(0,L(y+1)) <  1/2$  for all $y\in [0,1]$ and $L> \bar{L}$. Thanks to~\eqref{redefine-UL} applied with $\xi:=L(y+1)-\zeta_L(y+1)$ and $y+1$ instead of $y$, it follows that $\phi_L(L(y+1)-\zeta_L(y+1)+\zeta_L(y),y)=\phi_L(L(y+1)-\zeta_L(y+1)+\zeta_L(y),y+1)=U_L(0,L(y+1))<1/2$. Hence, the monotonicity of $\phi_L$ in its first variable implies that
$$L(y+1)-\zeta_L(y+1)>0 \,\,\hbox{ for all } y\in [0,1]\hbox{ and }L> \bar{L}. $$
Consequently, by the periodicity of $\zeta_L(y)$, we obtain $\zeta_L(y)/L<y+1$ for all $y\in [0,1]$ and $L> \bar{L}$. 

In a similar way, by using the fact that  $\lim_{x\to-\infty}U_L(0,x)=1$ and the monotonicity of $U_L$ in its first variable, we can prove that 
$$L(y-2)-\zeta_L(y-2)<0 \,\,\hbox{ for all } y\in [0,1]\hbox{ and }L> \bar{L},$$
hence $\zeta_L(y)/L>y-2$  for all $y\in [0,1]$ and $L> \bar{L}$. Combining the above, we immediately obtain that $-2<\zeta_L(y)/L<2$ for all $y\in [0,1]$ and $L> \bar{L}$. This together with the continuity of $\zeta_L(y)/L$ in $[L_*,+\infty)\times [0,1]$ implies that $ \zeta_L(y)/L$ is bounded in $[L_*,+\infty)\times [0,1]$. The proof of Theorem~\ref{wave-length-1} is thus complete.\hfill$\Box$


\section{Convergence of $c_L$ as $L\to+\infty$: proof of Theorem~\ref{speed-limit}}\label{sec3}

This section is devoted to the proof of Theorem~\ref{speed-limit}. When $L$ is large, we consider the Cauchy problem associated with equation~\eqref{eqL} under the large-time and large-space scaling, $t\to Lt$, $x\to Lx$. More precisely, letting $v(t,x)=u(Lt,Lx)$ with $u(t,x)$ being a solution of~\eqref{eqL}, the rescaled Cauchy problem reads
\begin{equation}\label{eq-scale}
\left\{\baa{ll}
\displaystyle v_t= \frac{1}{L}(a(x)v_x)_x+Lf(x,v),&  t>0,\,x\in\R,\vspace{5pt}\\
v(0,x)=g(x), & x\in\R.
\eaa \right.
\end{equation}
We consider initial conditions $g\in C(\R,[0,1])$ which are independent of $L$ and front-like, in the sense that
\begin{equation}\label{i-condition-2}
\liminf_{x\to-\infty}\,(g(x)-b(x))>0 ,\quad  \limsup_{x\to+\infty}\,(g(x)-b(x)) <0,
\end{equation}
where $x\mapsto b(x)\in (0,1)$ is the zero of $f(x,\cdot)$ provided by (A1). For each $L>0$, denote by~$v_L(t,x)$ the solution of the Cauchy problem~\eqref{eq-scale}. By the strong maximum principle, we have $0<v_L(t,x)<1$  for all $t>0$ and $x\in\R$. Moreover, it is easily seen from~\eqref{de-pulsating} that for each $L\geq L_*$, with $L_*>0$ as in Theorem~\ref{th-known}, $\phi_L(L(x-c_Lt),x)$ is a solution of the first equation of~\eqref{eq-scale}, and it is defined for all $(t,x)\in\R^2$.

As we mentioned in Section~\ref{intro}, we use G\"artner's result~\cite{g} to show the convergence of $c_L$ as $L\to+\infty$. Recall that the function $c:\R\to \R$, $y\mapsto c(y)$ denotes the unique front speeds of the family of homogeneous equations~\eqref{eq-homo}. The following lemma is a direct application of~\cite[Corollary, p.~140]{g} to our one-dimensional spatially periodic equation~\eqref{eq-scale} (see also Xin's review paper~\cite[Theorem~3.7]{x3} in the case where $g$ is compactly supported):

\begin{lem}\label{singular-limit}{\rm{\cite{g}}}
Let $g\in C(\R,[0,1])$ satisfy~\eqref{i-condition-2} and let $\Gamma_0=\{x\in\R: \, g(x)=b(x)\}$. For any $L>0$, let $v_L(t,x)$ be the solution of~\eqref{eq-scale} with initial condition $g$ $($$g$ is independent of $L$$)$. Then, as $L\to+\infty$,
\begin{equation*}
\left\{\begin{array}{ll}
v_L(t,x) \to 1 &  \mbox{locally uniformly in }  \,  \{(t,x)\in (0,+\infty)\times\R: \rho(x,\Gamma_0)<t \}, \vspace{5pt}\\
v_L(t,x) \to 0  &  \mbox{locally uniformly in }  \, \{(t,x)\in (0,+\infty)\times\R: \rho(x,\Gamma_0)>t \},
\end{array}\right.
\end{equation*}
where $\rho(x,\Gamma_0):=\inf\{\rho(x,x_2):x_2\in\Gamma_0\}$ and  $\rho(\cdot,\cdot)$ denotes the signed distance function defined by
$$\rho(x_1,x_2):=\int_{x_2}^{x_1}  c^{-1}(y) dy  \,\, \hbox{ for }x_1\in\R,\,x_2\in\R.  $$
\end{lem}

Recall that $c_*$ is the constant defined in~\eqref{limit-speed}. Since the function $y\mapsto c(y)$ is 1-periodic, the above result suggests that as $L\to +\infty$, the solution $v_L(t,x)$ propagates at a mean speed equal to~$c_*$. This however does not imply directly the convergence of the speeds $c_L$ to $c_*$, since the profiles of the pulsating fronts, say at time $0$, depend on $L$, unlike the initial condition $g$ in Lemma~\ref{singular-limit}. To circumvent this difficulty, we establish in Lemma~\ref{sub-super} some further comparisons between the solutions $v_L(t,x)$ and the pulsating fronts $\phi_L(L(x-c_Lt),x)$ for $L\ge L_*$, which easily lead to the proof of Theorem~\ref{speed-limit}, on the basis of the uniform estimates proved in Section~\ref{sec2}.\footnote{An alternate approach could be to use~\cite[Theorem~4.1]{g} with $L$-dependent initial conditions $\phi_L(L\cdot,\cdot)$ combined with Theorem~\ref{wave-length-1}. We here use Lemma~\ref{singular-limit} and the comparisons established in Lemma~\ref{sub-super}, since the method of proof based on sub- and super-solutions will serve as an archetype to get more involved comparisons in Section~\ref{sc-front-limit-1}.} To do so, let $\delta_0$ and $\gamma_0$ be the positive constants provided by the assumption (A2). Since $f(\cdot,0)=f(\cdot,1)=0$ in $\R$, there is $\epsilon_0 \in (0,\delta_0/2]$ such that
\begin{equation}\label{p-derivative-u}
\partial_u f(x,u) \leq -\frac{\gamma_0}{2}\,\,\hbox{ for all }u\in [-2\epsilon_0,2\epsilon_0]\cup [1-2\epsilon_0,1+2\epsilon_0]\hbox{ and }x\in\R.
\end{equation}
We now choose a special initial condition $g\in C(\R,[0,1])$, independent of $L\geq L_*$, such that
\begin{equation}\label{i-condition-3}
\liminf_{x\to-\infty}\,g(x)> 1-\epsilon_0,\quad  \limsup_{x\to+\infty}\,g(x) < \epsilon_0,\quad\hbox{and}\quad \Gamma_0=\{0\}.
\end{equation}
Notice that it is for instance possible to choose $g$ in such a way that $g=1$ in $(-\infty,-A]$ and $g=0$ in $[A,+\infty)$, for some large~$A$.

\begin{lem}\label{sub-super}
For any $L\geq L_*$, with $L_*>0$ as in Theorem~$\ref{th-known}$, let $(\phi_L,c_L)$ be a pulsating front of~\eqref{eqL} and let $v_L(t,x)$ be the solution of~\eqref{eq-scale} with the initial condition $g$ satisfying~\eqref{i-condition-3}. Then there exist constants $C_{\pm}\geq 0$ and $K_0>0$ $($all independent of $L$$)$ such that
\begin{equation}\label{ineq-sub-super}
\left\{\begin{array}{l}
v_L(t,x) \geq \phi_L(L(x-c_Lt+C_-)+K_0\epsilon_0,x)-\epsilon_0  \me^{-\gamma_0Lt/2} \vspace{5pt}\\
v_L(t,x) \leq \phi_L(L(x-c_Lt-C_+)-K_0\epsilon_0,x)+\epsilon_0  \me^{-\gamma_0Lt/2}
\end{array}\right. \hbox{ for all }t\geq 0\hbox{ and }x\in\R.
\end{equation}
\end{lem}

\begin{proof}
We use a Fife-McLeod~\cite{fm} type sub- and super-solutions method to show this lemma. We only give the construction of a super-solution, as the analysis for a sub-solution is analogous.

First of all, since the initial function $g$ satisfies~\eqref{i-condition-3}, it follows from Theorem~\ref{wave-length-1}~(i) and~(iii) that there exists a constant $C_+>0$ independent of $L$ such that, for all $L\ge L_*$ and $x\in\R$,
\be\label{phiLC+}
\phi_L\left(L(x-C_+),x\right)+\epsilon_0 =  \phi_L\Big(L\Big(x-C_+-\frac{\zeta_L(x)}{L}\Big)+\zeta_L(x),x\Big) +\epsilon_0\geq g(x).
\ee

For any $L\ge L_*$, we then set
$$\overline{v}_L(t,x):= \phi_{L}\big(L(x-c_Lt-C_+)+\eta_L(t),x\big)+q_L(t) \,\,\hbox{ for }t\geq 0\hbox{ and }x\in\R,$$
where $t\mapsto q_L(t)$ and $t\mapsto\eta_L(t)$ are some $C^1([0,+\infty))$ functions satisfying
\begin{equation}\label{condition-eta-q}\left\{\baa{lll}
q_L(0)=\epsilon_0, & q_L'(t) < 0<q_L(t) & \hbox{for all }t\ge0,\vspace{3pt}\\
\eta_L(0)=0, & \eta_L'(t)<0 & \hbox{for all }t\ge0.\eaa\right.
\end{equation}
By choosing suitable functions $\eta_L$ and $q_L$ later, we will show that $\overline{v}_L$ is a super-solution of~\eqref{eq-scale}.

Notice first that~\eqref{phiLC+}-\eqref{condition-eta-q} imply that $\overline{v}_L(0,\cdot)\ge g=v_L(0,\cdot)$ in $\R$ for each $L\ge L_*$. Now, for $(t,x)\in (0,+\infty)\times\R$, we define
$$N_L(t,x):=\partial_t\overline{v}_L(t,x)- \frac{1}{L}\partial_x(a(x)\partial_x\overline{v}_L(t,x))-Lf(x,\overline{v}_L(t,x)). $$
Since $(t,x)\mapsto \phi_L(L(x-c_Lt-C_+),x)$ is an entire solution of the first equation of~\eqref{eq-scale}, a straightforward calculation gives, for all $(t,x)\in (0,+\infty)\times\R$,
$$\baa{rcl}
N_L(t,x) &\!\!\! =\!\!\! & \eta_L'(t)\partial_{\xi}\phi_L(L(x-c_Lt-C_+)+\eta_L(t),x)+q_L'(t) \vspace{3pt}\\
& & \quad + Lf\big(x, \phi_{L}(L(x-c_Lt-C_+)+\eta_L(t),x)\big) -Lf(x,\overline{v}_L(t,x)).\eaa$$
Since $\partial_{\xi}\phi_L <0$ in $\R^2$, we get $\eta_L'(t)\partial_{\xi}\phi_L(L(x-c_Lt-C_+)+\eta_L(t),x)>0$ and $0<q_L(t)\leq \epsilon_0$ for all $t\geq 0$ and $x\in\R$, provided $q_L$ and $\eta_L$ fulfill~\eqref{condition-eta-q}. Then, on the one hand, for any pair $(t,x) \in (0,+\infty)\times\R$ such that $ \phi_{L}(L(x-c_Lt-C_+)+\eta_L(t),x) \in (0,\epsilon_0]\times [1-\epsilon_0,1)$, it follows from~\eqref{p-derivative-u} that
$$N_L(t,x)\geq q_L'(t) +Lf(x, \phi_{L}(L(x-c_Lt-C_+)+\eta_L(t),x))-Lf(x,\overline{v}_L(t,x))\geq q_L'(t)+\frac{\gamma_0L}{2}q_L(t).$$
On the other hand, by Theorem~\ref{wave-length-1}~(ii), there is $\beta>0$ (independent of $L\ge L_*$) such that $\partial_{\xi} \phi_{L}(\xi,y)\le-\beta$ for all $L\ge L_*$ and for all $(\xi,y)\in\R^2$ with $\epsilon_0 \leq \phi_L(\xi,y) \leq 1-\epsilon_0$. It follows that, if $ \phi_{L}(L(x-c_Lt-C_+)+\eta_L(t),x) \in (\epsilon_0,1-\epsilon_0)$, then
\begin{equation*}
\left.\begin{array}{ll}
N_L(t,x) \!\!\!&\geq  -\eta_L'(t)\beta+ q_L'(t) +Lf(x, \phi_{L}(L(x-c_Lt-C_+)+\eta_L(t),x))-Lf(x,\overline{v}_L(t,x))\vspace{5pt}\\
&\geq  -\eta_L'(t)\beta+ q_L'(t)-LC_1q_L(t),
\end{array}\right.
\end{equation*}
where $C_1:=\|\partial_uf\|_{L^\infty(\R\times\R)}$. Let us now choose $q_L(t)$ and $\eta_L(t)$ such that
$$\left\{\baa{lll}
q_L(0)=\epsilon_0, & \displaystyle q_L'(t)+\frac{\gamma_0L}{2}q_L(t)=0 & \hbox{for }t\ge0,\vspace{3pt}\\
\eta_L(0)=0, & -\beta\eta_L'(t)+ q_L'(t)-LC_1q_L(t)=0 & \hbox{for }t\ge0.\eaa\right.$$
Namely, we set
$$q_L(t)= \epsilon_0 \me^{-\gamma_0Lt/2}\ \hbox{ and }\ \eta_L(t)= -\frac{\epsilon_0(\gamma_0+2C_1)}{\beta\gamma_0}(1- \me^{-\gamma_0Lt/2})\,\,\hbox{ for }t\geq 0.$$
These functions $q_L$ and $\eta_L$ satisfy~\eqref{condition-eta-q}. Consequently, $N_L(t,x)\geq 0$ for all $(t,x)\in (0,+\infty)\times\R$.

Finally, the comparison principle implies that, for any $L\geq L_*$, $\overline{v}_L(t,x)\geq v_L(t,x)$ for all $t\geq 0$ and $x\in\R$. Taking   $K_0=(\gamma_0+2C_1)/(\beta\gamma_0)$ and using the fact that $\phi_L(\xi,y)$ is decreasing in $\xi$, we obtain the second inequality of~\eqref{ineq-sub-super}. As we have mentioned above, similar arguments imply the first one. The proof of Lemma~\ref{sub-super} is thus compete.
\end{proof}

We are now ready to complete the

\begin{proof}[{\it Proof of Theorem~$\ref{speed-limit}$}]
Assume by contradiction that $c_L$ does not converge to $c_*$ as $L\to+\infty$. Since the family $(c_L)_{L\geq L_*}\subset (0,+\infty)$ is bounded by Theorem~\ref{th-known}~(i), one finds a sequence $(L_n)_{n\in\N} \subset [L_*,+\infty)$ with $L_n\to+\infty$ as $n\to+\infty$ and a real number $c_{\infty}\geq 0$ such that
$$ c_{L_n} \to c_{\infty} \,\hbox{ as }n\to+\infty,\quad  \hbox{and}\quad  c_{\infty}\neq c_*.  $$
Without loss of generality, we assume that $c_{\infty}<c_*$ (as sketched below, the case where $c_{\infty} >c_*$ can be treated similarly). Then we have $c_{L_n}<(c_\infty+c_*)/2<c_*$ for all large $n\in\N$. Now, we can choose a large time $t_*>0$ and a positive integer $x_*$ (both independently of $n$) such that
\begin{equation}\label{choose-t*}
c_{L_n}t_*+C_++1 < x_* < c_*t_* \,\,\hbox{ for all large }n\in\N,
\end{equation}
where $C_+\geq 0$ is the constant given by Lemma~\ref{sub-super}. Since each function $\phi_{L_n}(\xi,y)$ is decreasing in $\xi$, it follows from the notations of Lemma~\ref{sub-super} and the second inequality of~\eqref{ineq-sub-super} that
\begin{equation}\label{vLn}\left.\begin{array}{ll}
0< v_{L_n}(t_*,x_*) \!\!\!& \leq \phi_{L_n}({L_n}(x_*-c_{L_n}t_*-C_+)-K_0\epsilon_0,x_*)+\epsilon_0  \me^{-\gamma_0L_nt_*/2}\vspace{5pt}\\
& \leq \phi_{L_n}({L_n}-K_0\epsilon_0,x_*)+\epsilon_0  \me^{-\gamma_0L_nt_*/2}
\end{array}\right.
\end{equation}
for all large $n\in\N$. Moreover, since $x_*\in\N$ and each function $\zeta_{L_n}$ given by~\eqref{normmal-front-1} is $1$-periodic, we have $\zeta_{L_n}(x_*)=\zeta_{L_n}(0)=0$. Therefore, passing to the limit as $n\to+\infty$ in~\eqref{vLn}, we see from Theorem~\ref{wave-length-1} that $\phi_{L_n}({L_n}-K_0\epsilon_0,x_*) \to 0$ as $n\to+\infty$, hence $v_{L_n}(t_*,x_*) \to 0$ as~$n\to+\infty$.

On the other hand, we had chosen $g$ such that $\Gamma_0=\{x\in\R:g(x)=b(x)\}=\{0\}$. Since~$x_*$ is a positive integer and since the function $y\mapsto c(y)$ is positive and $1$-periodic, we have $\rho(x_*,\Gamma_0)=x_*/c_*$, whence $\rho(x_*,\Gamma_0) <t_*$ due to~\eqref{choose-t*}. It then follows from Lemma~\ref{singular-limit} that  $v_{L_n}(t_*,x_*) \to 1$ as  $n\to+\infty$, yielding a contradiction. Therefore, the case where $c_{\infty}<c_*$ is ruled out.

In the case where $c_{\infty}>c_*$, one can derive a similar contradiction by choosing a large time $t^*>0$ and a positive integer $x^*$ (both independently of $n$) such that
$$c_{L_n}t^*-C_--1 >x^* > c_*t^* \,\,\hbox{ for all large }n\in\N.$$

As a conclusion, $c_{L} \to c_*$ as $L\to+\infty$, and the proof of Theorem~\ref{speed-limit} is thus complete.
\end{proof}


\section{Convergence of the front profiles: proof of Theorem~\ref{front-limit}}\label{sc-front-limit}\label{sec4}

This section is devoted to the proof of Theorem~\ref{front-limit}, that is, the convergence of the front profiles~$\phi_L$ as $L\to+\infty$, under the assumptions (A1)-(A4). The key step is to determine the asymptotic behavior as $L\to+\infty$ of the solutions of Cauchy problem of the following equation:
\begin{equation}\label{Cauchy-y}
\partial_t z_L= a\partial_{xx}z_L+ f\left(y+\frac{x}{L},z_L\right) \,\, \hbox{ for }t>0,\ x\in\R,
\end{equation}
where $y\in \R$ is arbitrary. We point out that $a$ is here a positive constant, thanks to~(A4). We recall that, for each $y\in\R$, the couple $(\psi(\cdot,y), c(y))$ denotes the front profile and front speed of the homogeneous equation~\eqref{traveling-wave}-\eqref{normal-homo-wave}. We first present some properties on $(\psi(\cdot,y), c(y))$ in Section~\ref{sc-pre}. In Section~\ref{sc-front-limit-1}, we choose the initial condition of~\eqref{Cauchy-y} sufficiently close to the homogeneous front $\psi(\cdot,y)$, and then construct a pair of sub- and super-solutions of~\eqref{Cauchy-y}. This will ensure that, at a certain time $t$, $z_L(t,\cdot+L)$ is still close to $\psi(\cdot,y)$ provided that $L$ is sufficiently large (see Proposition~\ref{compare-v-psi} in Section~\ref{sc-front-limit-2}). We complete the proof of Theorem~\ref{front-limit} in Section~\ref{sc-front-limit-4}.


\subsection{Preliminaries on homogeneous fronts}\label{sc-pre}

In this subsection, we present some properties regarding the homogeneous traveling front $(\psi(\cdot,y),c(y))$ of~\eqref{traveling-wave}-\eqref{normal-homo-wave} with respect to the parameter $y$. These properties will be used in our construction of sub- and super-solutions later and are also of interest in themselves.

\begin{pro}\label{homo-property}
Let {\rm (A1)-(A3)} hold. Then, the following statements hold true:
\begin{itemize}
\item [{\rm (i)}] the function $y\mapsto c(y)$ is positive, $1$-periodic, and the function $(\xi,y)\mapsto\psi(\xi,y)\in(0,1)$ is $1$-periodic in $y\in\R$ and decreasing in $\xi\in\R$;
\item [{\rm (ii)}]  $\lim_{\xi\to +\infty} \psi(\xi,y)=0$ and $\lim_{\xi\to -\infty} \psi(\xi,y)=1$ uniformly in $y\in\R$; furthermore, there exist positive constants $\mu_1$, $\mu_2$, $M$, $C_1$ and $C_2$ $($all are independent of $y$$)$ such that
\begin{equation}\label{psi-exp}
\left\{\baa{lll}
0< \psi(\xi,y) \leq C_1 \me^{-\mu_1 \xi}  & \hbox{for all  }\xi\geq M & \!\!\!\hbox{and }y\in\R,\vspace{5pt}\\
0< 1-\psi(\xi,y) \leq C_2 \me^{\mu_2 \xi}  & \hbox{for all  }\xi\leq -M & \!\!\!\hbox{and }y\in\R;
\eaa\right.
\end{equation}
\item [{\rm (iii)}] for any $\delta\in (0,1/2]$, there exists $\gamma=\gamma(\delta)>0$ $($independent of $y$$)$ such that
\begin{equation*}
\partial_{\xi} \psi(\xi,y) \leq -\gamma \,\hbox{ for all }(\xi,y)\in \R^2\hbox{ such that }\delta \leq \psi(\xi,y) \leq 1-\delta;
\end{equation*}
\item [{\rm (iv)}] the function $y\mapsto c(y)$ is of class $C^1(\R)$, and the function $(\xi,y) \mapsto \psi(\xi,y)$ is of class~$C^{2;1}_{\xi;y}(\R^2)$, and satisfies
\begin{equation}\label{esti-py-psi}
\sup_{\xi\in\R,\,y\in\R} \left| \partial_y \psi(\xi,y)\right| <+\infty.
\end{equation}
\end{itemize}
\end{pro}

The proof of Proposition~\ref{homo-property} is quite lengthy and is therefore postponed to the appendix in Section~\ref{sec5}. Let us emphasize that this proposition holds without the extra assumption (A4).


\subsection{Sub- and super-solutions}\label{sc-front-limit-1}

In the remaining part of Section~\ref{sec4}, we always assume that (A1)-(A4) hold. To present our sub- and super-solutions, we need some notations. Since the function $z\mapsto c(z)$ is $1$-periodic and of class $C^1(\R)$ by Proposition~\ref{homo-property}~(iv), one infers that, for each $y\in\R$ and $L>0$, the following ODE problem
\begin{equation}\label{ode-speed}
X'_{y,L}(t)=c\left(y+\frac{X_{y,L}(t)}{L}\right) \,\hbox{ for }t\ge0,\quad X_{y,L}(0)=0,
\end{equation}
admits a unique solution $X_{y,L}:[0,+\infty)\to\R$, with $X_{y,L}=X_{y+1,L}$ in $[0,+\infty)$ for each $y\in\R$ and~$L>0$. Furthermore, each function $X_{y,L}$ is increasing in $[0,+\infty)$ and $X_{y,L}(t)\to+\infty$ as $t\to+\infty$, since $\min_\R c>0$. Let then $T_L>0$ be the unique time such that
$$X_{y,L}(T_L)=L.$$
Since the function~$c$ is $1$-periodic, it is then checked by integrating the function $t\mapsto X_{y,L}'(t)/c(y+X_{y,L}(t)/L)=1$ over $[0,T_L]$ that $T_L$ does not depend on $y$ and is equal to
\begin{equation}\label{define-T}
T_L=\frac{L}{c_*}
\end{equation}
where $c_*>0$ is the constant defined by~\eqref{limit-speed}. In particular, $T_L\to+\infty$ as $L\to+\infty$ (uniformly in $y\in\R$). Notice also that, again by $1$-periodicity of the function $c$, one has
\be\label{XkTL}
X_{y,L}(kT_L)=kL\ \hbox{ for all }k\in\N.
\ee

Our super-solution is stated in the following lemma.

\begin{lem}\label{super-solu}
There exists $\epsilon_0\in(0,\delta_0)$, with $\delta_0\in(0,1/2)$ as in assumption~{\rm{(A2)}}, such that, for every $\epsilon\in(0,\epsilon_0]$, there exists $L_{1,\epsilon}>0$ such that, for every $y\in\R$ and $L\geq L_{1,\epsilon}$, the function $v^+_{\epsilon,y,L}:[0,+\infty)\times\R\to\R$ defined by
$$v^+_{\epsilon,y,L}(t,x):= \psi\left(x-X_{y,L}(t)+\eta_{\epsilon,L}(t),\,y+\frac{X_{y,L}(t)}{L}\right)+q_{\epsilon,L}(t) $$
is a super-solution of~\eqref{Cauchy-y} for $t\ge0$ and $x\in\R$, where $q_{\epsilon,L}$ and $\eta_{\epsilon,L}$ are $C^1([0,+\infty))$ functions independent of $y\in\R$ satisfying
\begin{equation}\label{condition-q}\left\{\baa{lll}
q_{\epsilon,L}(0)= \epsilon, & q_{\epsilon,L}'(t)<0<q_{\epsilon,L}(t) & \hbox{for all } t\ge0,\vspace{5pt}\\
\eta_{\epsilon,L}(0)=0, & \eta_{\epsilon,L}'(t) < 0 & \hbox{for all } t\ge0.\eaa\right.
\end{equation}
\end{lem}

\begin{proof}
First of all, by the assumptions (A1)-(A2) and the $C^1(\R\times[0,1])$ smoothness of $f$ and its periodicity in $x$, together with~\eqref{extf}, there exist $\delta_1\in(0,1/2)$ and $\gamma_1>0$ such that
\begin{equation}\label{choose-delta1}
\partial_uf(x,u) \leq -\gamma_1 \,\,\hbox{ for all }x\in\R\hbox{ and }u\in [-\delta_1,\delta_1]\cap [1-\delta_1,1+\delta_1].
\end{equation}
Without loss of generality, one can assume that $\delta_1\le\delta'_0$, with $\delta'_0>0$ as in assumption~(A4). Moreover, by Proposition~\ref{homo-property} (iv) and the boundedness of the function $z\mapsto c(z)$, there is a constant $C_1>0$ such that, for all $y\in\R$ and $L>0$,
\be\label{estX}
\left| \frac{X_{y,L}'(t)}{L}\partial_y \psi(\xi,z) \right| \leq \frac{C_1}{L}\,\,\hbox{ for all $t\ge0$ and $(\xi,z)\in\R\times\R$}.
\ee
It further follows from Proposition~\ref{homo-property}~(i)-(ii) that there exist $M_1>0$ and
$$\epsilon_0\in (0,\min\{\delta_1/2,\delta_0\})\subset(0,1/4),$$
with $\delta_0\in(0,1/2)$ as in assumption~(A2), such that, for all $(\xi,z)\in\R\times\R$,
 \begin{equation*}
\left\{\baa{ll}
\displaystyle 0< \psi(\xi,z)\leq\frac{\delta_1}{2} & \hbox{ if }\xi\geq M_1,\vspace{6pt}\\
\displaystyle 1-\frac{\delta_1}{2}\leq \psi(\xi,z)<1 & \hbox{ if }\xi\leq -M_1,\vspace{6pt}\\
 \displaystyle 2\epsilon_0\leq \psi(\xi,z)\le1-2\epsilon_0 & \hbox{ if }-M_1<\xi< M_1,
\eaa \right.
\end{equation*}
where the last inequality is actually a consequence of the continuity of $\psi:\R^2\to(0,1)$ and its periodicity in the second variable.

In the arguments below, $\epsilon\in(0,\epsilon_0]$ is arbitrary. To show Lemma~\ref{super-solu}, it suffices to find suitable $C^1([0,+\infty))$ functions $q_{\epsilon,L}$ and $\eta_{\epsilon,L}$ satisfying~\eqref{condition-q} such that, if $L$ is sufficiently large and independently of~$y$, there holds
\begin{equation*}
\displaystyle N(t,x):= \partial_t v^+_{\epsilon,y,L}(t,x)- a\partial_{xx}v^+_{\epsilon,y,L}(t,x)- f\left(y+\frac{x}{L},v^+_{\epsilon,y,L}(t,x)\right)\geq 0
\end{equation*}
for all $(t,x)\in[0,+\infty)\times\R$. Since $\psi$ is a solution of~\eqref{traveling-wave} and since here $a$ is constant by assumption~(A4), it is straightforward to check that, for any $y\in\R$, $L>0$ and any $C^1([0,+\infty))$ functions $q_{\epsilon,L}$ and $\eta_{\epsilon,L}$, one has
\begin{equation*}
N(t,x)=\frac{X_{y,L}'(t)}{L}\partial_y \psi +\eta_{\epsilon,L}'(t)\partial_{\xi}\psi+q_{\epsilon,L}'(t)+f\left(y+\frac{X_{y,L}(t)}{L},\psi\right)-f\left(y+\frac{x}{L},v^+_{\epsilon,y,L}(t,x)\right)
\end{equation*}
for all $(t,x)\in[0,+\infty)\times\R$, where $\psi=\psi(x-X_{y,L}(t)+\eta_{\epsilon,L}(t),y+X_{y,L}(t)/L)$, and $\partial_{\xi}\psi$, $\partial_y\psi$ stand for the partial derivatives of $\psi$ with respect to the first variable and the second variable, respectively, evaluated at the same point $(x-X_{y,L}(t)+\eta_{\epsilon,L}(t),y+X_{y,L}(t)/L)$. We complete the proof by considering three cases: (a) $x-X_{y,L}(t)+\eta_{\epsilon,L}(t)\geq M_1$, (b)~$x-X_{y,L}(t)+\eta_{\epsilon,L}(t)\leq -M_1$, (c)~$-M_1<x-X_{y,L}(t)+\eta_{\epsilon,L}(t)< M_1$.

In case~(a), with $q_{\epsilon,L}$ required to satisfy~\eqref{condition-q}, we have
$$0<\psi\Big(x-X_{y,L}(t)+\eta_{\epsilon,L}(t),y+\frac{X_{y,L}(t)}{L}\Big)< v^+_{\epsilon,y,L}(t,x)\le\frac{\delta_1}{2}+\epsilon\le\frac{\delta_1}{2}+\epsilon_0\le\delta_1.$$
Since $\psi$ is decreasing in its first variable and since $\eta_{\epsilon,L}$ is required to satisfy~\eqref{condition-q}, we have
$$\eta_{\epsilon,L}'(t)\,\partial_{\xi}\psi\Big(x-X_{y,L}(t)+\eta_{\epsilon,L}(t),y+\frac{X_{y,L}(t)}{L}\Big)\ge0.$$
Remember that $f(x,u)$ is independent of $x$ for $u\in[0,\delta_1]\subset[0,\delta'_0]$ by assumption~(A4).  It then follows from~\eqref{choose-delta1}-\eqref{estX} that
$$\baa{rcl}
N(t,x) & \!\!\!\geq\!\!\! & \displaystyle -\frac{C_1}{L}+q_{\epsilon,L}'(t)+f\left(0,\psi\left(x-X_{y,L}(t)+\eta_{\epsilon,L}(t),y+\frac{X_{y,L}(t)}{L}\right)\right)-f(0,v^+_{\epsilon,y,L}(t,x))\vspace{3pt}\\
& \!\!\!\geq\!\!\! & \displaystyle-\frac{C_1}{L}+q_{\epsilon,L}'(t)+\gamma_1q_{\epsilon,L}(t).\eaa$$
Let us choose
\be\label{defL1eps}
L_{1,\epsilon}=\frac{2C_1}{\gamma_1\epsilon}>0
\ee
(notice that $L_{1,\epsilon}$ is independent of $y\in\R$) and the function $q_{\epsilon,L}$ such that
\be\label{defqepsL}
q_{\epsilon,L}(0)=\epsilon\ \hbox{ and }\ -\frac{C_1}{L}+q_{\epsilon,L}'(t)+\gamma_1q_{\epsilon,L}(t)=0\hbox{ for }t\ge0,
\ee
namely,
\begin{equation}\label{function-q}
q_{\epsilon,L}(t)= \frac{C_1}{L\gamma_1}+\left(\epsilon-\frac{C_1}{L\gamma_1}\right)\me^{-\gamma_1t} \,\,\hbox{ for }t\geq 0.
\end{equation}
It is clear that, for any $L\geq L_{1,\epsilon}$ and for any $y\in\R$, the function $q_{\epsilon,L}$ satisfies~\eqref{condition-q} and $N(t,x)\geq 0$ for all $(t,x)\in[0,+\infty)\times\R$ such that $x-X_{y,L}(t)+\eta_{\epsilon,L}(t)\geq M_1$, provided $\eta_{\epsilon,L}$ satisfies~\eqref{condition-q} too.

Proceeding similarly as above, we can conclude that, for any $L\geq L_{1,\epsilon}$ and $y\in\R$, $N(t,x)\geq 0$ for all $(t,x)\in[0,+\infty)\times\R$ such that $x-X_{y,L}(t)+\eta_{\epsilon,L}(t)\leq-M_1$, as soon as $\eta_{\epsilon,L}$ satisfies~\eqref{condition-q}.

It remains to find a suitable function $\eta_{\epsilon,L}$ such that $N(t,x)\geq 0$ in case~(c), i.e., $-M_1<x-X_{y,L}(t)+\eta_{\epsilon,L}(t)< M_1$. In this case, we have
$$2\epsilon_0\leq \psi\left(x-X_{y,L}(t)+\eta_{\epsilon,L}(t),y+\frac{X_{y,L}(t)}{L}\right) \leq 1-2\epsilon_0$$
and $2\epsilon_0 \leq v^+_{\epsilon,y,L}(t,x)\leq 1-\epsilon_0$ (since $0<q_{\epsilon,L}(t)\leq \epsilon \leq \epsilon_0$). It then follows from Proposition~\ref{homo-property}~(iii) that there exists $\beta_1>0$ (independent of $\epsilon$, $y$, $L$, $t$ and $x$) such that $\partial_{\xi}\psi(x-X_{y,L}(t)+\eta_{\epsilon,L}(t),y+X_{y,L}(t)/L) \leq -\beta_1$. Noticing that the function $f$ is of class $C^1(\R\times[0,1])$ and periodic in $x$, one finds some constants $K_1>0$ and $K_2>0$ (independent of $\epsilon$, $y$, $L$, $t$ and $x$) such that
$$\left\{\baa{ll}
|f(x_1,u)-f(x_2,u)| \leq K_1|x_1-x_2| & \hbox{ for all }x_1,\,x_2\in \R\hbox{ and }u\in [0,1],\vspace{3pt}\\
|f(z,u_1)-f(z,u_2)|\le K_2|u_1-u_2|& \hbox{ for all }z\in\R\hbox{ and }u_1,\,u_2\in[0,1].\eaa\right.$$
Therefore, with $\eta_{\epsilon,L}$ required to satisfy~\eqref{condition-q} (then being nonpositive in $[0,+\infty)$), we have
\begin{equation*}
\left.\baa{ll}
\displaystyle N(t,x) \!\!\!& \displaystyle \geq -\frac{C_1}{L}-\beta_1 \eta_{\epsilon,L}'(t)+q_{\epsilon,L}'(t)+ f\left(y+\frac{X_{y,L}(t)}{L},\psi\Big(x-X_{y,L}(t)+\eta_{\epsilon,L}(t),y+\frac{X_{y,L}(t)}{L}\Big)\right)\vspace{6pt}\\
& \ \ \ \ \displaystyle-f\left(y+\frac{x}{L},v^+_{\epsilon,y,L}(t,x)\right) \vspace{6pt}\\
& \geq \displaystyle -\frac{C_1}{L}-\beta_1 \eta_{\epsilon,L}'(t)+q_{\epsilon,L}'(t)-K_1\frac{|X_{y,L}(t)-x|}{L}-K_2q_{\epsilon,L}(t),\vspace{6pt}\\
& \geq \displaystyle-\frac{C_1}{L}-\beta_1 \eta_{\epsilon,L}'(t)+q_{\epsilon,L}'(t)-K_1\frac{M_1-\eta_{\epsilon,L}(t)}{L}-K_2q_{\epsilon,L}(t).\eaa \right.
\end{equation*}
Then, with $q_{\epsilon,L}(t)$ given by~\eqref{defqepsL}-\eqref{function-q}, by choosing $\eta_{\epsilon,L}$ such that
$$\eta_{\epsilon,L}(0)=0\quad\hbox{and}\quad  -\beta_1 \eta_{\epsilon,L}'(t) +\frac{K_1}{L} \eta_{\epsilon,L}(t)-(\gamma_1+K_2)q_{\epsilon,L}(t)=\frac{K_1M_1}{L}  \,\,\hbox{ for }t\ge0,$$
we have $N(t,x)\geq 0$ for $(t,x)\in[0,+\infty)\times\R$ such that $-M_1<x-X_{y,L}(t)+\eta_{\epsilon,L}(t)< M_1$. It is straightforward to compute that
\begin{equation}\label{function-eta}
\eta_{\epsilon,L}(t)=-\frac{\gamma_1+K_2}{\beta_1}\left[(C_2+C_{3,\epsilon,L}) \me^{K_1t/(L\beta_1)} -  C_{3,\epsilon,L} \me^{-\gamma_1 t} -C_2 \right] \,\,\hbox{ for }t\geq 0,
\end{equation}
where
$$C_2= \frac{C_1\beta_1}{\gamma_1K_1}+\frac{M_1\beta_1}{K_2+\gamma_1} \quad\hbox{and}\quad C_{3,\epsilon,L}=\frac{\epsilon-C_1/(L\gamma_1)}{\gamma_1+K_1/(L\beta_1)}.$$
It is clear that $C_2>0$ and $C_{3,\epsilon,L}>0$ if $L\geq L_{1,\epsilon}$, whence the function $\eta_{\epsilon,L}$ given in~\eqref{function-eta} satisfies~\eqref{condition-q}.

Combining the above, we can conclude that for any $\epsilon\in(0,\epsilon_0]$, $y\in\R$ and $L\geq L_{1,\epsilon}$, one has $N(t,x)\geq 0$ for all $(t,x)\in[0,+\infty)\times\R$. This ends the proof of Lemma~\ref{super-solu}.
\end{proof}

The following lemma gives the sub-solution of problem~\eqref{Cauchy-y}.

\begin{lem}\label{sub-solu}
Let $\epsilon_0\in(0,\delta_0)$ be given by Lemma~$\ref{super-solu}$ and, for any $\epsilon\in(0,\epsilon_0]$,
let $L_{1,\epsilon}>0$ be given by~\eqref{defL1eps} as in Lemma~$\ref{super-solu}$. Then, for every $\epsilon\in(0,\epsilon_0]$, $y\in\R$ and $L\geq L_{1,\epsilon}$, the function $v^-_{\epsilon,y,L}:[0,+\infty)\times\R\to\R$ defined by
$$v^-_{\epsilon,y,L}(t,x):= \psi\left(x-X_{y,L}(t)-\eta_{\epsilon,L}(t),\,y+\frac{X_{y,L}(t)}{L}\right)-q_{\epsilon,L}(t) $$
is a sub-solution of~\eqref{Cauchy-y} for $t\ge0$ and $x\in\R$, where $q_{\epsilon,L}$ and $\eta_{\epsilon,L}$ are $C^1([0,+\infty))$ functions satisfying~\eqref{condition-q}, given by~\eqref{function-q} and~\eqref{function-eta}, respectively.
\end{lem}

\begin{proof}
The proof is analogous to that of Lemma~\ref{super-solu}; therefore, we omit the details.
\end{proof}

Before going further on, we give two remarks about the functions $q_{\epsilon,L}$ and $\eta_{\epsilon,L}$, which will be useful in the proof of our main result later.

\begin{rem}\label{rm-bound-eta} 
For any $\epsilon\in(0,\epsilon_0]$ and $L\geq L_{1,\epsilon}=2C_1/(\gamma_1\epsilon)$, let $T_L=L/c_*>0$ be the time provided by~\eqref{define-T}. For any given $k\in \N$ and $\tau>0$, the function $\eta_{\epsilon,L}$ given by~\eqref{function-eta} is bounded in $[0,kT_L+\tau]$ uniformly with respect to $\epsilon\in(0,\epsilon_0]$ and $L\geq L_{1,\epsilon}$. More precisely, we can find some constants $0<A_1<A_2$ $($independent of $\epsilon$, $L$, $k$ and $\tau$$)$ such that
$$A_1- A_2\me^{K_1k/(c_*\beta_1)+K_1\tau\gamma_1\epsilon_0/(2C_1\beta_1)}\le A_1- A_2\me^{K_1k/(c_*\beta_1)+K_1\tau/(L\beta_1)}  \leq \eta_{\epsilon,L}(t) \leq 0$$
for all $0\leq t \leq kT_L+\tau$.
\end{rem}

\begin{rem}\label{limit-eta-q}
We also point out that, for any $\epsilon\in(0,\epsilon_0]$, the functions $t\mapsto q_{\epsilon,L}(t)$ and $t\mapsto\eta_{\epsilon,L}(t)$ converge locally uniformly in $t\geq 0$ as $L\to+\infty$. More precisely, there holds
$$\lim_{L\to+\infty} q_{\epsilon,L}(t) = \epsilon\,\me^{-\gamma_1t} \quad\hbox{and}\quad \lim_{L\to+\infty} \eta_{\epsilon,L}(t)=-\epsilon\,\frac{\gamma_1+K_1}{\gamma_1\beta_1}\,(1-\me^{-\gamma_1t}),$$
locally uniformly in $t\in[0,+\infty)$.
\end{rem}

Now, for any $\epsilon\in(0,\epsilon_0]$, with $\epsilon_0\in(0,\delta_0)$ provided by Lemmas~\ref{super-solu}-\ref{sub-solu}, we consider any family of continuous functions $v^0_{\epsilon,y}:\R\to[0,1]$ such that
\begin{equation}\label{i-condition-1}
v^0_{\epsilon,y+1}\equiv v^0_{\epsilon,y}\hbox{ in }\R,\ \ \displaystyle v^0_{\epsilon,y}(0)=\frac{1}{2}\ \hbox{ and }\|v^0_{\epsilon,y}- \psi(\cdot,y)\|_{L^{\infty}(\R)} \leq \epsilon\ \hbox{ for all }y\in\R.\footnote{A typical example is given by: $v^0_{\epsilon,y}=\psi(\cdot,y)$ for $y\in\R$, but other functions $v^0_{\epsilon,y}$ will be used in the proof of Theorem~\ref{front-limit} in Section~\ref{sc-front-limit-4}.}
\end{equation}
Then, for any $y\in\R$ and $L>0$, we denote by $v_{\epsilon,y,L}:[0,+\infty)\times\R\to[0,1]$ the solution of~\eqref{Cauchy-y} with the initial condition~$v^0_{\epsilon,y}$. Since the function $f$ is $1$-periodic in its first variable, it follows that the functions $v_{\epsilon,y,L}$ are $1$-periodic with respect to the parameter $y\in\R$. It is also easily seen that
\be\label{ineqvpm0}
\max\big\{v^-_{\epsilon,y,L}(0,x),0\big\}\leq  v^0_{\epsilon,y}(x) \leq\min\big\{v^+_{\epsilon,y,L}(0,x),1\big\} \,\,\hbox{ for all }x\in\R,
\ee
with $v^{\pm}_{\epsilon,y,L}$ are in Lemmas~\ref{super-solu}-\ref{sub-solu}. Therefore, the next lemma is an immediate consequence of the comparison principle.

\begin{lem}\label{estimate-T-v}
For any $\epsilon\in(0,\epsilon_0]$, with $\epsilon_0\in(0,\delta_0)$ as in Lemmas~$\ref{super-solu}$-$\ref{sub-solu}$, any $y\in\R$, and any $L\geq L_{1,\epsilon}$ with $L_{1,\epsilon}>0$ as in~\eqref{defL1eps}, there holds
$$\max\big\{v^-_{\epsilon,y,L}(t,x),0\big\}\leq v_{\epsilon,y,L}(t,x) \leq\min\big\{v^+_{\epsilon,y,L}(t,x),1\big\}\ \hbox{ for all $t\geq 0$ and $x\in\R$}.$$
\end{lem}


\subsection{Asymptotic behavior of the Cauchy problem as $L\to+\infty$}\label{sc-front-limit-2}

In this subsection and the following one, we always let $\epsilon_0\in(0,\delta_0)$ be provided by Lemmas~\ref{super-solu}-\ref{sub-solu}. Now, for any $\epsilon\in(0,\epsilon_0]$, we consider the asymptotic behavior as $L\to+\infty$ of the solutions~$v_{\epsilon,y,L}$ of~\eqref{Cauchy-y} with initial conditions $v^0_{\epsilon,y}$ satisfying~\eqref{i-condition-1}, at a time $\tilde{T}_{\epsilon,y,L}$ defined by
\begin{equation}\label{define-tildeT}
\tilde{T}_{\epsilon,y,L}:=\inf\Big\{t>0: \, v_{\epsilon,y,L}(t,L)=\frac{1}{2}\Big\}.
\end{equation}
Notice immediately that, since the functions $v_{\epsilon,y,L}$ are $1$-periodic with respect to $y$, so are the quantities $\tilde{T}_{\epsilon,y,L}$'s. Furthermore, for each $y\in\R$ and $L\ge L_{1,\epsilon}$, one has $v_{\epsilon,y,L}(0,L)\le v^+_{\epsilon,y,L}(0,L)=\psi(L,y)+\epsilon$ by~\eqref{condition-q} and~\eqref{ineqvpm0}, hence
$$\limsup_{L\to+\infty}\Big(\sup_{y\in\R}v_{\epsilon,y,L}(0,L)\Big)\le\epsilon\le\epsilon_0<\delta_0<\frac12$$
since $\psi(+\infty,y)=0$ uniformly in $y\in\R$ by Proposition~\ref{homo-property}~(ii). Therefore, $\tilde{T}_{\epsilon,y,L}>0$ for all $L$ large enough, uniformly in $y\in\R$. On the other hand, for each $y\in\R$ and $L\ge L_{1,\epsilon}$, with $T_L$ as in~\eqref{define-T}, one has $v_{\epsilon,y,L}(2T_L,L)\ge v^-_{\epsilon,y,L}(2T_L,L)\ge\psi(-L-\eta_{\epsilon,L}(2T_L),y+2)-\epsilon$ by Lemma~\ref{estimate-T-v} together with~\eqref{XkTL}-\eqref{condition-q}. Hence,
$$\liminf_{L\to+\infty}\Big(\inf_{y\in\R}v_{\epsilon,y,L}(2T_L,L)\Big)\ge1-\epsilon\ge1-\epsilon_0>1-\delta_0>\frac12$$
by Remark~\ref{rm-bound-eta} and since $\psi(-\infty,y)=1$ uniformly in $y\in\R$ by Proposition~\ref{homo-property}~(ii). As a consequence, there is
$$L_{2,\epsilon}\ge L_{1,\epsilon}>0$$
such that
\be\label{TTtilde}
0<\tilde{T}_{\epsilon,y,L}\le 2T_L<+\infty\ \hbox{ for all $y\in\R$ and $L\ge L_{2,\epsilon}$},
\ee
and then
\be\label{v12}
v_{\epsilon,y,L}(\tilde{T}_{\epsilon,y,L},L)=\frac12.
\ee

The following proposition shows that when $L$ is sufficiently large, the profile $v_{\epsilon,y,L}(t,\cdot+L)$ at time $t=\tilde{T}_{\epsilon,y,L}$ is close to that of the initial condition $v_{\epsilon,y,L}(0,\cdot)=v^0_{\epsilon,y}$.

\begin{pro}\label{compare-v-psi}
For any $\epsilon\in(0,\epsilon_0]$, with $\epsilon_0\in(0,\delta_0)$ as in Lemmas~$\ref{super-solu}$-$\ref{sub-solu}$, there exists $L_{3,\epsilon}\geq L_{2,\epsilon}$, with $L_{2,\epsilon}>0$ as in~\eqref{TTtilde}, such that, for any $y\in\R$ and $L\geq L_{3,\epsilon}$, there holds
$$\left\|v_{\epsilon,y,L}(\tilde{T}_{\epsilon,y,L},\cdot+L)- \psi(\cdot,y)\right\|_{L^{\infty}(\R)} \leq\frac{\epsilon}{2},$$
where $\tilde{T}_{\epsilon,y,L}$ is given by~\eqref{define-tildeT}, and $v_{\epsilon,y,L}$ solves~\eqref{Cauchy-y} with initial conditions $v^0_{\epsilon,y}$ satisfying~\eqref{i-condition-1}.
\end{pro}

For the proof of Proposition~\ref{compare-v-psi}, let us first show two lemmas which are concerned with the comparison of $T_L$ and $\tilde{T}_{\epsilon,y,L}$.

\begin{lem}\label{compare-T-tT}
For any fixed $\epsilon\in(0,\epsilon_0]$, with $\epsilon_0\in(0,\delta_0)$ as in Lemmas~$\ref{super-solu}$-$\ref{sub-solu}$, let $T_L$ and $\tilde{T}_{\epsilon,y,L}$ be given by~\eqref{define-T} and~\eqref{define-tildeT}, respectively, for $y\in\R$ and $L\ge L_{2,\epsilon}$. Then, there is a constant $M_\epsilon\ge0$ such that
$$\sup_{y\in \R,\,L\ge L_{2,\epsilon}}\,|T_L-\tilde{T}_{\epsilon,y,L}|\le M_\epsilon,$$
for all families of initial conditions $v^0_{\epsilon,y}\in C(\R,[0,1])$ satisfying~\eqref{i-condition-1}.
\end{lem}

\begin{proof}
We fix $\epsilon\in(0,\epsilon_0]$ throughout the proof. We first prove the existence of some real numbers $L^+_{2,\epsilon}\in[L_{2,\epsilon},+\infty)$ and $M^+_{\epsilon}\ge0$ such that
\be\label{limsup1}
\sup_{y\in \R,\,L\ge L^+_{2,\epsilon}}\,(\tilde{T}_{\epsilon,y,L}-T_L)\le M^+_{\epsilon}
\ee
for all families of initial conditions $v^0_{\epsilon,y}\in C(\R,[0,1])$ satisfying~\eqref{i-condition-1}. Assume by contradiction, there exist sequences $(y_n)_{n\in\N}\subset \R$, $(L_n)_{n\in\N} \subset [L_{2,\epsilon},+\infty)$, and initial conditions $(v^0_{\epsilon,y_n})_{n\in\N}\subset C(\R,[0,1])$ satisfying~\eqref{i-condition-1}, such that $L_n\to+\infty$ and $\tilde{T}_{\epsilon,y_n,L_n}-T_{L_n}\to+\infty$ as $n\to+\infty$. Then for each $n\in\N$, by Lemma~\ref{estimate-T-v}, we have
\begin{equation}\label{bound-vLn}
\left\{\baa{l}
\displaystyle v_{\epsilon,y_n,L_n}(t,x)\geq\psi\Big(x-X_{y_n,L_n}(t)-\eta_{\epsilon,L_n}(t),\,y_n+\frac{X_{y_n,L_n}(t)}{L_n}\Big)-q_{\epsilon,L_n}(t)\vspace{6pt}\\
\displaystyle v_{\epsilon,y_n,L_n}(t,x)\leq\psi\Big(x-X_{y_n,L_n}(t)+\eta_{\epsilon,L_n}(t),\,y_n+\frac{X_{y_n,L_n}(t)}{L_n}\Big)+q_{\epsilon,L_n}(t)\eaa\right.
\ee
for all $t\geq 0$ and $x\in\R$.

We already know from~\eqref{TTtilde} that $\tilde{T}_{\epsilon,y_n,L_n}\le2T_{L_n}$ for all $n\in\N$. Next, for any $n\in\N$, choosing $t=\tilde{T}_{\epsilon,y_n,L_n}$ and $x=L_n$ in the first inequality of~\eqref{bound-vLn} yields
$$\baa{rcl}
\displaystyle\frac{1}{2}= v_{\epsilon,y_n,L_n}(\tilde{T}_{\epsilon,y_n,L_n},L_n)  & \!\!\!\geq\!\!\! & \displaystyle\psi\Big(L_n-X_{y_n,L_n}(\tilde{T}_{\epsilon,y_n,L_n})-\eta_{\epsilon,L_n}(\tilde{T}_{\epsilon,y_n,L_n}),\,y_n+ \frac{X_{y_n,L_n}(\tilde{T}_{\epsilon,y_n,L_n})}{L_n}\Big)\vspace{3pt}\\
& \!\!\!\!\!\! & -q_{\epsilon,L_n}(\tilde{T}_{\epsilon,y_n,L_n}).\eaa$$
Since we have assumed $\tilde{T}_{\epsilon,y_n,L_n}-T_{L_n}\to+\infty$ as $n\to+\infty$, we may assume without loss of generality that $\tilde{T}_{\epsilon,y_n,L_n}>T_{L_n}=L_n/c_*$ for all $n\in\N$ (hence, $\tilde{T}_{\epsilon,y_n,L_n}\to+\infty$ as $n\to+\infty$). Then, we have
$$X_{y_n,L_n}(\tilde{T}_{\epsilon,y_n,L_n})=\underbrace{X_{y_n,L_n}(T_{L_n})}_{=L_n}+\int_{T_{L_n}}^{\tilde{T}_{\epsilon,y_n,L_n}} c\left(y_n+ \frac{X_{y_n,L_n}(t)}{L_n}\right)dt \geq L_n + c^-(\tilde{T}_{\epsilon,y_n,L_n}-T_{L_n}),$$
where $c^-=\min_{z\in\R} c(z)>0$. Remember that the function $(\xi,z)\mapsto\psi(\xi,z)$ is decreasing in $\xi\in\R$. It follows that
\be\label{psiqeps}
\frac{1}{2} \geq   \psi\Big(-c^-(\tilde{T}_{\epsilon,y_n,L_n}-T_{L_n})-\eta_{\epsilon,L_n}(\tilde{T}_{\epsilon,y_n,L_n}),y_n+ \frac{X_{y_n,L_n}(\tilde{T}_{\epsilon,y_n,L_n})}{L_n} \Big)-q_{\epsilon,L_n}(\tilde{T}_{\epsilon,y_n,L_n}).
\ee
Since $\tilde{T}_{\epsilon,y_n,L_n}\to+\infty$ as $n\to+\infty$, we see from~\eqref{function-q} that $q_{\epsilon,L_n}(\tilde{T}_{\epsilon,y_n,L_n})\to 0$ as $n\to+\infty$. Furthermore, since $0<\tilde{T}_{\epsilon,y_n,L_n}\leq 2T_{L_n}$ by~\eqref{TTtilde}, one infers from Remark~\ref{rm-bound-eta} that the sequence $(\eta_{\epsilon,L_n}(\tilde{T}_{\epsilon,y_n,L_n}))_{n\in\N}$ is bounded. Passing to the limit as $n\to+\infty$ in~\eqref{psiqeps} and using Proposition~\ref{homo-property}~(ii), we obtain $1/2\geq 1$, which is impossible.
Thus,~\eqref{limsup1} is proved.

Similarly as above, one can prove the existence of some real numbers $L^-_{2,\epsilon}\in[L_{2,\epsilon},+\infty)$ and $M^-_{\epsilon}\ge0$ such that $\sup_{y\in \R,\,L\ge L^-_{2,\epsilon}}\,(T_L-\tilde{T}_{\epsilon,y,L})\le M^-_{\epsilon}$ for all families of initial conditions $v^0_{\epsilon,y}\in C(\R,[0,1])$ satisfying~\eqref{i-condition-1}, by using this time in the previous paragraph the second inequality of~\eqref{bound-vLn} instead of the first one, together with $q_{\epsilon,L_n}\le\epsilon<1/2$ in $[0,+\infty)$ by~\eqref{condition-q}. 

Finally, since $0<\tilde{T}_{\epsilon,y,L}\le 2T_L=2L/c_*$ for all $y\in\R$ and $L\ge L_{2,\epsilon}$ by~\eqref{define-T} and~\eqref{TTtilde}, the desired conclusion of Lemma~\ref{compare-T-tT} follows, with $M_\epsilon:=\max\{M^+_\epsilon,M^-_\epsilon,2L^+_{2,\epsilon}/c_*,2L^-_{2,\epsilon}/c_*\}$.
\end{proof}

\begin{lem}\label{property-X}
For any fixed $\epsilon\in(0,\epsilon_0]$, with $\epsilon_0\in(0,\delta_0)$ as in Lemmas~$\ref{super-solu}$-$\ref{sub-solu}$, let $T_L$ and $\tilde{T}_{\epsilon,y,L}$ be given by~\eqref{define-T} and~\eqref{define-tildeT}, respectively, for $y\in\R$ and $L\ge L_{2,\epsilon}$. Then there holds
\be\label{lim1}
\lim_{L\to+\infty}\frac{X_{y,L}(t+\tilde{T}_{\epsilon,y,L})}{X_{y,L}(T_L)}=1
\ee
and
\be\label{lim2}
\lim_{L\to+\infty} \,\big|X_{y,L}(t+\tilde{T}_{\epsilon,y,L})-X_{y,L}(T_L)-c(y)(t+\tilde{T}_{\epsilon,y,L}-T_L)\big|=0
\ee
locally uniformly with respect to $t\in\R$, and uniformly with respect to $y\in\R$ and with respect to the families of initial conditions $v^0_{\epsilon,y}\in C(\R,[0,1])$ satisfying~\eqref{i-condition-1}.
\end{lem}

\begin{proof}
For any $y\in\R$, $L\geq L_{2,\epsilon}$, and $t\in\R$, one has $X_{y,L}(T_L)=L$ and
$$\left|\frac{X_{y,L}(t+\tilde{T}_{\epsilon,y,L})}{X_{y,L}(T_L)}-1\right|= \frac{  \left|\displaystyle \int^{t+\tilde{T}_{\epsilon,y,L}}_{T_L} c\left(y+ \frac{X_{y,L}(s)}{L}\right)ds\right|}{X_{y,L}(T_L)} \leq  \frac{ c^+|\tilde{T}_{\epsilon,y,L}-T_L+t| }{L},$$
where $c^+=\max_{z\in\R} c(z)$. Lemma~\ref{compare-T-tT} then immediately gives~\eqref{lim1}.

It remains to show~\eqref{lim2}. Let $y\in \R$, $L\geq L_{2,\epsilon}$ and $t\in\R$ be arbitrary. Without loss of generality, we may assume that $t+\tilde{T}_{\epsilon,y,L}\geq T_L$ (the case $t+\tilde{T}_{\epsilon,y,L}\leq T_L$ can be treated identically).  Since
$$X_{y,L}(t+\tilde{T}_{\epsilon,y,L})-X_{y,L}(T_L)-c(y)(t+\tilde{T}_{\epsilon,y,L}-T_L) = \int_{T_L}^{t+\tilde{T}_{\epsilon,y,L}} \left(  \displaystyle c\Big(y+ \frac{X_{y,L}(s)}{L} \Big)- c(y) \right)  ds,  $$
and since the function $z\mapsto c(z)$ is $1$-periodic and of class $C^1$ by Proposition~\ref{homo-property}~(iv), we have $c(y)=c(y+1)$ and
$$\baa{l}\left|X_{y,L}(t+\tilde{T}_{\epsilon,y,L})-X_{y,L}(T_L)-c(y)(t+\tilde{T}_{\epsilon,y,L}-T_L)\right|\vspace{3pt}\\
\qquad\qquad\qquad\qquad\qquad\qquad\displaystyle\leq\|c'\|_{L^\infty(\R)}\times|t+\tilde{T}_{\epsilon,y,L}-T_L|\times\max_{s\in  [T_L,\,t+\tilde{T}_{\epsilon,y,L}] }\left| \frac{X_{y,L}(s)}{X_{y,L}(T_L)}-1\right|.\eaa$$
Lemma~\ref{compare-T-tT} and~\eqref{lim1} then yield~\eqref{lim2}.
\end{proof}

Now, we are ready to give the

\begin{proof}[Proof of Proposition~$\ref{compare-v-psi}$]
Fix any $\epsilon\in(0,\epsilon_0]$, and assume by contradiction that the conclusion is not true. Then, there exist sequences $(y_n)_{n\in\N}\subset \R$, $(L_n)_{n\in\N}\subset [L_{2,\epsilon},+\infty)$, $(x_n)_{n\in\N}\subset \R$, and initial conditions $(v^0_{\epsilon,y_n})_{n\in\N}\subset C(\R,[0,1])$ satisfying~\eqref{i-condition-1} such that $L_n\to+\infty$ as $n\to+\infty$, and
\begin{equation}\label{vL-close-psi}
\left|v_{\epsilon,y_n,L_n}(\tilde{T}_{\epsilon,y_n,L_n},x_n)- \psi(x_n-L_n,y_n)\right| > \frac{\epsilon}{2}\ \hbox{ for all }n\in\N.
\end{equation}
Remember that, for each $n\in\N$, the inequalities stated in~\eqref{bound-vLn} hold true. By $1$-periodicity of the functions $v_{\epsilon,y,L}$ and $\psi(\cdot,y)$ with respect to $y$, we may assume without loss of generality that~$(y_n)_{n\in\N}\subset [0,1]$.

We first claim that the sequence $(x_n-L_n)_{n\in\N}$ is bounded. Suppose by way of contradiction that $x_n-L_n\to+\infty$ as $n\to+\infty$ (as we will sketch below, the case where $x_n-L_n\to-\infty$ as~$n\to+\infty$ can be treated analogously). Choosing $t=\tilde{T}_{\epsilon,y_n,L_n}$ (notice that $\tilde{T}_{\epsilon,y_n,L_n}\to+\infty$ as~$n\to+\infty$ by~\eqref{define-T} and Lemma~\ref{compare-T-tT}) and $x=x_n$ in the second inequality of~\eqref{bound-vLn}, we have
$$\baa{rcl}
0 & \!\!\!\le\!\!\! & v_{\epsilon,y_n,L_n}(\tilde{T}_{\epsilon,y_n,L_n},x_n)\vspace{3pt}\\
& \!\!\!\le\!\!\! & \displaystyle\psi\Big(x_n-L_n+X_{y_n,L_n}(T_{L_n})-X_{y_n,L_n}(\tilde{T}_{\epsilon,y_n,L_n})+\eta_{\epsilon,L_n}(\tilde{T}_{\epsilon,y_n,L_n}),\,y_n+ \frac{X_{y_n,L_n}(\tilde{T}_{\epsilon,y_n,L_n})}{L_n}\Big)\vspace{3pt}\\
& & +q_{\epsilon,L_n}(\tilde{T}_{\epsilon,y_n,L_n}).\eaa$$
Due to Lemma~\ref{compare-T-tT} and the definition of $X_{y_n,L_n}(t)$, the sequence $(X_{y_n,L_n}(T_{L_n})-X_n(\tilde{T}_{\epsilon,y_n,L_n}))_{n\in\N}$ is bounded. Moreover, by Remark~\ref{rm-bound-eta} and~\eqref{TTtilde}, the sequence $(\eta_{\epsilon,L_n}(\tilde{T}_{\epsilon,y_n,L_n}))_{n\in\N}$ is also bounded. Thus, since $\psi(\xi,y)\to 0$ as $\xi\to+\infty$ uniformly in $y\in\R$ by Proposition~\ref{homo-property}~(ii), we have
$$\psi\Big(x_n-L_n+X_{y_n,L_n}(T_{L_n})-X_{y_n,L_n}(\tilde{T}_{\epsilon,y_n,L_n})+\eta_{\epsilon,L_n}(\tilde{T}_{\epsilon,y_n,L_n}),y_n+ \frac{X_{y_n,L_n}(\tilde{T}_{\epsilon,y_n,L_n})}{L_n}\Big) \to 0$$
as $n\to+\infty$. This together with $q_{\epsilon,L_n}(\tilde{T}_{\epsilon,y_n,L_n})\to 0$ as $n\to+\infty$ (from~\eqref{function-q} and Lemma~\ref{compare-T-tT}) yields $v_{\epsilon,y_n,L_n}(\tilde{T}_{\epsilon,y_n,L_n},x_n)\to0$ as~$n\to+\infty$, hence $v_{\epsilon,y_n,L_n}(\tilde{T}_{\epsilon,y_n,L_n},x_n)-\psi(x_n-L_n,y_n) \to 0$ as~$n\to+\infty$, which is a contradiction with~\eqref{vL-close-psi}.

Similarly, if $x_n-L_n\to-\infty$ as $n\to+\infty$, by using the first inequality of~\eqref{bound-vLn}, one can conclude that $v_{\epsilon,y_n,L_n}(\tilde{T}_{\epsilon,y_n,L_n},x_n)\to1$ and $\psi(x_n-L_n,y_n) \to 1$ as $n\to+\infty$, which is also a contradiction with~\eqref{vL-close-psi}. Therefore, the sequence $(x_n-L_n)_{n\in\N}$ is bounded.

Next, for each $n\in\N$, we set $w_n(t,x)=v_{\epsilon,y_n,L_n}(t+\tilde{T}_{\epsilon,y_n,L_n},x+L_n)$ for $(t,x)\in[-\tilde{T}_{\epsilon,y_n,L_n},+\infty)\times\R$. Clearly, $w_n(0,0)=1/2$ and $w_n$ solves
$$\partial_t w_n= a\partial_{xx} w_n+ f\left(y_n+\frac{x}{L_n},w_n\right) \,\,\hbox{ for }t>-\tilde{T}_{\epsilon,y_n,L_n}\hbox{ and }x\in\R.$$
It is also easily seen from~\eqref{bound-vLn} that, for any $t\in\R$, $x\in\R$ and any $n\in\N$ such that $t+\tilde{T}_{\epsilon,y_n,L_n}\ge0$, $w_n(t,x)$ satisfies
$$\baa{rcl}
w_{n}(t,x) & \!\!\!\geq\!\!\! & \displaystyle\psi\Big(x\!+\!X_{y_n,L_n}(T_{L_n})\!-\!X_{y_n,L_n}(t\!+\!\tilde{T}_{\epsilon,y_n,L_n})\!-\!\eta_{\epsilon,L_n}(t\!+\!\tilde{T}_{\epsilon,y_n,L_n}),y_n\!+\!\frac{X_{y_n,L_n}(t\!+\!\tilde{T}_{\epsilon,y_n,L_n})}{L_n} \Big)\vspace{3pt}\\
& \!\!\!\!\!\! & -\,q_{\epsilon,L_n}(t+\tilde{T}_{\epsilon,y_n,L_n})\eaa$$
and
$$\baa{rcl}
w_{n}(t,x) & \!\!\!\leq\!\!\! & \displaystyle\psi\Big(x\!+\!X_{y_n,L_n}(T_{L_n})\!-\!X_{y_n,L_n}(t\!+\!\tilde{T}_{\epsilon,y_n,L_n})\!+\!\eta_{\epsilon,L_n}(t\!+\!\tilde{T}_{\epsilon,y_n,L_n}),y_n\!+\!\frac{X_{y_n,L_n}(t\!+\!\tilde{T}_{\epsilon,y_n,L_n})}{L_n} \Big)\vspace{3pt}\\
& \!\!\!\!\!\! & +\,q_{\epsilon,L_n}(t+\tilde{T}_{\epsilon,y_n,L_n}).\eaa$$
Thanks to Lemma~\ref{compare-T-tT} and the boundedness of the sequence $(x_n-L_n)_{n\in\N}$, one can find some $\tau_{\infty}\in\R$, $\xi_{\infty}\in\R$ and $y_{\infty}\in [0,1]$, such that, up to extraction of some subsequence,
$$\tilde{T}_{\epsilon,y_n,L_n}-T_{L_n} \to \tau_{\infty},\quad\ x_n-L_n \to \xi_{\infty}\quad\hbox{and}\quad y_n\to y_{\infty}\quad\hbox{as }n\to+\infty.$$
Furthermore, by standard parabolic estimates and possibly up to extraction of a further subsequence, there is a function $w_{\infty}\in C^{1;2}_{t;x}(\R^2)$ such that $w_n\to w_{\infty}$ in $C^{1;2}_{t;x;loc}(\R^2)$ as $n\to+\infty$. Clearly, $w_{\infty}$ is an entire solution of
\begin{equation}\label{bis-homo}
\partial_t w_\infty= a\partial_{xx} w_\infty+ f\left(y_{\infty},w_\infty\right) \,\,\hbox{ for }t\in\R\hbox{ and }x\in\R.
\end{equation}
Notice that for any $t\in\R$, $\lim_{n\to+\infty}q_{\epsilon,L_n}(t+\tilde{T}_{\epsilon,y_n,L_n})=0$ and there exists a constant $C>0$ independent of $t$ such that
\begin{equation}\label{esti-eta-n}
-C\leq \liminf_{n\to+\infty}\,\eta_{\epsilon,y_n}(t+\tilde{T}_{\epsilon,y_n,L_n}) \leq \limsup_{n\to+\infty} \eta_{\epsilon,y_n}(t+\tilde{T}_{\epsilon,y_n,L_n})  \leq 0
\end{equation}
(indeed, since the sequence $(T_{L_n}-\tilde{T}_{\epsilon,y_n,L_n})_{n\in\N}$ is bounded by Lemma~\ref{compare-T-tT} and since $T_{L_n}\to+\infty$ as~$n\to+\infty$, it follows that, for any $t\in\R$, one has $0\leq t+\tilde{T}_{\epsilon,y_n,L_n}\leq 2T_{L_n}$ for all large $n$; then~\eqref{esti-eta-n} follows immediately from Remark~\ref{rm-bound-eta}). Passing to the limit as $n\to+\infty$ in the above inequalities on $w_n(t,x)$, it follows from Lemmas~\ref{compare-T-tT} and~\ref{property-X}, together with the continuity of the map $z\mapsto c(z)$ and the $1$-periodicity of~$\psi(\xi,y)$ with respect to $y$, that
$$\psi\left(x-c(y_{\infty})(t+\tau_{\infty})+C, y_{\infty}\right) \leq w_{\infty}(t,x) \leq  \psi\left(x-c(y_{\infty})(t+\tau_{\infty})-C, y_{\infty}\right)$$
for all $(t,x)\in\R^2$. Namely, $w_{\infty}$ is an entire solution of~\eqref{bis-homo} which is trapped between two shifts of the corresponding traveling front $\psi(x-c(y_\infty)t,y_{\infty})$. Remember that the reaction $f(y_\infty,\cdot)$ is of the  bistable type from~(A1). It then follows from~\cite[Theorem~3.1]{bh07} that there exists $x_0\in\R$ such that $w_{\infty}(t,x)\equiv \psi(x-c(y_{\infty})t+x_0,y_{\infty})$ for all $(t,x)\in\R^2$. Furthermore, notice that~$w_n(0,0)=1/2$ for each $n\in\N$, hence $w_{\infty}(0,0)=1/2$. This together with the normalization condition~\eqref{normal-homo-wave} implies that $x_0$ has to be $0$. Therefore, we have
$$w_{\infty}(t,x)\equiv \psi(x-c(y_{\infty})t,y_{\infty})$$
for all $(t,x)\in\R^2$.

On the other hand, it follows from~\eqref{vL-close-psi} that $|w_n(0,x_n-L_n)-\psi(x_n-L_n,y_n)| > \epsilon/2$ for all $n\in\N$. Taking the limit as $n\to+\infty$ yields $|w_{\infty}(0,\xi_{\infty})-\psi(\xi_{\infty},y_{\infty})| \geq \epsilon/2$, which is a contradiction. The proof of Proposition~\ref{compare-v-psi} is thus complete.
\end{proof}

In addition to the previous observations, the last step before doing the proof of Theorem~\ref{front-limit} in Section~\ref{sc-front-limit-4} is the global stability of the pulsating fronts $\phi_L$, as stated below.

\begin{pro}\label{global-stability}
For any $\epsilon\in(0,\epsilon_0]$ with $\epsilon_0\in(0,\delta_0)$ as in Lemmas~$\ref{super-solu}$-$\ref{sub-solu}$, and for any $L\geq L_*$ with $L_*>0$ as in Theorem~$\ref{th-known}$, let $(v_{\epsilon,y,L})_{y\in\R}$ be the family the solutions of~\eqref{Cauchy-y} with initial conditions $(v^0_{\epsilon,y})_{y\in\R}$ satisfying~\eqref{i-condition-1}. Then, there exists a $1$-periodic real-valued function $y\mapsto\xi_{\epsilon,y,L}$ such that
\begin{equation}\label{stability-fm}
\lim_{t\to+\infty}\left\|v_{\epsilon,y,L}(t,\cdot) -\phi_L\left(\cdot-c_Lt+\xi_{\epsilon,y,L},\frac{\cdot}{L}+y\right)\right\|_{L^{\infty}(\R)} =0\ \hbox{ for every }y\in\R.
\end{equation}
\end{pro}

\begin{proof}
Let $\epsilon\in(0,\epsilon_0]\subset(0,\delta_0)$ and $L\geq L_*$ be fixed. For each $y\in\R$, the function
$$u_{\epsilon,y,L}:(t,x)\mapsto u_{\epsilon,y,L}(t,x):=v_{\epsilon,y,L}(t,x-Ly)$$
is a solution of
$$\partial_tu_{\epsilon,y,L}=a\partial_{xx}u_{\epsilon,y,L}+f(x/L,u_{\epsilon,y,L})\ \hbox{ for $t>0$ and $x\in\R$},$$
with initial condition $u_{\epsilon,y,L}(0,\cdot)=v^0_{\epsilon,y}(\cdot-Ly)\in C(\R,[0,1])$ such that $\|u_{\epsilon,y,L}(0,\cdot)-\psi(\cdot-Ly,y)\| \leq \epsilon\le\epsilon_0<\delta_0$, by~\eqref{i-condition-1}. It then follows from Theorem~\ref{th-known}~(iii) that there exists a unique $\tilde{\xi}_{\epsilon,y,L}\in\R$ such that
$$\lim_{t\to+\infty} \left\|u_{\epsilon,y,L}(t,\cdot) -\phi_L\left(\cdot-c_Lt+\tilde{\xi}_{\epsilon,y,L},\frac{\cdot}{L}\right)\right\|_{L^{\infty}(\R)} =0.$$
Choosing $\xi_{\epsilon,y,L}=\tilde{\xi}_{\epsilon,y,L}+Ly$, we see from the above that~\eqref{stability-fm} holds.

Lastly, since $v_{\epsilon,y,L}$ and $\phi_L(\cdot,y)$ are $1$-periodic in $y$, and since $\phi_L(\xi,y)$ is decreasing with respect to~$\xi$, it follows that $\xi_{\epsilon,y,L}$ in~\eqref{stability-fm} is unique and necessarily $1$-periodic in $y$. The proof of Proposition~\ref{global-stability} is thus complete.
\end{proof}


\subsection{Proof of Theorem~\ref{front-limit}}\label{sc-front-limit-4}

For clarity, we proceed with several steps.

\medskip
\noindent{\it Step 1: approximation of the front profiles $\phi_L$ by $\psi(\cdot,y)$.} In this step, we fix any $\epsilon\in(0,\epsilon_0]$, with $\epsilon_0\in(0,\delta_0)$ as in Lemmas~\ref{super-solu}-\ref{sub-solu} and $\delta_0\in(0,1/2)$ as in assumption~(A2). Let
$$L_{4,\epsilon}:=\max\{L_*,L_{3,\epsilon}\},$$
with~$L_*>0$ and $L_{3,\epsilon}>0$ provided by Theorem~\ref{th-known} and Proposition~\ref{compare-v-psi}, respectively.

We shall prove in this step that there exists a real number $L_{5,\epsilon}\in[L_{4,\epsilon},+\infty)$ such that, for every $L\ge L_{5,\epsilon}$, there is a $1$-periodic real-valued function $y\mapsto\xi^*_{\epsilon,y,L}$ such that
\begin{equation}\label{phiL-psi-close}
\left\|\phi_L\left(\cdot+\xi^*_{\epsilon,y,L},\frac{\cdot}{L}+y\right)-\psi(\cdot,y) \right\|_{L^{\infty}(\R)} \leq \epsilon\ \hbox{ for every }y\in\R\hbox{ and }L\ge L_{5,\epsilon}.
\end{equation}

To do so, for any $y\in\R$ and $L\geq L_{4,\epsilon}\ (\ge L_{3,\epsilon}\ge L_{2,\epsilon})$, let $v_{\epsilon,y,L}$ be the solution of~\eqref{Cauchy-y} with initial condition~$v^0_{\epsilon,y}$ satisfying~\eqref{i-condition-1} (for instance, say, take here $v^0_{\epsilon,y}=\psi(\cdot,y)$). By Proposition~\ref{compare-v-psi} and~\eqref{v12}, we have
$$v_{\epsilon,y,L}(\tilde{T}_{\epsilon,y,L},L)=\frac{1}{2} \quad\hbox{and}\quad \left\|v_{\epsilon,y,L}\big(\tilde{T}_{\epsilon,y,L},\cdot+L\big)-\psi(\cdot,y) \right\|_{L^\infty(\R)} \leq \frac{\epsilon}{2}\le\epsilon,$$
where $\tilde{T}_{\epsilon,y,L}$ is defined by~\eqref{define-tildeT}-\eqref{TTtilde}. Since $f$ is $1$-periodic with respect to its first variable, the function $(t,x)\mapsto v_{\epsilon,y,L}(t+\tilde{T}_{\epsilon,y,L},x+L)$ is still a solution of~\eqref{Cauchy-y}, now with the initial condition $v_{\epsilon,y,L}(\tilde{T}_{\epsilon,y,L},\cdot+L)\in C(\R,[0,1])$. Furthermore, since both $v_{\epsilon,y,L}$ and $\tilde{T}_{\epsilon,y,L}$ are $1$-periodic with respect to $y$, one has $v_{\epsilon,y+1,L}(\tilde{T}_{\epsilon,y+1,L},\cdot+L)\equiv v_{\epsilon,y,L}(\tilde{T}_{\epsilon,y,L},\cdot+L)$ in $\R$ for each $y\in\R$. By using Proposition~\ref{compare-v-psi} again, we infer that
$$\left\|v_{\epsilon,y,L}\big(\tilde{T}^2_{\epsilon,y,L},\cdot+2L\big)-\psi(\cdot,y)\right\|_{L^\infty(\R)} \leq\frac{\epsilon}{2}\le\epsilon,$$
where $\tilde{T}^2_{\epsilon,y,L}:=\min\big\{t>\tilde{T}_{\epsilon,y,L}: \, v_{\epsilon,y,L}(t,2L)=1/2\big\}$ (the quantities $\tilde{T}^2_{\epsilon,y,L}$ are well defined real numbers, as are the $\tilde{T}_{\epsilon,y,L}$'s, from the same arguments as in~\eqref{define-tildeT}-\eqref{v12}). Then, a simple induction argument implies that, for any $k\in\N$,
\begin{equation}\label{vL-kL}
\left\|v_{\epsilon,y,L}\big(\tilde{T}^k_{\epsilon,y,L},\cdot+kL,y\big)-\psi(\cdot,y) \right\|_{L^\infty(\R)} \leq \frac{\epsilon}{2},
\end{equation}
where $\tilde{T}^k_{\epsilon,y,L}=\min\big\{t>\tilde{T}^{k-1}_{\epsilon,y,L}: \, v_{\epsilon,y,L}(t,kL)=1/2 \big\}$. Moreover, due to Lemma~\ref{compare-T-tT}, there holds
$$\sup_{y\in\R,\,L\ge L_{4,\epsilon},\,k\in\N}\,\big|\tilde{T}^{k+1}_{\epsilon,y,L}-\tilde{T}^{k}_{\epsilon,y,L}-T_L\big|<+\infty,$$
where $T_L$ is defined by~\eqref{define-T}. Since $T_L\to+\infty$ as $L\to+\infty$ (independently of $y\in\R$), there exists $L_{5,\epsilon}\in[L_{4,\epsilon},+\infty)$ such that $\tilde{T}^k_{\epsilon,y,L}\to+\infty$ as $k\to+\infty$ uniformly in $y\in\R$ and $L\geq L_{5,\epsilon}$.

On the other hand, for each $L\ge L_{5,\epsilon}\ (\ge L_{4,\epsilon}\ge L_*)$, it follows from Proposition~\ref{global-stability} that there exists a $1$-periodic real-valued function $y\mapsto\xi_{\epsilon,y,L}$ such that
$$\lim_{t\to+\infty} \left\|v_{\epsilon,y,L}(t,\cdot) -\phi_L\left(\cdot-c_Lt+\xi_{\epsilon,y,L},\frac{\cdot}{L}+y\right)\right\|_{L^{\infty}(\R)} =0$$
for every $y\in\R$. Therefore, for each $L\ge L_{5,\epsilon}$ and $y\in[0,1)$, there is $k_{\epsilon,y,L}\in\N$ such that
$$\left\|v_{\epsilon,y,L}\big(\tilde{T}^{k_{\epsilon,y,L}}_{\epsilon,y,L},\cdot+k_{\epsilon,y,L}L\big) -\phi_L\left(\cdot+\xi^*_{\epsilon,y,L},\frac{\cdot}{L}+y\right)\right\|_{L^{\infty}(\R)} \leq \frac{\epsilon}{2},$$
where $\xi^*_{\epsilon,y,L}:=k_{\epsilon,y,L}L-c_L\tilde{T}^{k_{\epsilon,y,L}}_{\epsilon,y,L}+\xi_{\epsilon,y,L}$ (we here use the $1$-periodicity of $\phi$ in its second variable).  Combining this with~\eqref{vL-kL} at $k=k_{\epsilon,y,L}$, we obtain that
\begin{equation}\label{phipsi}
\left\|\phi_L\big(\cdot+\xi^*_{\epsilon,y,L},\frac{\cdot}{L}+y\big)-\psi(\cdot,y) \right\|_{L^{\infty}(\R)} \leq \epsilon.
\ee
We can then extend the function $y\mapsto\xi^*_{\epsilon,y,L}$ as a $1$-periodic function, so that~\eqref{phipsi} holds for all $y\in\R$, since $\phi_L(\xi,y)$ and $\psi(\xi,y)$ are $1$-periodic with respect to $y\in\R$. The proof of Step~1 is thus complete.

\medskip
\noindent{\it Step 2: proof of
\be\label{claim2}
\lim_{\epsilon\to 0}\ \Big(\sup_{y\in\R,\,L\geq L_{5,\epsilon}} \big|\xi^*_{\epsilon,y,L}-\zeta_L(y)\big|\Big)=0,
\ee
with $\zeta_L$ given by~\eqref{normmal-front-1}.} Indeed, for each $\epsilon\in(0,\epsilon_0]$, since $\psi(0,y)=\phi_L(\zeta_L(y),y)=1/2$ for all $y\in\R$ and $L\geq L_{5,\epsilon}$,~\eqref{phiL-psi-close} implies in particular that
\begin{equation}\label{xiL-zetaL-close}
\sup_{y\in\R,\,L\geq L_{5,\epsilon}} \Big|\phi_L\left(\xi^*_{\epsilon,y,L},y\right)-\underbrace{\phi_L(\zeta_L(y),y)}_{=1/2} \Big|\leq \epsilon\le\epsilon_0<\frac14.
\end{equation}
Then, by Theorem~\ref{wave-length-1}~(i)-(ii), it follows that there exists a constant $B_1>0$ such that
\begin{equation}\label{xiL-zetaL-bound}
|\xi^*_{\epsilon,y,L}-\zeta_L(y)| \leq B_1 \,\,\hbox{ for all }\epsilon\in(0,\epsilon_0],\ y\in\R\hbox{ and }L\geq L_{5,\epsilon},
\end{equation}
and also that~\eqref{claim2} holds.

\medskip
\noindent{\it Step 3: $\limsup_{L\to+\infty}|L-c_LT_L|<+\infty$, with $T_L$ given by~\eqref{define-T}}.\footnote{With this property, we here recover that $c_L\to c_*$ as $L\to+\infty$.} Notice that for any $y\in\R$ and $L\geq L_{5,\epsilon_0}\ge L_{4,\epsilon_0}=\max\{L_*,L_{3,\epsilon_0}\}$, the function
$$(t,x)\mapsto\phi_L\Big(x-c_Lt+\xi^*_{\epsilon_0,y,L},\frac{x}{L}+y\Big)=U_L\Big(t-\frac{\xi^*_{\epsilon_0,y,L}}{c_L}+\frac{Ly}{c_L},x+Ly\Big)$$
is the solution of~\eqref{Cauchy-y} with initial condition $\phi_L(\cdot+\xi^*_{\epsilon_0,y,L},\cdot/L+y)$. Since $L\ge L_{5,\epsilon_0}\ge L_{1,\epsilon_0}$, it follows from~\eqref{phiL-psi-close} and Lemmas~\ref{super-solu}-\ref{sub-solu} (with their notations) that
\begin{equation}\label{comp-vpm-phiL}
v^-_{\epsilon_0,y,L}(t,x) \leq \phi_L\left(x-c_Lt+\xi^*_{\epsilon_0,y,L},\frac{x}{L}+y\right)\leq  v^+_{\epsilon_0,y,L}(t,x) \,\,\hbox{ for all }t\geq 0\hbox{ and }x\in\R.
\end{equation}
In particular, choosing $t=T_L$ and $x=L$ yields
$$\psi(-\eta_{\epsilon_0,L}(T_L),y)-q_{\epsilon_0,L}(T_L) \leq \phi_L(L-c_LT_L+\xi^*_{\epsilon_0,y,L},y)\leq \psi(\eta_{\epsilon_0,L}(T_L),y)+q_{\epsilon_0,L}(T_L).$$
Since $\lim_{L\to+\infty}q_{\epsilon_0,L}(T_L)=0$ by~\eqref{function-q} and since $\sup_{L\ge L_{5,\epsilon_0}}|\eta_{\epsilon_0,L}(T_L)|<+\infty$ by Remark~\ref{rm-bound-eta}, it then follows from Proposition~\ref{homo-property} (namely, the continuity of $\psi:\R^2\to(0,1)$ and its $1$-periodicity in its second variable) that there exist $\kappa_0 \in (0,1/2)$ and $\bar{L}\ge L_{5,\epsilon_0}$ such that
$$\kappa_0 \leq \phi_L\left(L-c_LT_L+\xi^*_{\epsilon_0,y,L},y\right) \leq 1-\kappa_0\,\,\hbox{ for all $y\in\R$ and $L\geq \bar{L}$}.$$
Hence, by Theorem~\ref{wave-length-1}~(i), there exists $B_2>0$ such that
$$\big|L-c_LT_L+\xi^*_{\epsilon_0,y,L}-\zeta_L(y)\big| \leq B_2\,\,\hbox{ for all $y\in\R$ and $L\geq\bar{L}$}.$$
This together with~\eqref{xiL-zetaL-bound} implies that  $|L-c_LT_L|\leq B_1+B_2$ for all $L\geq\bar{L}$, which yields the desired result.

\medskip
\noindent{\it Step 4: proof the convergence~\eqref{eq-front-limit}}. We fix any $\sigma>0$. We have to show the existence of $L^*_\sigma>0$ such that for any $L\geq L^*_\sigma$,
\begin{equation}\label{aim-sigma}
\sup_{y\in \R} \left\|\phi_L \left(\cdot+ \zeta_L(y),  \frac{\cdot}{L}+y \right) -\psi(\cdot,y) \right\|_{L^{\infty}(\R)}\leq \sigma.
\end{equation}
To do so, we first recall that the map $L\mapsto c_L$ is continuous and positive from Theorem~\ref{th-known}~(i) and~(iv), and that $c_L\to c_*>0$ by Theorem~\ref{speed-limit} (or by Step~3 of the present proof), hence $\sup_{L\ge L_*}c_L^{-1}<+\infty$. It then follows from standard parabolic estimates that $M\!:=\!\sup_{L\ge L_*}\!\|\partial_\xi\phi_L\|_{L^\infty(\R^2)}\!=\!\sup_{L\ge L_*}\!c_L^{-1}\|\partial_tU_L\|_{L^\infty(\R^2)}\!\in\![0,+\infty)$. Next, by Step~2, there is $\epsilon_\sigma\in(0,\min\{\epsilon_0,\sigma/2\}]$ such that
$$M\times|\xi^*_{\epsilon_\sigma,y,L}-\zeta_L(y)|\le\frac{\sigma}{2}\ \hbox{ for all }y\in\R\hbox{ and }L\ge L_{5,\epsilon_\sigma}.$$
Hence,
$$\left\|\phi_L\left(\cdot+ \zeta_L(y),  \frac{\cdot}{L}+y \right)-\phi_L\left(\cdot+\xi^*_{\epsilon_\sigma,y,L},  \frac{\cdot}{L}+y \right)\right\|_{L^{\infty}(\R)}\leq\frac{\sigma}{2}\ \hbox{ for all }y\in\R\hbox{ and }L\ge L_{5,\epsilon_\sigma}.$$
Then,~\eqref{phiL-psi-close} (applied with $\epsilon_\sigma$) and the inequality $\epsilon_\sigma\le\sigma/2$ give~\eqref{aim-sigma} for all $L\ge L^*_\sigma:=L_{5,\epsilon_\sigma}>0$. This ends the proof of Step 4.

\medskip
\noindent{\it Step 5: the convergence~\eqref{asym-zeta} holds locally uniformly in $x\in\R$.} Let us first show that, for any given $A>0$,
\begin{equation}\label{xiL-zetaL}
\sup_{y\in\R,\,L\geq L_{5,\epsilon},\,x\in[0,A]} \left|  \xi^*_{\epsilon,y,L}+L\int_{y}^{y+x/L} \left( 1-\frac{c_L}{c(s)}\right)ds-\zeta_L\left(y+\frac{x}{L}\right)   \right|\ \to0\ \hbox{ as }\epsilon\to0.
\end{equation}
For any $\epsilon\in(0,\epsilon_0]$, $y\in\R$, $L\geq L_{5,\epsilon}$ and $x\in[0,A]$, define
\be\label{deftyLx}
t_{y,L,x}:=L\int_{y}^{y+x/L}c^{-1}(s)ds.
\ee
Clearly, $t_{y,L,x}\ge0$. It is easily seen from the definition of $X_{y,L}(t)$ in~\eqref{ode-speed} that $X_{y,L}(t_{y,L,x})=x$. Choosing $x$ and $t=t_{y,L,x}$ in~\eqref{comp-vpm-phiL} (these inequalities actually hold as in Step~3 with $\epsilon_0$ replaced by any $\epsilon\in(0,\epsilon_0]$, and for any $L\ge L_{5,\epsilon}$) gives
\begin{equation}\label{choose-x-t}
\left.\baa{l}
\displaystyle \psi\left(-\eta_{\epsilon,L}(t_{y,L,x}),\,y+\frac{X_{y,L}(t_{y,L,x})}{L} \right)-q_{\epsilon,L}(t_{y,L,x})\vspace{3pt}\\
\qquad\qquad\qquad\leq\displaystyle \phi_{L} \left(L\int_{y}^{y+x/L} \left( 1-\frac{c_L}{c(s)}\right)ds  +\xi^*_{\epsilon,y,L}, \frac{x}{L}+y\right)\vspace{3pt}\\
\qquad\qquad\qquad\displaystyle  \leq \psi\left(\eta_{\epsilon,L}(t_{y,L,x}),\,y+\frac{X_{y,L}(t_{y,L,x})}{L} \right)+q_{\epsilon,L}(t_{y,L,x}). \eaa\right.
\end{equation}
Moreover, by the normalization conditions~\eqref{normal-homo-wave} and~\eqref{normmal-front-1}, we have
$$\psi\left(0,y+\frac{X_{y,L}(t_{y,L,x})}{L}\right)=\phi_L\left(\zeta_L\left(\frac{x}{L}+y\right),\frac{x}{L}+y\right)=\frac{1}{2}.$$
It then follows that
\begin{equation}\label{phiL-vpm-t0}
\baa{l}
\displaystyle \left| \phi_{L}\left(L\int_{y}^{y+x/L} \left( 1-\frac{c_L}{c(s)}\right)ds+\xi^*_{\epsilon,y,L}, \frac{x}{L}+y\right)  - \phi_L\left(\zeta_L\left(\frac{x}{L}+y\right),\frac{x}{L}+y\right) \right|   \vspace{6pt}\\
\qquad\leq  \displaystyle\max\left\{\left| \psi\left(\eta_{\epsilon,L}(t_{y,L,x}),\,y+\frac{X_{y,L}(t_{y,L,x})}{L} \right)-\psi\left(0,\,y+\frac{X_{y,L}(t_{y,L,x})}{L} \right) \right|,\right.\vspace{3pt}\\
\qquad\qquad\ \ \ \ \ \,\displaystyle\left.\left| \psi\left(-\eta_{\epsilon,L}(t_{y,L,x}),\,y+\frac{X_{y,L}(t_{y,L,x})}{L} \right)-\psi\left(0,\,y+\frac{X_{y,L}(t_{y,L,x})}{L} \right) \right|\right\}\vspace{3pt}\\
\qquad\ \ +q_{\epsilon,L}(t_{y,L,x}).\eaa
\end{equation}
By definition of $t_{y,L,x}$, we have $t_{y,L,x}/L \to 0$ as $L\to+\infty$ uniformly in $y\in\R$ and $x\in[0,A]$. Furthermore, thanks to~\eqref{condition-q} and~\eqref{function-eta}, and since $L_{5,\epsilon}\ge L_{1,\epsilon}\to+\infty$ as $\epsilon\to0$ by~\eqref{defL1eps}, one infers that
$$\sup_{y\in\R,\,L\ge L_{5,\epsilon},\,x\in[0,A]}\ \Big(\underbrace{q_{\epsilon,L}(t_{y,L,x})}_{\ge0}+|\eta_{\epsilon,L}(t_{y,L,x})|\Big)\ \to 0\ \hbox{ as }\epsilon\to0.$$
Passing to the limit as $\epsilon\to 0$ in~\eqref{phiL-vpm-t0} and using the boundedness of $\partial_\xi\psi$ in $\R^2$ (which itself follows from standard elliptic estimates), one infers that
$$\sup_{y\in\R,\,L\geq L_{5,\epsilon},\,x\in[0,A]} \left| \phi_{L}\!\left(L\int_{y}^{y+x/L}\!\!\!\left( 1-\frac{c_L}{c(s)}\right)ds+\xi^*_{\epsilon,y,L}, \frac{x}{L}+y\right)\!\!-\underbrace{\phi_L\left(\zeta_L\left(\frac{x}{L}+y\right),\frac{x}{L}+y\right)}_{=1/2}\right|\to0$$
as $\epsilon\to0$. This together with Theorem~\ref{wave-length-1}~(i)-(ii) implies~\eqref{xiL-zetaL} (the proof is actually similar to that of Step 2; therefore, we omit the details). Combining~\eqref{claim2} and~\eqref{xiL-zetaL}, we immediately obtain that~\eqref{asym-zeta} holds in the case where $x\in[0,A]$.

It remains to consider the case where $x\in[-A,0]$. Clearly, the convergence~\eqref{xiL-zetaL} holds with~$x$ replaced by~$A$ and by~$x+A\in[0,A]$ as well. Since $\xi^*_{\epsilon,y,L}$ is independent of the variable $x$, it then follows that
$$\sup_{y\in \R,\,L\ge L_{5,\epsilon},\,x\in[-A,0]} \left| \zeta_L\left(y+\frac{x+A}{L}\right) - \zeta_L\left(y+\frac{A}{L}\right)-L\int_{y+A/L}^{y+(x+A)/L} \left( 1-\frac{c_L}{c(s)}\right)ds   \right|\ \to0$$
as $\epsilon\to0$. Replacing $y$ by $y+A/L$, we obtain~\eqref{asym-zeta} in the case where $x\in[-A,0]$. The proof of Step 5 is thus compete.

\medskip
\noindent{\it Step 6: proof of the estimate~\eqref{zetaL-zetaLx}.} By~\eqref{claim2} and the continuity of the map $(L,y)\mapsto\zeta_L(y)$ in~$[L_*,+\infty)\times\R$ and its $1$-periodicity in $y$, in order to show~\eqref{zetaL-zetaLx}, it suffices to prove the existence of $\epsilon_1\in(0,\epsilon_0]$ such that
\begin{equation}\label{xiL-zetaL-close-ge}
\sup_{\epsilon\in(0,\epsilon_1],\, y\in\R,\,L\geq L_{5,\epsilon}} \left\| \xi^*_{\epsilon,y,L}+L\int_{y}^{y+\cdot/L} \left( 1-\frac{c_L}{c(s)}\right)ds-\zeta_L\left(y+\frac{\cdot}{L}\right)  \right\|_{L^{\infty}([-L,L])} <+\infty.
\end{equation}

Consider any $\epsilon\in(0,\epsilon_0]$, $y\in\R$, $L\geq L_{5,\epsilon}$, and $x\in[0,L]$. Let $t_{y,L,x}$ be as in~\eqref{deftyLx}. There holds $0\le t_{y,L,x}\le T_L$ and $X_{y,L}(t_{y,L,x})=x$, since $X_{y,L}(T_L)=L\ge x$ and $X_{y,L}$ is increasing. It is also easily seen that the inequalities~\eqref{choose-x-t} remain valid. Moreover, by~\eqref{condition-q} and Remark~\ref{rm-bound-eta}, we have $0<q_{\epsilon,L}(t_{y,L,x}) \leq \epsilon$ and $-B_3\leq\eta_{\epsilon,L}(T_L) \leq  \eta_{\epsilon,L}(t_{y,L,x})\le0$,  where $B_3$ is a positive constant independent of $\epsilon$, $y$, $L$ and $x$. Then, by the continuity of the map $\psi:\R^2\to(0,1)$ and its $1$-periodicity with respect to its second variable, there exist $\epsilon_1\in(0,\epsilon_0]$ and $\kappa_1\in(0,1/2)$ such that
$$\kappa_1 \leq \phi_{L}\left(L\int_{y}^{y+x/L} \left( 1-\frac{c_L}{c(s)}\right)ds+\xi^*_{\epsilon,y,L}, \frac{x}{L}+y\right) \leq 1-\kappa_1$$
for all $\epsilon\in(0,\epsilon_1]$, $y\in\R$, $L\ge L_{5,\epsilon}$, and $x\in[0,L]$. It further follows from Theorem~\ref{wave-length-1} that there exists $B_4\in[0,+\infty)$ such that
\begin{equation}\label{xi-zeta-int}
\sup_{\epsilon\in(0,\epsilon_1],\,y\in\R,\, L\geq L_{5,\epsilon},\,x\in[0,L]} \left| \xi^*_{\epsilon,y,L}+L\int_{y}^{y+x/L} \left( 1-\frac{c_L}{c(s)}\right)ds-\zeta_L\left(y+\frac{x}{L}\right) \right | \leq B_4.
\end{equation}
This means that~\eqref{xiL-zetaL-close-ge} holds true with $[-L,L]$ replaced by $[0,L]$.

Next, for any $x\in [-L,0]$, we have $x+L \in [0,L]$. Since $z\mapsto\zeta_L(z)$ and $z\mapsto c(z)$ are $1$-periodic, it then follows from~\eqref{xi-zeta-int} and the formula $T_L=L/c_*$ that
\begin{equation*}
\left.\baa{ll}
& \displaystyle\left |\xi^*_{\epsilon,y,L} -\zeta_L\left(y+\frac{x}{L}\right) +L\int_{y}^{y+x/L} \left( 1-\frac{c_L}{c(s)}\right)ds \right|  \vspace{6pt}\\
=  \!\!\!&\displaystyle \left |\xi^*_{\epsilon,y,L} -\zeta_L\left(y+\frac{x+L}{L}\right) +L\int_{y}^{y+(x+L)/L} \left( 1-\frac{c_L}{c(s)}\right)ds -(L-c_LT_L) \right| \vspace{6pt}\\
\leq   \!\!\!&\displaystyle B_4+ |L-c_LT_L|,
\eaa\right.
\end{equation*}
for all $\epsilon\in(0,\epsilon_1]$, $y\in\R$, $L\ge L_{5,\epsilon}$, and $x\in[-L,0]$. By the conclusion of Step~3 and the continuity of $c_L$ and $T_L$ with respect to $L$, the last quantity is then bounded uniformly with respect to $\epsilon\in(0,\epsilon_1]$, $y\in\R$, $L\ge L_{5,\epsilon}$, and $x\in[-L,0]$. Combining the above, we obtain~\eqref{xiL-zetaL-close-ge}, hence~\eqref{zetaL-zetaLx} is proved. This ends the proof of Theorem~\ref{front-limit}.\hfill$\Box$


\section{Appendix}\label{sec5}

The appendix is devoted to the proof of Proposition~\ref{homo-property}. Statement (i) follows easily from the assumption that $a(y)$ and $f(y,u)$ are $1$-periodic in $y$ and the classical theory on homogeneous bistable traveling waves established in \cite{aw}. It is also not difficult to prove statements (ii)-(iii) by using phase-plane arguments. As we failed to find a direct proof in the literature, we just sketch the main arguments below for the sake of completeness.

\begin{proof}[{\it Proof of Proposition~$\ref{homo-property}$~{\rm(ii)-(iii)}}] 
For each $y\in\R$, with writting $q(\cdot,y)=\psi(\cdot,y)$, the equation in~\eqref{traveling-wave} is equivalent to the following first-order system 
\begin{equation}\label{first-ode}
  \frac{dq(\xi,y)}{d\xi}=p(\xi,y),\quad \frac{dp(\xi,y)}{d\xi}=-\frac{c(y)}{a(y)}\,p(\xi,y)-\frac{f(y,q(\xi,y))}{a(y)} \,\, \hbox{ for }\, \xi\in\R.
\end{equation} 
Since $p(\xi,y)<0$ for all $\xi\in\R$ and $y\in\R$, there is a one-to-one corresponding between the trajectories $T_y(\xi)=(q(\xi,y),p(\xi,y))$ of \eqref{first-ode} and the integral curves $I_y(q)=(q,P(q,y))$ of the equation 
\begin{equation}\label{phase}
	\frac{dP(q,y)}{dq}=-\frac{c(y)}{a(y)}-\frac{f(y,q)}{a(y)\,P(q,y)} 
\end{equation}
in the semi-strip $S=(0,1)\times (-\infty,0)$ of the $(q,p)$-plane. It follows directly from~\cite[Proposition~2.3]{g} that $c(y)$ is continuous in $y\in\R$, and hence, the right hand side of \eqref{phase} is continuous with respect to $(q,P,y)\in S\times \R$. This together with the fact that for each $y\in\R$,  $I_y$ is the unique curve lying in $S$ and connecting $(0,0)$ and $(1,1)$, implies that the function $(q,y)\in (0,1)\times \R\mapsto P(q,y)$ is continuous. Now, coming back to equation \eqref{first-ode} and recalling that $q(0,y)=1/2$ by \eqref{normal-homo-wave}, we obtain that $q(\xi,y)$ is continuous in $(\xi,y)\in \R^2$. Then, applying standard elliptic estimates to the equation in \eqref{traveling-wave}, one can conclude that the map $y\mapsto \psi(\cdot,y)$ from $\R$ to $C^2_{loc}(\R)$ is continuous. Using this continuity and arguing by contradiction, one can derive that $\psi(\xi,y)$ tends to its limits and $\partial_\xi\psi(\xi,y)$ tends to $0$ as $\xi\to\pm \infty$  uniformly in $y\in\R$. This immediately gives the first part of~(ii), and~(iii) as well. Furthermore, together with \cite[Lemmas 3.2-3.3]{po} (see also \cite{aw}), the above convergences also imply that 
$$\frac{\partial_{\xi}\psi(\xi,y)}{\psi(\xi,y)}\to -\lambda_1(y) \, \hbox{ as }\,\xi \to +\infty \quad\hbox{and}\quad \frac{\partial_{\xi}\psi(\xi,y)}{1-\psi(\xi,y)}\to -\lambda_2(y)  \, \hbox{ as }\,\xi \to -\infty, $$  
where the convergences are taken uniformly in $y\in\R$, and $\lambda_1(y)$, $\lambda_2(y)$ are positive constants given by 
$$\lambda_1(y)=\frac{c(y)+\sqrt{c^2(y)-4 a(y)\partial_uf(y,0)}}{2a(y)} \quad\hbox{and}\quad \lambda_2(y)=\frac{-c(y)+\sqrt{c^2(y)-4 a(y)\partial_uf(y,1)}}{2a(y)}.  $$
Thanks to the assumption (A2) and the fact that the functions $y\mapsto a(y)$ and $y\mapsto c(y)$ are continuous positive and 1-periodic, choosing $\mu_1>0$ and $\mu_2>0$ such that
$$0<\mu_1< \min_{y\in\R} \frac{c(y)+\sqrt{c^2(y)+4\gamma_0 a(y)}}{2a(y)}\quad \hbox{and}\quad 0<\mu_2< \min_{y\in\R} \frac{-c(y)+\sqrt{c^2(y)+4\gamma_0 a(y)}}{2a(y)}, $$ 
one obtains \eqref{psi-exp}. This gives the second part of (ii). The proof of Proposition~$\ref{homo-property}$~(ii)-(iii) is thus complete.  
\end{proof}

\begin{rem}\label{uniform-esti-psi}
Before proceeding with the proof of statement~{\rm{(iv)}}, we collect some easy corollaries of the estimate~\eqref{psi-exp}, which will be used frequently later. First of all, by~\eqref{psi-exp} and standard elliptic estimates applied to~\eqref{traveling-wave} and its derivative, we can find a constant $\tilde{C}_1>0$ $($independent of $y\in\R$$)$ such that
\begin{equation}\label{esit-psi-l1}
\|1-\psi(\cdot,y)\|_{L^2((-\infty,0))}+\|\psi(\cdot,y)\|_{L^2((0,+\infty))}+ \| \partial_{\xi} \psi(\cdot,y)\|_{H^2(\R)} \leq \tilde{C}_1
\end{equation}
for all $y\in\R$. Furthermore, for $\xi\in\R$ and $y\in\R$, letting
\begin{equation}\label{adjoint-ker}
w^*(\xi,y):={\rm exp}\left(\frac{c(y)}{a(y)}\xi\right)\partial_{\xi}\psi(\xi,y),
\end{equation}
one observes from the proof of~\eqref{psi-exp} and standard elliptic estimates that there is a constant $\tilde{C}_2>0$ such that
\be\label{estpartpsi}
|\partial_\xi\psi(\xi,y)|\le\tilde{C}_2\min\{e^{-((c(y)+\sqrt{c^2(y)+4\gamma_0a(y)})/(2a(y)))\,\xi},1\}\ \hbox{ for all $(\xi,y)\in\R^2$}.
\ee
Therefore, $w^*(\xi,y)$ decays no slower than $\me^{-\tilde{\mu}_1\xi}$ as $\xi\to +\infty$, where $\tilde{\mu}_1$ is given by
$$\tilde{\mu}_1= \min_{y\in\R} \frac{-c(y)+\sqrt{c^2(y)+4\gamma_0 a(y)}}{2a(y)}>0,$$
and, as above, there is a constant  $\tilde{C}_3>0$ such that $\left\|w^*(\cdot,y)\right\|_{H^2(\R)} \leq \tilde{C}_3$ for all $y\in\R$.
\end{rem}

Now, we turn to the proof of statement (iv) of Proposition~\ref{homo-property}. The proof shares some similarities with the proof of~\cite[Theorem~1.2]{dhz1} which established the existence of pulsating fronts of~\eqref{eqL} when $L$ is small and the convergence of those fronts as $L\to 0$ by using the implicit function theorem (see also~\cite[Theorem~1.1]{he} for the study of pulsating fronts in perforated domains). Here, we will apply the implicit function theorem under a similar setting to show the $C^1$-smoothness of the homogeneous fronts $(\psi(\cdot,y),c(y))$ with respect to $y\in\R$.

Let us first introduce a family of auxiliary operators.  In the sequel, we fix a real number $\beta>0$. For any $c\in\R$ and $y\in\R$, we define
$$M_{c,y}(v)=a(y)v''+cv'-\beta v \ \ \hbox{for }v\in H^2(\R).$$
Clearly, each $M_{c,y}$ maps $H^2(\R)$ into $L^2(\R)$. In the following lemma, we present some basic properties of this operator.

\begin{lem}\label{prop-M}
Fix $\beta>0$. We have
\begin{itemize}
\item[{\rm (i)}] For any $c\in\R$ and $y\in\R$, the operator $M_{c,y}: H^2(\R)\to L^2(\R)$ is invertible. Furthermore, for every $A>0$, there exists a constant $C>0$ $($depending on $\beta$ and $A$, but independent of~$y$$)$ such that for any $c\in[-A,A]$, $y\in\R$ and $g\in L^2(\R)$,
\begin{equation}\label{es-inverse-M}
\|M_{c,y}^{-1}(g)\|_{H^2(\R)} \leq  C \|g\|_{L^2(\R)}.
\end{equation}
\item[{\rm (ii)}] For any $g\in L^2(\R)$,
\begin{equation}\label{Mn-converge}
M_{c_n,y_n}^{-1}(g_n) \to M_{c,y}^{-1}(g)\quad  \hbox{in } H^2(\R)
\end{equation}
as $n\to+\infty$ for all sequences $(g_n)_{n\in\N} \subset L^2(\R)$, $(c_n)_{n\in\N}\subset \R$, $(y_n)_{n\in\R}\subset \R$ such that $\|g_n- g\|_{L^2(\R)} \to 0$, $c_n\to c$ and $y_n\to y$ as $n\to+\infty$. Furthermore, the above convergence is uniform in $(g,c,y)\in B_A\times\R$ for any $A>0$, where $B_A$ is the ball given by $B_A =\{ (g,c)\in L^2(\R)\times\R:  \|g\|_{L^2(\R)}+ |c| \leq A  \}$.
\end{itemize}
\end{lem}

\begin{proof}
The proof of the invertibility of $M_{c,y}$ follows from similar arguments used in the proof of~\cite[Lemma~3.1]{dhz1}; therefore, we omit the details. To obtain the estimate~\eqref{es-inverse-M}, let $v=M_{c,y}^{-1}(g)$. Integrating $M_{c,y}(v)=g$ against $v$ over $\R$ gives $\int_{\R} a(y)(v')^2+\beta v^2=-\int_{\R}gv$, whence
$$\int_{\R} a(y)(v')^2+\frac{\beta}{2} v^2 \leq \frac{1}{2\beta}\int_{\R}g^2.$$
This together with  $M_{c,y}(v)=g$ implies that
$$\|v \|_{L^2(\R)} \leq \frac{1}{\beta} \,\|g\|_{L^2(\R)},\quad \|v' \|_{L^2(\R)} \leq \sqrt{\frac{1}{2\beta \min_{x\in\R}a(x)}} \,\|g\|_{L^2(\R)}    $$
and
$$\|v''\|_{L^2(\R)}  \leq \frac{1}{\min_{x\in\R}a(x)} \left(\sqrt{\frac{c^2}{2\beta \min_{x\in\R}a(x)}} +1\right) \|g\|_{L^2(\R)}.$$
Thus, part~(i) of Lemma~\ref{prop-M} follows.

Next, we prove the convergence~\eqref{Mn-converge}.  For each $n\in\N$, let $u_n=M_{c_n,y_n}^{-1}(g_n)$ and $u=M_{c,y}^{-1}(g)$.  Let $A>0$ be fixed and $(g,c,y)\in B_A\times\R$. Without loss of generality, we may assume that $(g_n,c_n,y_n)\in B_{2A}\times\R$ for all $n\in\N$.  By the proof of~\eqref{es-inverse-M}, we have $\|u\|_{H^2(\R)} \leq  C_1 \|g\|_{L^2(\R)}$, where $C_1$ is a positive constant depending only on $\beta$ and $A$. Notice that
$$M_{c_n,y_n}(u_n-u)=(a(y_n)-a(y))u''+(c-c_n)u'+(g_n-g)$$
for all $n\in\N$. By the proof of~\eqref{es-inverse-M} again, we find a positive constant $C_2$ depending only on $\beta$ and $A$ such that  $\|u_n-u\|_{H^2(\R)} \leq C_2 \|(a(y_n)-a(y))u''+(c_n-c)u'+(g_n-g)\|_{L^2(\R)}$, whence
$$\|u_n-u\|_{H^2(\R)} \leq C_1C_2\Big(\max_{x\in\R}|a'(x)|\Big)|y_n-y|\|g\|_{L^2(\R)}+C_1C_2|c_n-c|\|g\|_{L^2(\R)}+C_2\|g_n-g\|_{L^2(\R)},$$
for all $n\in\N$. This implies that~\eqref{Mn-converge} holds  uniformly in $(g,c,y)\in B_A\times\R$. The proof of Lemma~\ref{prop-M} is thus complete.
\end{proof}

In addition to $\beta>0$, we also consider in the sequel an arbitrary real number $y_0$. Now, for any $(v,c,y)\in L^2(\R)\times\R\times\R$, we define
$$K(v,c,y):\xi\mapsto K(v,c,y)(\xi):=a(y)\partial_{\xi\xi} \psi(\xi,y_0)+c\partial_{\xi}\psi(\xi,y_0)+\beta v(\xi)+ f(y,v(\xi)+\psi(\xi,y_0)).$$
Clearly, $K(0,c(y_0),y_0)(\xi)=0$ for all $\xi\in\R$, by~\eqref{traveling-wave}. By Remark~\ref{uniform-esti-psi}, we know that $\psi(\cdot,y_0) \in L^2((0,+\infty))$, $1-\psi(\cdot,y_0) \in L^2((-\infty,0))$, and $\partial_{\xi} \psi(\cdot,y_0)\in H^2(\R)$.  Moreover, since the function~$f(x,u)$ satisfies~(A1) and is globally Lipschitz continuous with respect to $u$ uniformly in~$x\in\R$, it follows that the function $\xi\mapsto f(y,v(\xi)+\psi(\xi,y_0))$ belongs to $L^2(\R)$. Therefore, for any~$(v,c,y)\in L^2(\R)\times\R\times\R$, we have $K(v,c,y) \in L^2(\R)$.

Now, for any $(v,c,y)\in H^2(\R)\times\R\times\R$, we set $G(v,c,y)=(G_1,G_2)(v,c,y)$ with
$$G_1(v,c,y)= v+M_{c,y}^{-1}(K(v,c,y)) \quad\hbox{and}\quad G_2(v,c,y)=v(0).  $$
It is easily seen from Lemma~\ref{prop-M}~(i) that the function $G$ maps $H^2(\R)\times\R\times\R$ into $H^2(\R)\times\R$. Remember that $(\psi(\cdot,y),c(y))$ is the unique solution of~\eqref{traveling-wave}-\eqref{normal-homo-wave}. It is straightforward to check that for any $y\in\R$,  we have $\psi(\cdot,y)-\psi(\cdot,y_0)\in H^2(\R)$ by~\eqref{esit-psi-l1}, and
$$G(\psi(\cdot,y)-\psi(\cdot,y_0),c(y),y)=(0,0).$$
On the other hand, thanks to Lemma~\ref{prop-M}~(ii) and the assumption that the function $f(x,u)$ is of class $C^1$ in $(x,u)\in\R^2$, and $f(x,u)$, $\partial_u f(x,u)$ are globally Lipschitz continuous in $u$ uniformly in $x\in\R$, by similar arguments to those used in the proof of \cite[Lemma 3.4]{dhz1}, one can check that the function $G: H^2(\R)\times\R\times\R \to H^2(\R)\times\R $ is continuously Fr\'echet differentiable (the verification is simpler in our case, since no singularity occurs) and that the derivatives are given by
\begin{equation*}
\partial_{(v,c,y)}G(v,c,y)(\tilde{v},\tilde{c},\tilde{y})=
\left(\baa{l}
\partial_{v}G_1(v,c,y)(\tilde{v})+\partial_{c}G_1(v,c,y)(\tilde{c})+\partial_{y}G_1(v,c,y)(\tilde{y}) \vspace{5pt}\\
\partial_{v}G_2(v,c,y)(\tilde{v})+\partial_{c}G_2(v,c,y)(\tilde{c})+\partial_{y}G_2(v,c,y)(\tilde{y})
\eaa\right)
\end{equation*}
for all $(\tilde{v},\tilde{c},\tilde{y}) \in H^2(\R)\times\R\times\R$, where
\begin{equation*}
\left\{\baa{l}
\partial_{v}G_1(v,c,y)(\tilde{v})=\tilde{v}+M_{c,y}^{-1}\big[\partial_uf(y,v+\psi(\cdot,y_0))\,\tilde{v}+\beta\tilde{v} \big], \vspace{5pt}\\
\partial_{c}G_1(v,c,y)(\tilde{c})= -\tilde{c}\,M_{c,y}^{-1}\Big\{ \partial_{\xi} \big[ M_{c,y}^{-1}(K(v,c,y))-\psi(\cdot,y_0) \big]\Big\},  \vspace{5pt}\\
\partial_{y}G_1(v,c,y)(\tilde{y})=-\tilde{y}\,M_{c,y}^{-1}\Big\{ a'(y) \partial_{\xi\xi} \big[M_{c,y}^{-1}(K(v,c,y))-\psi(\cdot,y_0)\big]-\partial_xf(y,v+\psi(\cdot,y_0))\Big\},
\eaa\right.
\end{equation*}
and
$$\partial_{v}G_2(v,c,y)(\tilde{v})=\tilde{v}(0),\quad \partial_{c}G_2(v,c,y)(\tilde{c})=0,\quad \partial_{y}G_2(v,c,y)(\tilde{y})=0.$$

We will apply the implicit function theorem for the function $G:H^2(\R)\times\R\times\R \to H^2(\R)\times\R$ and show the $C^1$-smoothness of the functions $y\mapsto c(y)$ and $y\mapsto \psi(\cdot,y)$ at $y=y_0$. To do so, we need the following lemma.

\begin{lem}\label{frechet}
The operator $Q_{y_0}= \partial_{(v,c)}G(0,c(y_0),y_0): H^2(\R)\times \R \to H^2(\R)\times\R$ is invertible. Furthermore, there exists a constant $C>0$ independent of $y_0$ such that for any $(\tilde{g},\tilde{d})\in H^2(\R)\times\R$,
\begin{equation}\label{esti-inverse}
\|Q_{y_0}^{-1}(\tilde{g},\tilde{d})\|_{H^2(\R)\times\R} \leq C \|(\tilde{g},\tilde{d}) \|_{H^2(\R)\times\R},
\end{equation}
where the space $H^2(\R)\times\R$ is endowed with norm $ \|(\tilde{g},\tilde{d}) \|_{H^2(\R)\times\R}=\|\tilde{g}\|_{H^2(\R)}+|\tilde{d}|$.
\end{lem}

\begin{proof}
Let $H$ be the linearization of~\eqref{traveling-wave} with respect to $\psi$ at $(\psi(\cdot,y),c(y))$ with $y=y_0$. Namely, we define
$$H(v)=a(y_0)v''+c(y_0)v'+\partial_u f(y_0,\psi(\cdot,y_0))v \ \ \hbox{for }v\in H^2(\R). $$
By \cite{he}, the operator $H$ and the adjoint operator $H^*$, given by $H^*(v)=a(y_0)v''-c(y_0)v'+\partial_u f(y_0,\psi(\cdot,y_0))v$ for $v\in H^2(\R)$, have algebraically simple eigenvalue $0$. Moreover, the kernel ${\rm ker}(H)$ is equal to $\R\partial_{\xi}\psi(\cdot,y_0)$ and ${\rm ker}(H^*)$ is equal to $\R w^*$, with $w^*:=w^*(\cdot,y_0)$ given by~\eqref{adjoint-ker}.  Moreover, it is also known from \cite{he} that the range of $H$ is closed in $L^2(\R)$. Based on these properties, one can conclude that the operator $Q_{y_0}: H^2(\R)\times \R \to H^2(\R)\times\R$ is invertible. The proof follows the same lines as those used in \cite[Lemma~3.4]{dhz1} and~\cite[Lemma~2.4]{he}; therefore, we do not repeat the details here.

Next, we show the estimate~\eqref{esti-inverse} (with a constant $C$ independent of $y_0$). To do so, by linearity, it is sufficient to show that
$$\sup_{y_0\in\R,\,(\tilde{g},\tilde{d})\in\mathcal{S}}\|Q_{y_0}^{-1}(\tilde{g},\tilde{d})\|_{H^2(\R)\times\R}<+\infty,$$ 
where $\mathcal{S}:=\{(\tilde{g},\tilde{d})\in H^2(\R)\times\R:\|(\tilde{g},\tilde{d})\|_{H^2(\R)\times\R}=1\}$. For any $y_0\in\R$ and $(\tilde{g},\tilde{d})\in\mathcal{S}$, set $(\tilde{v},\tilde{c})=Q_{y_0}^{-1}(\tilde{g},\tilde{d})$ and $\tilde{w}=\tilde{v}-\tilde{g}$. It then follows from the definition of $Q_{y_0}$ that
$$\tilde{w}= -M_{c(y_0),y_0}^{-1}\big[\partial_uf(y_0,\psi(\cdot,y_0))\,\tilde{v}+\beta\tilde{v}\big]-\tilde{c} M_{c(y_0),y_0}^{-1}\big[\partial_{\xi}\psi(\cdot,y_0) \big],$$
whence, owing to the definition of $H$,
$$H(\tilde{w})= -\partial_uf(y_0,\psi(\cdot,y_0))\,\tilde{g}-\beta\tilde{g}-\tilde{c} \partial_{\xi}\psi(\cdot,y_0).$$
Testing this equation with $w^*(\cdot,y_0)\in {\rm ker}(H^*)$ (which is given by~\eqref{adjoint-ker}), we obtain
$$\tilde{c} \int_{\R} \partial_{\xi}\psi(\xi,y_0)\,w^*(\xi,y_0)\,d\xi=-\int_{\R}\big[\partial_uf(y_0,\psi(\xi,y_0))\,\tilde{g}(\xi)+\beta\tilde{g}(\xi)\big]\,w^*(\xi,y_0)\,d\xi. $$
From Remark~\ref{uniform-esti-psi}, the map $y\mapsto\int_{\R} \partial_{\xi}\psi(\xi,y)w^*(\xi,y) d\xi$ is positive, periodic, continuous in~$\R$, and~$\|w^*(\cdot,y)\|_{L^2(\R)}$ is bounded uniformly in $y\in\R$. Therefore, there is a constant $C_1>0$ (independent of $y_0\in\R$ and $(\tilde{g},\tilde{d})\in\mathcal{S}$) such that
\begin{equation}\label{esti-tilde-c}
|\tilde{c}| \leq C_1 \|\tilde{g}\|_{L^2(\R)}\le C_1.
\end{equation}

It remains to estimate $\|\tilde{v}\|_{H^{2}(\R)}$. Let $\tilde{w}=\tilde{w}_1+\tilde{w}_2$, where $\tilde{w}_1$ is orthogonal to~$\partial_{\xi}\psi(\cdot,y_0)$ in~$L^2(\R)$ and $\tilde{w}_2\in \R\partial_{\xi}\psi(\cdot,y_0)$. Since $\tilde{w}_2 \in {\rm Ker}(H)$, we have
$$H(\tilde{w}_1)=H(\tilde{w})= -\partial_uf(y_0,\psi(\cdot,y_0))\,\tilde{g}-\beta\tilde{g}-\tilde{c} \partial_{\xi}\psi(\cdot,y_0).$$
We first claim that $ \|\tilde{w}_1\|_{L^2(\R)}$ is bounded uniformly with respect to $y_0\in\R$ and $(\tilde{g},\tilde{d})\in\mathcal{S}$. Assume by contradiction that this is not true. Then, thanks to~\eqref{esti-tilde-c} and the $1$-periodicity with respect to $y_0$, there exist a sequence $(y_n)_{n\in\N}\subset [0,1]$, a sequence $(\tilde{g}_n,\tilde{d}_n)_{n\in\N}\subset\mathcal{S}$, a bounded sequence $(\tilde{c}_n)_{n\in\N}\subset \R$  and  a sequence of $(p_n)_{n\in\N}\subset H^2(\R)$ with $\|p_n\|_{L^2(\R)} \to+\infty$ as $n\to+\infty$ such that~$p_n$ is orthogonal to $\partial_{\xi}\psi(\cdot,y_n)$ for each $n\in\N$ and
$$H_n(p_n):=a(y_n)p_n''+c(y_n)p_n'+\partial_u f(y_n,\psi(\cdot,y_n))p_n= -\partial_uf(y_n,\psi(\cdot,y_n))\,\tilde{g}_n-\beta\tilde{g}_n-\tilde{c}_n \partial_{\xi}\psi(\cdot,y_n).$$
For each $n\in\N$, writing
$$q_n=\frac{p_n}{\|p_n\|_{L^2(\R)}} \quad \hbox{and} \quad g_n=\frac{ -\partial_uf(y_n,\psi(\cdot,y_n))\,\tilde{g}_n-\beta\tilde{g}_n-\tilde{c}_n \partial_{\xi}\psi(\cdot,y_n)}{\|p_n\|_{L^2(\R)}},$$
we have $\|g_n\|_{L^2(\R)}\to 0$ as $n\to+\infty$ by~\eqref{esit-psi-l1}, and  $M_{c(y_n),y_n}(q_n)=-\partial_u f(y_n,\psi(\cdot,y_n))\,q_n-\beta q_n+g_n$ is uniformly bounded in $L^2(\R)$. Up to extraction of a sequence, we can assume that $y_n\to y_*\in[0,1]$ as $n\to+\infty$. It follows from Lemma~\ref{prop-M} and the boundedness of the sequence $(c(y_n))_{n\in\N}$ that the sequence~$(q_n)_{n\in\N}$ is bounded in $H^2(\R)$. Hence, up to extraction of another subsequence, the sequence~$(q_n)_{n\in\N}$ converges in $H^2(\R)$ weakly and in $C^1_{loc}(\R)$ to some $q_*\in H^2(\R)$. This implies that $a(y_*)q_*''+c(y_*)q_*'+\partial_u f(y_*,\psi(\cdot,y_*))q_* =0$. Furthermore, $q_*$ is orthogonal to~$\partial_{\xi}\psi(\cdot, y_*)$ in~$L^2(\R)$, by~\eqref{estpartpsi} and the continuity of $\partial_{\xi}\psi(\cdot, y)$ in $y$. As a consequence, $q_*=0$.

Therefore, $q_n\to0$ in $C^1_{loc}(\R)$, hence in $L^2_{loc}(\R)$, as $n\to+\infty$. In order to get a contradiction, we further show that this convergence holds in $L^2(\R)$. Indeed, by Proposition~\ref{homo-property}~(ii) and the assumption~(A2), there exists a constant $M'>0$ (independent of $n\in\N$) such that~$\partial_u f(y_n, \psi(\xi,y_n)) \leq -\gamma_0/2$ for all $\xi\in  (-\infty,-M']\cup [M',+\infty)$ and $n\in\N$. Integrating the equation $H_n(q_n)=g_n$ against $q_n$ over~$(M',+\infty)$, we get
$$-a(y_n)\int_{M'}^{+\infty} (q_n')^2(\xi)\,d\xi +  \int_{M}^{+\infty} \partial_u f(y_n,\psi(\xi,y_n)) q_n^2(\xi)\,d\xi \to 0 \,\,\hbox{ as } n\to+\infty.$$
This implies in particular that $\|q_n\|_{L^2((M',+\infty))}\to 0$ as $n\to+\infty$. The same analysis over~$(-\infty,-M')$ gives $\|q_n\|_{L^2((-\infty,-M'))}\to 0$ as $n\to+\infty$. Finally, the sequence $(q_n)_{n\in\N}$ tends to $0$ in $L^2(\R)$, which is a contradiction with the fact that $\|q_n\|_{L^2(\R)}=1$ for each $n\in\N$. We can thus conclude that $\|\tilde{w}_1\|_{L^2(\R)}$ is bounded uniformly with respect to $y_0\in\R$ and $(\tilde{g},\tilde{d})\in\mathcal{S}$.

Furthermore, since $M_{c(y_0),y_0}(\tilde{w}_1)=-\partial_uf(y_0,\psi(\cdot,y_0))\,(\tilde{w}_1+\tilde{g})-\beta(\tilde{w}_1+\tilde{g})-\tilde{c}\partial_{\xi}\psi(\cdot,y_0)$, it follows from~\eqref{esit-psi-l1},~\eqref{esti-tilde-c} and Lemma~\ref{prop-M}, together with the bondedness of $c(y_0)$ with respect to $y_0\in\R$, that there exists a positive constant $C_2$ (independent of $y_0$ and $(\tilde{g},\tilde{d})\in\mathcal{S}$) such that
\begin{equation}\label{esti-tildew-1}
\|\tilde{w}_1\|_{H^2(\R)} \leq  C_2.
\end{equation}
By the Sobolev inequality, $\max_{\xi\in\R} |\tilde{w}_1(\xi)| \leq C_3\|\tilde{w}_1\|_{H^2(\R)} \leq C_3C_2$ for some constant $C_3>0$ independent of $y_0$ and $(\tilde{g},\tilde{d})\in\mathcal{S}$.  In particular, we have $|\tilde{w}_1(0)| \leq C_3C_2$. Similarly, we have $|\tilde{g}(0)| \leq C_3 \|\tilde{g}\|_{H^2(\R)}\le C_3$. Finally, remembering that $\tilde{v}(0)=\tilde{d}\in[-1,1]$ and $\tilde{v}-\tilde{g}=\tilde{w}=\tilde{w}_1+\tilde{w}_2$, we have $\tilde{w}_2(0)=\tilde{d}-\tilde{g}(0)-\tilde{w}_1(0)$. This together with Proposition~\ref{homo-property}~(iii),~\eqref{normal-homo-wave},~\eqref{esit-psi-l1} and the fact that $\tilde{w}_2\in \R\partial_{\xi}\psi(\cdot,y_0)$ implies that
$$\|\tilde{w}_2\|_{H^2(\R)} = \left|\frac{\tilde{w}_2(0)}{\partial_{\xi}\psi(0,y_0)}\right| \|\partial_{\xi}\psi(\cdot,y_0)\|_{H^2(\R)} \leq  C_4$$
for some constant $C_4>0$ independent of $y_0$ and $(\tilde{g},\tilde{d})\in\mathcal{S}$. Combining this with~\eqref{esti-tilde-c}-\eqref{esti-tildew-1} and the fact that $\tilde{v}=\tilde{w}+\tilde{g}$,  we obtain that $\|Q_{y_0}^{-1}(\tilde{g},\tilde{d})\|_{H^2(\R)\times\R}\le C_5$ for some constant $C_5>0$ independent of $y_0$ and $(\tilde{g},\tilde{d})\in\mathcal{S}$. This completes the proof of Lemma~\ref{frechet}.
\end{proof}

Based on the above preparations, we are now ready to complete the

\begin{proof}[{\it Proof of Proposition~$\ref{homo-property}$~{\rm{(iv)}}}]
For any $y_0\in\R$, we define
$$T(y)=(\psi(\cdot,y)-\psi(\cdot,y_0), c(y))\,\,\hbox{ for } y\in\R.$$
Clearly, $T$ maps $\R$ into $H^2(\R)\times\R$, and $G(T(y),y)=(0,0)$ for all $y\in\R$, by~\eqref{traveling-wave} and the definitions of $T$ and $G$. Applying the implicit function theorem to the function $G: H^2(\R)\times\R\times\R \to H^2(\R)\times\R$, one infers that, for each $y_0\in\R$, there is $\delta>0$ such that the function $T:(y_0-\delta,y_0+\delta)\to H^2(\R)\times \R$ is continuously Fr\'echet differentiable, namely of class~$C^1$. Denote the derivative operator at $y\in(y_0-\delta,y_0+\delta)$ by $A_{y}:\R \to H^2(\R)\times \R$. Thus, the function $y\mapsto A_y$ is continuous from $(y_0-\delta,y_0+\delta)$ to $\mathcal{L}(\R,H^2(\R)\times\R)$. At $y=y_0$, we have $A_{y_0}(\tilde{y})=-Q^{-1}_{y_0}(\partial_y G(0,c(y_0),y_0)(\tilde{y}))$ for every $\tilde{y} \in \R$, that is,
$$A_{y_0}(\tilde{y})=-\tilde{y}\,Q^{-1}_{y_0}\left(M^{-1}_{c(y_0),y_0}\left(a'(y_0)\partial_{\xi\xi}\psi(\cdot,y_0)+\partial_xf(y_0,\psi(\cdot,y_0))\right),0\right).$$
Since the function $f$ is of class $C^1$ in $\R^2$ and satisfies~(A1), it follows from~\eqref{esit-psi-l1} that there exists a constant $C_1>0$ independent of $y_0$ such that $\|a'(y_0)\partial_{\xi\xi}\psi(\cdot,y_0)+\partial_xf(y_0,\psi(\cdot,y_0)) \|_{L^2(\R)} \leq C_1$. Remember that the periodic function $y_0\mapsto c(y_0)$ is bounded. Then, by Lemmas~\ref{prop-M} and~\ref{frechet}, we find a constant $C_2>0$, independent of $y_0\in\R$, such that $\|A_{y_0}(\tilde{y})\|_{H^2(\R)\times\R} \leq C_2 |\tilde{y}|$ for all $\tilde{y}\in\R$, that is,
\begin{equation}\label{esti-F-A}
\|A_{y_0}\|_{H^2(\R)\times\R} \leq C_2,
\end{equation}
where we identify $A_{y_0}\in\mathcal{L}(\R,H^2(\R)\times\R)$ to an element of $H^2(\R)\times\R$, with a slight abuse of notation.

On the other hand, notice that, for each $y_0\in\R$,
$$T(y_0+\tilde{y})-T(y_0) =(\psi(\cdot,y_0+\tilde{y})-\psi(\cdot,y_0), c(y_0+\tilde{y})-c(y_0)) = A_{y_0}\tilde{y}+\omega(y_0,\tilde{y})\ \hbox{ for all }\tilde{y}\in\R,$$
where $\omega(y_0,\tilde{y}) \in H^2(\R)\times\R$ satisfies $\|\omega(y_0,\tilde{y})\|_{H^2(\R)\times\R}=o(|\tilde{y}|)$ as $\tilde{y}\to 0$. Hence, for each $y_0\in\R$,
\be\label{der}
\left\| \left( \frac{\psi(\cdot,y_0+\tilde{y})-\psi(\cdot,y_0)}{\tilde{y}}, \frac{c(y_0+\tilde{y})-c(y_0)}{\tilde{y}}  \right)- A_{y_0}  \right\|_{H^2(\R)\times\R} \longrightarrow 0 \,\,\hbox{ as } \tilde{y}\mathop{\longrightarrow}^{\neq}0.
\ee
Writing $A_{y_0}=(A_{y_0}^1,A_{y_0}^2) \in H^2(\R)\times\R$ for each $y_0\in\R$, we see from~\eqref{der} and the previous paragraph that the functions $y_0\mapsto c(y_0)$ and $y_0\mapsto \psi(\cdot,y_0)$ are of class $C^1$ in $\R$ with derivatives $\partial_{y} \psi(\cdot,y_0)=A_{y_0}^1$ and $c'(y_0)=A_{y_0}^2$ at each $y_0\in\R$. It further follows from~\eqref{esti-F-A} that~$\|\partial_y\psi(\cdot,y_0)\|_{H^2(\R)}\leq C_2$ for all $y_0\in\R$, with $C_2$ independent of $y_0$. Using the Sobolev inequa\-lity, we immediately obtain~\eqref{esti-py-psi}. Finally, by standard elliptic estimates, the function $(\xi,y)\mapsto\psi(\xi,y)$ is of class $C^{2;1}_{\xi;y}(\R^2)$. This completes the proof of Proposition~\ref{homo-property}.
\end{proof}

\begin{rem}
Consider now a reaction-diffusion equation
\begin{equation}\label{eqadvection}
v_t=(a_L(x)v_x)_x+q_L(x)v_x+f_L(x,v),\ \ t\in\R,\ x\in\R
\end{equation}
with an advection term
$$q_L(x)=q\Big(\frac{x}{L}\Big),$$
where $q:\R\to\R$ is a given $1$-periodic function of class $C^1$. Under the same assumptions on the coefficients $a$ and $f$ as in Section~$\ref{intro}$, it follows that, for each $y\in\R$, there is a unique pair $(\Psi(\cdot,y),\gamma(y))\in C^2(\R)\times\R$ solving the $y$-frozen problem
$$\left\{\baa{l}
a(y)\partial_{\xi\xi}\Psi(\xi,y)+\gamma(y)\partial_{\xi}\Psi(\xi,y)+q(y)\partial_\xi\Psi(\xi,y)+f(y,\Psi(\xi,y))=0 \,\hbox{ for }\xi\in\R,\vspace{5pt}\\
\Psi(-\infty,y)=1,\ \Psi(+\infty,y)=0,\eaa\right.$$
where $\Psi(\xi,y)$ is unique up to shifts in $\xi$, and $\gamma(y)=c(y)-q(y)$ by uniqueness of the speed $c(y)$ for problem~\eqref{traveling-wave} $($actually, by uniqueness, $\Psi(\cdot,y)$ is a shift of the profile $\psi(\cdot,y)$ solving~\eqref{traveling-wave}$)$. Notice that the map $\gamma$ is $1$-periodic in $\R$. We conjecture that, if $\gamma(y)\neq0$ for all $y\in\R$, then problem~\eqref{eqadvection} admits non-stationary pulsating fronts $\Phi_L(x-\gamma_Lt,x/L)$ connecting $0$ and $1$ for all $L$ large enough, and that their speeds $\gamma_L$ converge as $L\to+\infty$ to the harmonic mean of the function $\gamma$ over $[0,1]$, that is,
$$\gamma_L\to\left(\int_0^1\gamma^{-1}(y)dy\right)^{-1}\ \hbox{ as }L\to+\infty.$$
However, the arguments of the paper do not extend as such to this case, especially the intersection-number arguments of Section~$\ref{sec2}$ using the comparison with stationary solutions of $y$-frozen pro\-blems, not to mention that the results of~\cite{g} used in Section~$\ref{sec3}$ did not cover the case of bistable equations with advection terms. We point out that even the existence of pulsating fronts $\Phi_L(x-\gamma_Lt,x/L)$ connecting $0$ and $1$ for all $L$ large enough is not known in this case. Furthermore, when $q$ has large enough oscillations, the speeds $\gamma(y)$ can vanish for some $y$'s and blocking phenomena associated with the existence of stationary fronts may occur, as in~\cite{dhz2,hfr,k,x2,xz}, making $\gamma_L$ vanish for large $L$. All these questions are very relevant and left open for further studies.

Coming back to~\eqref{eqL}, assume now that, instead of steady $1$, the equation~\eqref{eqL} admits a unique stable non-constant positive $L$-periodic steady state $x\mapsto p_L(x)$. After rescaling the unknown function ${u=p_Lv}$, the equation~\eqref{eqL} becomes of the type~\eqref{eqadvection} with an $L$-periodic advection term $q_L(x)=2a_L(x)p'_L(x)/p_L(x)$ and a new reaction term $g_L(x,v)=(f(x/L,p_L(x)v)-f(x/L,p_L(x))v)/p_L(x)$ which is $L$-periodic in $x$, instead of $f_L(x,v)$. Even if the constant $1$ is now a steady state of this new equation and even if one assumes that these new reaction terms $g_L(x,\cdot)$ have a bistable structure over $[0,1]$ for all $x$, we are led to the same issues as in the previous paragraph, together with the additional complication that the functions $p_L$ may not be $L$-rescalings of a common function~$p_1$, and so for the coefficients $q_L$ and $g_L$. Nevertheless, under suitable assumptions on the profiles of the functions~$u\mapsto f(x,u)$, by working directly on the equations having a stable non-constant positive steady state and adopting strategies analogous to those employed in \cite{dhz1} and the present paper, we conjecture that the existence of non-stationary pulsating fronts for large $L$~\cite[Theorem~1.5]{dhz1} and the convergence of wave speed as $L\to+\infty$ $($Theorem~$\ref{speed-limit}$ in this paper$)$ can be generalized to the case with stable non-constant positive steady state. The verifications are also left for further studies.    
\end{rem}

\section*{Acknowledgements} 
We are grateful to the anonymous referee for careful reading and helpful suggestions.

\section*{Declarations of interest} 

None.


\end{document}